\documentclass[12pt,oneside,reqno]{amsart}
\usepackage[colorlinks,linkcolor=black,citecolor=blue]{hyperref}
\usepackage{amssymb}
\usepackage[dvips]{graphics}
\usepackage{txfonts}
\usepackage{bbm}
\usepackage{cases}
\usepackage{amsmath}
\usepackage{graphicx}
\usepackage{mathrsfs}
\usepackage{stmaryrd}
\usepackage{amsfonts}
\newcommand{\comment}[1]{}
\usepackage{enumerate,amsmath,amssymb,amsthm}

\numberwithin{equation}{section}

\newcommand{\be}{\begin{eqnarray}}
\newcommand{\ee}{\end{eqnarray}}
\newcommand{\ce}{\begin{eqnarray*}}
\newcommand{\de}{\end{eqnarray*}}
\newtheorem{theorem}{Theorem}[section]
\newtheorem{lemma}[theorem]{Lemma}
\newtheorem{remark}[theorem]{Remark}
\newtheorem{definition}[theorem]{Definition}
\newtheorem{proposition}[theorem]{Proposition}
\newtheorem{Examples}[theorem]{Example}
\newtheorem{corollary}[theorem]{Corollary}

\def\eps{\varepsilon}

\def\e{\mathrm{e}}

\def\a{\alpha}

\def\p{\partial}
\def\d{\delta}

\def\[{{\Big[}}
\def\]{{\Big]}}
\def\<{{\langle}}
\def\>{{\rangle}}
\def\({{\Big(}}
\def\){{\Big)}}

\def\bx{{\mathbf{x}}}

\def\dif{{\mathord{{\rm d}}}}

\def\min{{\mathord{{\rm min}}}}

\def\no{\nonumber}
\def\={&\!\!=\!\!&}

\def\cC{{\mathcal C}}

\def\cF{{\mathcal F}}

\def\sL{{\mathcal L}}

\def\cN{{\mathcal N}}

\def\sS{{\mathcal S}}

\def\mA{{\mathbb A}}
\def\mB{{\mathbb B}}

\def\mE{{\mathbb E}}

\def\mN{{\mathbb N}}

\def\mP{{\mathbb P}}

\def\mR{{\mathbb R}}

\def\bB{{\mathbf B}}
\def\bA{{\mathbf A}}

\def\bE{{\mathbf E}}
\def\1{{\mathbf{1}}}

\def\sA{{\mathscr A}}
\def\sB{{\mathscr B}}

\def\sF{{\mathscr F}}

\def\sI{{\mathscr I}}

\def\sL{{\mathscr L}}

\def\sS{{\mathscr S}}
\def\sT{{\mathscr T}}

\def\sZ{{\mathscr Z}}\def\E{\mathbb E}

\def\geq{\geqslant}
\def\leq{\leqslant}

\def\div{\mathord{{\rm div}}}

\def\eps{\varepsilon}

\def\e{\mathrm{e}}

\def\a{\alpha}

\def\p{\partial}
\def\d{\delta}

\def\[{{\Big[}}
\def\]{{\Big]}}
\def\<{{\langle}}
\def\>{{\rangle}}
\def\({{\Big(}}
\def\){{\Big)}}

\def\bx{{\mathbf{x}}}

\def\dif{{\mathord{{\rm d}}}}

\def\min{{\mathord{{\rm min}}}}

\def\no{\nonumber}
\def\={&\!\!=\!\!&}
\def\bt{\begin{theorem}}
\def\et{\end{theorem}}
\def\bl{\begin{lemma}}
\def\el{\end{lemma}}
\def\br{\begin{remark}}
\def\er{\end{remark}}
\def\bx{\begin{Examples}}
\def\ex{\end{Examples}}
\def\bd{\begin{definition}}
\def\ed{\end{definition}}
\def\bp{\begin{proposition}}
\def\ep{\end{proposition}}
\def\bc{\begin{corollary}}
\def\ec{\end{corollary}}
\def\bpf{\begin{proof}}
\def\epf{\end{proof}}

\def\geq{\geqslant}
\def\leq{\leqslant}

\def\div{\mathord{{\rm div}}}

\def\bH{{\mathbf H}}

\def\bA{{\mathbf A}}

\def\N{\mathbb N}  
   
\def\<{\langle} \def\>{\rangle}  \def\gg{\gamma}
    
\def\d{\text{\rm{d}}}   \def\D{\scr D}
   
 \def\beq{\begin{equation}}  
 
\def\e{\text{\rm{e}}}    
  
 \def\P{\mathbb P}

\def\D{\Delta}
\allowdisplaybreaks

\begin{document}
\title{Dirichlet problem for supercritical non-local operators}
\author{Xicheng Zhang and Guohuan Zhao}
\address{Xicheng Zhang:
School of Mathematics and Statistics, Wuhan University,
Wuhan, Hubei 430072, P.R.China\\
Email: XichengZhang@gmail.com
 }
\address{Guohuan Zhao:
Applied Mathematics, Chinese Academy of Science,
Beijing, 100081, P.R.China\\
Email: zhaoguohuan@gmail.com
 }
\thanks{
Research of X. Zhang is partially supported by NNSFC grant of China (No. 11731009)  and the DFG through the CRC 1283 
``Taming uncertainty and profiting from randomness and low regularity in analysis, stochastics and their applications''. 
Research of G. Zhao is partially supported by National Postdoctoral Program for Innovative Talents (BX201600183) of China.}
\maketitle
\begin{abstract}
Let $D$ be a bounded $C^2$-domain. Consider the following Dirichlet initial-boundary problem of nonlocal operators with a drift:
$$
\partial_t u={\mathscr L}^{(\alpha)}_\kappa u+b\cdot \nabla u+f\ \mathrm{in}\ \mathbb R_+\times D,\ \ u|_{\mathbb R_+\times D^c}=0,\ u(0,\cdot)|_{D}=\varphi,
$$
where $\alpha\in(0,2)$ and $\mathscr L^{(\alpha)}_\kappa$ is an $\alpha$-stable-like nonlocal operator with kernel function $\kappa(x,z)$ 
bounded from above and below by positive constants, 
and $b:\mathbb R^d\to\mathbb R^d$ is a bounded $C^\beta$-function with $\alpha+\beta>1$, $f: \mathbb R_+\times D\to\mathbb R$ is a $C^\gamma$-function in 
$D$ uniformly in $t$ with $\gamma\in((1-\alpha)\vee 0,\beta]$, $\varphi\in C^{\alpha+\gamma}(D)$. 
Under some H\"older assumptions on $\kappa$, we show the existence of 
a unique classical solution $u\in L^\infty_{loc}(\mathbb R_+; C^{\alpha+\gamma}_{loc}(D))\times C(\mathbb R_+; C_b(D))$ to the above problem.
Moreover, we establish the following probabilistic representation for $u$
$$
u(t,x)=\mathbb E_x \Big(\varphi(X_{t}){\bf 1}_{\tau_{D}>t}\Big)+\mathbb E_x\left(\int^{t\wedge\tau_{D}}_0f(t-s,X_s){\rm d} s\right),\ t\geq 0,\ x\in D,
$$
where $((X_t)_{t\geq 0},\mathbb P_x; x\in\mathbb R^d)$ is the Markov process associated with the operator $\mathscr L^{(\alpha)}_\kappa+b\cdot \nabla$,
and $\tau_D$ is the first  exit time of $X$ from $D$. In the sub and critical case $\alpha\in[1,2)$, the kernel function $\kappa$ can be rough in $z$.
In the supercritical case $\alpha\in(0,1)$, we classify the boundary points according to the sign of
$b(z)\cdot\vec{n}(z)$, where $z\in\partial D$ and $\vec{n}(z)$ is the unit outward normal vector.
Finally, we provide an example and simulate it by Monte-Carlo method to show our results.

\bigskip
\noindent \textbf{Keywords}: 
Dirichlet problem, Nonlocal operator, Schauder's estimate, Probabilistic representation, Maximum principle \\

\noindent
 {\bf AMS 2010 Mathematics Subject Classification:}  Primary: 35R09, 60J75; Secondary: 	60G52  
\end{abstract}

\tableofcontents

 \section{Introduction and main results}
 
\noindent1.1. {\bf  Introduction.} Let $D\subset\mR^d$ be a bounded $C^2$-domain. For $\alpha\in(0,2)$,
consider the following nonlocal elliptic Dirichlet problem:
\begin{align}\label{DIR}
L^{(\alpha)}_b u:=\Delta^{\frac{\alpha}{2}}u+b\cdot\nabla u=-f\ \mbox{ in $D$ and $u=0$ in $D^c$,}
\end{align}
where $\Delta^{\frac{\alpha}{2}}:=-(-\Delta)^{\frac{\alpha}{2}}$ is the usual fractional Laplacian operator and $b:\mR^d\to\mR^d$ is a bounded
H\"older continuous vector field. 
Notice that $\Delta^{\frac{\alpha}{2}}$ 
is a nonlocal integral operator. For $x_0\in D$ and $0<R<{\rm dist}(x_0,\p D)$, define
$$
u_R(x):=R^{-\alpha}u(Rx+x_0), b_R:=b(Rx+x_0),\ f_R(x):=f(Rx+x_0).
$$
By the scaling property of $\Delta^{\frac{\alpha}{2}}$, it is easy to see that
$$
\Delta^{\frac{\alpha}{2}}u_R+R^{\alpha-1}b_R\cdot\nabla u_R=-f_R \mbox{ in $B_R(x_0):=\{x\in\mR^d: |x-x_0|<R\}$}.
$$
In particular, for $\alpha\in(0,1)$, if $R\to 0$, then the drift term will blow up. So roughly speaking, 
the first order term plays a dominant role. 
In this sense we call $L^{(\alpha)}_b$ with $\alpha\in(0,1)$ the supercritical nonlocal operator. While
for $\alpha=1$, since $\Delta^{\frac{1}{2}}$ has the same order as $b\cdot\nabla$, we shall call $L^{(1)}_b$ the critical operator;
and for $\alpha\in(1,2)$, if $R\to 0$, then the drift term will go to zero and $\Delta^{\frac{\alpha}{2}}$ plays a dominant role, 
it is naturally called subcritical operator.
From the viewpoint of analysis, in the supercritical case, it is not possible to use the standard perturbation method to handle the drift term. 
This is the main source of the difficulties of studying supercritical operators.

\medskip

Let us also explain the difficulties in studying the Dirichlet problem of supercritical nonlocal operators from the probabilistic viewpoint. 
Let $(Z_t)_{t\geq 0}$ be a rotationally invariant and symmetric $\alpha$-stable process and $b$ a Lipschitz vector field. It is well known that for each $x\in\mR^d$, 
the following SDE admits a unique strong solution $X_t(x)$,
$$
 X_t =x+\int_0^t b(X_s)\dif s+ Z_t,
$$
which determines a family of strong Markov processes $\{X,\mP_x; x\in\mR^d\}$.
Let $u\in C^2_b(D)$ be a classical solution of \eqref{DIR}. By It\^o's formula, it is easy to see that
\begin{align}\label{REP0}
u(x)=\mE_x\left(\int^{\tau_D}_0 f(X_s)\dif s\right),\ x\in D,
\end{align}
where $\mE_x$ denotes the expectation with respect to $\mP_x$ 
and $\tau_D:=\{t\geq 0: X_t\notin D\}$ is the first exit time of $X$ from $D$. 
In particular, $u(x):=\mE_x\tau_D$ satisfies $L^{(\alpha)}_b u=-1$ in $D$. 
As discussed at the beginning, in the supercritical case, the boundary behavior of $u$ is determined by the first order term $b\cdot\nabla$.
We explain this point in the case of $d=1$ and $b\equiv 1$. 
The following proposition is proven in the appendix. 

\bp
Let $D:=(0,1)$ and $\alpha\in(0,2)$. It holds that for $\alpha\in[1,2)$,
$$
{\rm (i)}\ \mP_x(X_{\tau_D}=0\mbox{ or }1)=0,\ {\rm (ii)} \ \mE_x\tau_D\leq c_\alpha d_x^{\alpha/2},\ \ x\in D,
$$
where $d_x:=(x\wedge(1-x))_+$ is the distance of $x$ to $D^c$;
and for $\alpha\in(0,1)$,
$$
{\rm (iii)}\ \inf_{x\in(0,1/4)}\mE_x\tau_D>0, \ {\rm (iv)}\ \sup_{x\in D}\mP_x(X_{\tau_D}=0)=0, \ {\rm (v)}\ \sup_{x\in D}\mP_x(X_{\tau_D}=1)>0.
$$
\ep

For $\alpha\in[1,2)$, the conclusions (i) and (ii) are well-known (see \cite{Bo, Ch-Ki-So3, Ch-So, Ro-Se1}),
which implies that $X$ always jumps out $D$ without touching the boundary and the mean time of $X$ exiting from $D$ 
goes to zero as the starting point is close to the boundary. However, when $\alpha\in(0,1)$, conclusion (iii) 
in the above proposition means that the mean time of $X$ exiting from the interval $(0,1)$ 
has a strictly positive lower bound whatever the starting point $x$ is how close to the boundary point $0$. In particular, 
$L^{(\alpha)}_1 u|_D=-1$ can not have a {\it continuous} solution in $\mR$ when $\alpha\in(0,1)$.
Conclusions (iv) and (v) means that the position of $X$ exiting from the interval $(0,1)$ never hits the boundary point $0$, 
but possibly hits the boundary point $1$.
Notice that all these phenomenons are caused by a positive direction drift. In other words,
in the supercritical case, the drift will determine the boundary behavior of the solution to $L^{(\alpha)}_b u|_D=-f$.
Thus, in order to solve the Dirichlet problem \eqref{DIR} for $\alpha\in(0,1)$,
we need to make a better understanding about the effect of drift $b$ for the boundary behavior of the solution $u$. 
The theoretical analysis tells us that $\mE_x\tau_D\leq c(\1_{\alpha\in[1,2)}d_x^{\alpha/2}+\1_{\alpha\in(0,1)}\min(1/2,(1-x)_+))$
for each $x\in (0,1)$. The following figure 
exhibits the simulation result for $\mE_x\tau_D/(\1_{\alpha\in[1,2)}d_x^{\alpha/2}+\1_{\alpha\in(0,1)}\min(1/2,(1-x)_+))$ by using R-language, which coincides with the theoretical prediction.

\begin{figure}[h!]\small
\centering
\includegraphics
[height=0.5\textwidth,width=0.8\textwidth]{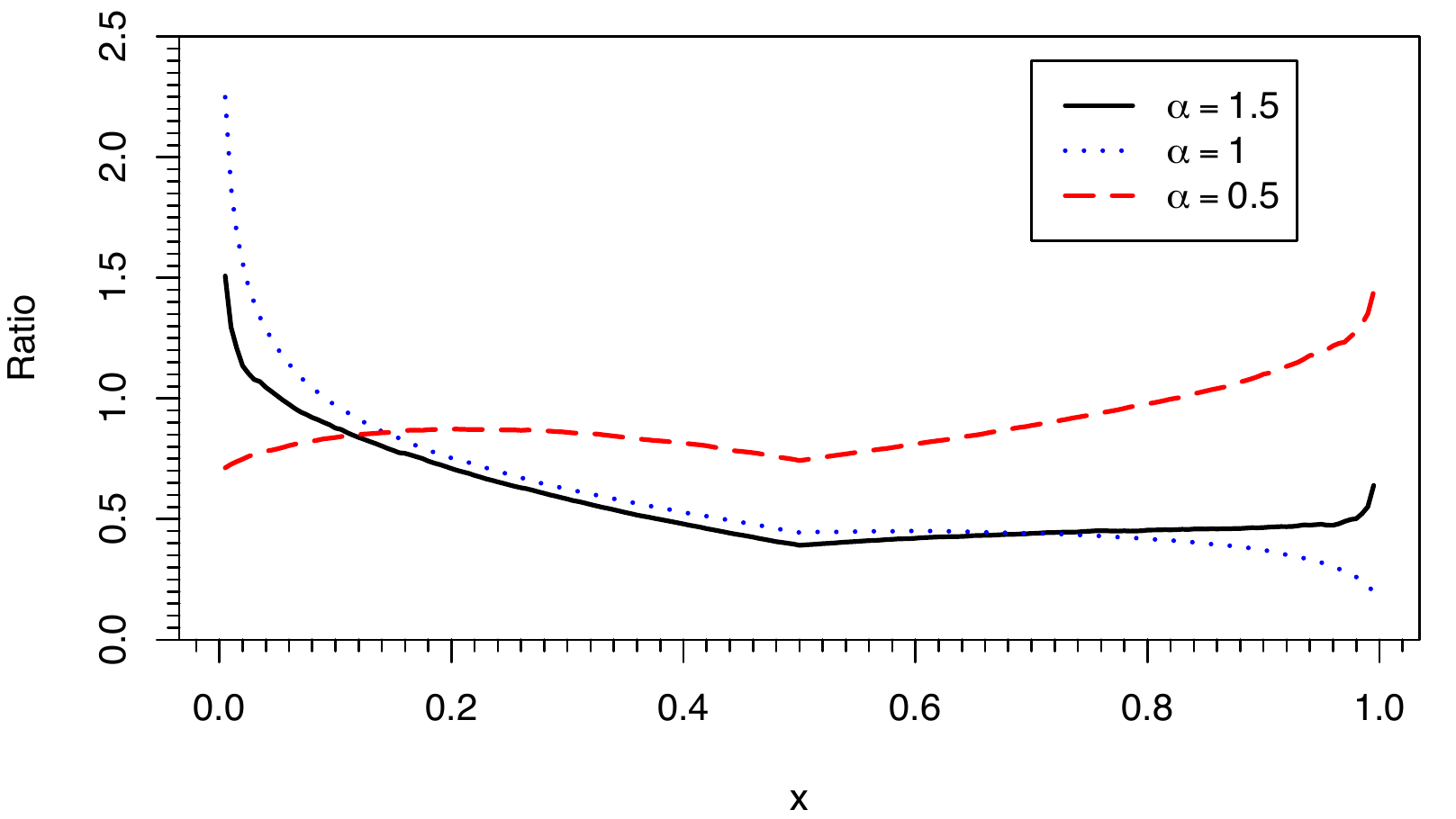}
\caption{Ratio of $\mE_x\tau_D$ to $\1_{\alpha\in[1,2)}d_x^{\alpha/2}+\1_{\{\alpha\in(0,1)\}}\min(1/2,(1-x)_+)$}\label{fige1c}
\end{figure}

 \medskip
In recent years, there is a great interest for both probabilists and analysts to study the non-local operators and related topics. 
Parts of the reasons lie in the facts that the nonlocal operators exhibit quite different features 
compared with local differential operators, and have many applications in mathematical finance, control, physics, image processing, and so on.
Up to now, there are a lot of deep works about nonlocal operators and related L\'evy processes. Let us only recall some of them related to our problem below. 
In \cite{Bo} and \cite{Ch-So}, 
the authors studied the potential theory of fractional Laplacian $\Delta^{\a/2}$, 
and the boundary Harnack principle is established therein.
Moreover, the sharp two-sided estimates of Green functions and Poisson kernels of $\Delta^{\frac{\alpha}{2}}$ in a bounded $C^{1,1}$-domain 
are also obtained in \cite{Ch-So}, see also \cite{Ch-Ki-So-Vo} for the study of boundary Harnack principle of operator $(\Delta+\Delta^{\frac{\a}{2}})|_{D}$. 
Sharp two-sided estimate of Dirichlet heat kernel of fractional Laplacian
was first proved by Chen, Kim and Song in \cite{Ch-Ki-So1}.
Later, it was extended to the operator $(\Delta^{\frac{\alpha}{2}}+b\cdot\nabla)|_D$ 
with $\alpha\in(1,2)$ in \cite{Ch-Ki-So3} and $(\Delta^{\frac{\a}{2}}+\Delta^{\frac{\beta}{2}})|_D$ in \cite{Ch-Ki-So2}.
The optimal boundary regularity of fractional Dirichlet Laplacian $\Delta^{\frac{\alpha}{2}}|_D$ was obtained by Ros-Oton and Serra in \cite{Bo}.
In the subcritical case $\a\in (1,2)$, the solvability and probabilistic representation of classical solutions
to elliptic Dirichlet problem \eqref{DIR} were studied recently by Arapostathis, Biswas and Caffarelli in \cite{Ar-Bi-Ca},
and more general nonlocal operators $\sL^{(\alpha)}_\kappa$ (see \eqref{Oper} below for a definition) are considered therein. 
In their work, besides requiring the H\"older regularity of kernel function $\kappa(x,z)$ in $x$,
some weak regularity is also imposed on the second variable $z$.
The $L^2$-estimates for nonlocal Dirichlet problems are established by energy or variational method in \cite{Fe-Ka-Vo} (see also the references therein). 
The extending problem in Sobolev spaces for nonlocal operators under minimal regularity of the exterior values is solved in \cite{Bo-Gr-Pi-Ru}. 
The global Schauder's estimates for  nonlocal operators are studied in \cite{Do-Ki} and \cite{Ba-Ka}.  
Moreover, H\"older interior estimates as well as the boundary behavior for linear and nonlinear nonlocal 
Dirichlet problems are obtained in recent works \cite{Ro-Se1,Ro-Se2, Ro-Se3, Ro-Va}, etc. 
However, none of the works  mentioned above handle the supercritical operator. 
To our best knowledge, the supercritical case was first studied by Silvestre \cite{Si2}. 
He obtained the a priori interior estimate for solutions to the following parabolic equation 
$$
\p_tu=\Delta^{\frac{\alpha}{2}}u+b\cdot\nabla u+f,
$$ 
where $\a\in (0,1)$ and $b\in L^\infty([0,T]; C^\beta)$ for some $\beta>1-\a$. 
The approach therein strongly depends on realizing the fractional Laplacian in $\mR^d$ as the boundary trace of an elliptic operator in upper half space 
of $\mR^{d+1}$. Extending this approach to the $\a$-stable-like operators  seems very hard if it is not impossible. 
We mention that similar global results are also proved in \cite{Ch-So-Zh} for more general L\'evy type operators in H\"older spaces and in \cite{Ch-Zh-Zh} for 
singular non-degenerate $\a$-stable operators in Besov spaces. 

\medskip

On the other hand, the probabilistic representation of Dirichlet problem can be dated back to the pioneering work of 
Kakutani \cite{Ka} for the harmonic functions in a domain.
A systematic probabilisitic treatment for the Dirichlet problem of Laplacian operator can be found in the monograph of Chung and Zhao \cite{Ch-Zh} (see also \cite{Fr}). For non-local operators, in the subcritical case, the probabilistic representation of nonlocal Dirichlet problem was proved in \cite{Ar-Bi-Ca}. 
It should be emphasized that probabilistic techniques have been extensively used in the studies of heat kernel estimates, H\"older estimates, Harnack inequalities
for nonlocal operators  in \cite{Ch-Zh1}, \cite{Ch-Zh14} and  \cite{So-Vo} (see also the references therein).

\medskip

Let $\alpha\in(0,2)$ and $\mR_+:=[0,\infty)$. 
In this paper we are interesting in solving he following Dirichlet problem of nonlocal parabolic equation:
\begin{align}\label{IBP0}
\left\{
\begin{aligned}
&\p_t u=\sL^{(\alpha)}_{\kappa} u+b\cdot\nabla u+f\ \mbox{ on }\ \mR_+\times D,\\
&u=0\ \mbox{ on }\ \mR_+\times D^c,\ \ u(0,\cdot)=\varphi\mbox{ on }  D,
\end{aligned}
\right.
\end{align}
where $\sL^{(\alpha)}_\kappa$ is a nonlocal operator defined by
\begin{align}\label{Oper}
\sL^{(\alpha)}_{\kappa} u(x):= \int_{\mR^d}\Xi^{(\alpha)} u(x,z)\kappa(x,z)|z|^{-d-\a}\dif z
\end{align}
with
\begin{align}\label{XI}
\Xi^{(\alpha)} u(x,z)&:=u(x+z)-u(x)-z^{(\a)}\cdot\nabla u(x),\ \ 
z^{(\alpha)}:=\1_{\{\a=1\}}z\1_{\{|z|\leq 1\}}+\1_{\{\alpha\in(1,2)\}}z,
\end{align}
and $\kappa(x,z):\mR^d\times\mR^d\to\mR_+$ satisfies that for some $\kappa_0>0$ and $\beta\in(0,1)$,
\begin{align}\label{Con1}
\left[
\begin{aligned}
&\kappa_0^{-1}\leq \kappa(x, z)\leq \kappa_0,\ |\kappa(x,z)-\kappa(x',z)|\leq\kappa_0|x-x'|^\beta\\
&\quad \1_{\a=1}\int_{r<|z|<R}z\cdot\kappa(x,z)\dif z=0,\ 0<r<R<\infty
\end{aligned}
\right].\tag{\bf H$^\beta_\kappa$}
\end{align}
The main contributions of this paper embody the following three aspects:

\begin{enumerate}[(i)]

\item In the whole space $D=\mR^d$,  under {\bf (H$^\beta_\kappa$)}
and $b\in C^\beta$ with $\alpha+\beta>1$, 
we establish the global  Schauder theory for general nonlocal parabolic equation \eqref{IBP0}
by using the Littlewood-Paley theory. To our knowledge, this is the first full result for 
nonlocal {\it parabolic} operators with {\it drifts} and 
{\em nonsymmetric rough kernels}, 
see Theorem \ref{Main2} below. 

\medskip

\item  In the sub and critical cases, we show the existence of a unique classical  solution for nonlocal Dirichlet problem \eqref{IBP0} with rough kernels.
Compared with \cite{Ar-Bi-Ca}, we do not make any regularity assumption on $\kappa(x,z)$ in the second variable $z$. It is noted that 
the interior and boundary  regularity theory for general stable L\'evy operators is established in \cite{Ro-Se2}, 
which does not cover the rough kernels as studied in this paper.

\medskip

\item In the supercritical case, by suitable boundary probabilistic estimates, we give a characterization that how the drift $b$ affects the boundary behavior of the solution. 
For a boundary point $z_0\in\p D$, let $\vec n(z_0)$ be the unit outward normal vector at point $z_0$. 
From the angle of probability, roughly to say, the sign of $b(z_0)\cdot\vec n(z_0)$ will determine whether the associated Markov process would touch the boundary 
when it exits from a bounded domain, or whether the solution would be continuous up to the boundary.
 
\end{enumerate}

\medskip 

\noindent1.2. {\bf  Statement of main results.}
We first introduce some spaces of real-valued H\"older functions in a domain.
Let $D\subset\mR^d$ be a domain.  
For an integer $k\geq 0$, denote by $C^k(D)$ the space
of all $k$-order continuous differentiable functions on $D$. For $\beta\in(0,1]$, we also denote by $C^{k,\beta}(D)$ the space of functions whose $k$-order derivatives are $\beta$-order locally H\"older continuous in $D$. For simplicity, we write for $\gamma>0$,
$$
C^\gamma(D):=C^{[\gamma],\gamma-[\gamma]}(D),\ \ C^{k,0}(D):=C^k(D),
$$
where $[\gamma]$ denotes the integer part of $\gamma$.
Let $D$ be a bounded domain. For $x,y \in D$, define 
$$
d_x:={\rm dist} (x,\p D), \ \ \ d_{x,y}:=\min\{d_x, d_y\}. 
$$
For $\theta\in\mR$, $k\in\{0\}\cup\mN$ and $0<\gamma\notin\mN$, define 
$$
[u]^{(\theta)}_{k;D}:=\sup_{x\in D} d^{k+\theta}_x|\nabla^k u(x)|,\ \ 
[u]^{(\theta)}_{\gamma;D}:=\sup_{x,y\in D} \left(d^{\gamma+\theta}_{x,y}\frac{|\nabla^{[\gamma]} u(x)-\nabla^{[\gamma]}u(y)|}{|x-y|^{\gamma-[\gamma]}}\right),
$$
where $\nabla^k$ denotes the $k$-order gradient.
For general $\gamma\geq 0$ with $\gamma+\theta\geq 0$,  we introduce the following Banach spaces for later use:
$$
\cC^{(\theta)}_{\gamma}(D):=\left\{u\in C^\gamma(D)\cap L^\infty(\mR^d): 
\|u\|^{(\theta)}_{\gamma;D}:=[u]^{(\theta)}_{0;D}+[u]^{(\theta)}_{\gamma;D}<\infty,\ u|_{D^c}=0\right\},
$$
and for  $T>0$, 
\begin{align}\label{SPAC1}
\mB^{(\theta)}_{\gamma;T}(D):=L^\infty([0,T]; \cC^{(\theta)}_\gamma(D)),\ \ 
\mB^{(\theta)}_{\gamma}(D):=\cap_{T>0}\mB^{(\theta)}_{\gamma; T}(D).
\end{align}
If the distance functions $d_x, d_{x,y}$ are not in the above definitions, we shall denote the corresponding notations by
$[\cdot]_{\gamma; D}$ and define
$$
\|u\|_{C^\gamma(D)}:=[u]_{0; D}+[u]_{\gamma; D}.
$$
In particular, if $D=\mR^d$, we shall simply write
$$
\|u\|_{C^\gamma}:=\|u\|_{C^\gamma(\mR^d)}=[u]_{0;\mR^d}+[u]_{\gamma; \mR^d}.
$$
We recall the following interpolation inequalities (see \cite{Tr}): Let $0\leq\beta<\gamma$ with $\beta\notin\mN$.
Let $D=\mR^d$ or $D$ be a bounded $C^\gamma$-domain.
For any $\eps>0$, there are constants $c_1=c_1(\beta,\gamma,D)$ and $c_2=c_2(\eps,\beta,\gamma, D)>0$ 
such that for all $u\in C^\gamma(D)$,
\begin{align}\label{Interpo2}
\|u\|_{\beta; D} \leq c_1\|u\|^{1-\frac{\beta}{\gamma}}_{0;D}\|u\|^{\frac{\beta}{\gamma}}_{\gamma;D}\leq c_2\|u\|_{0;D}+\eps\|u\|_{\gamma;D}.
\end{align}

In the following, for simplicity we write
$$
\sL^{(\alpha)}_{\kappa, b}:=\sL^{(\alpha)}_\kappa+b\cdot\nabla.
$$
\bd
We call a function $u\in L^\infty_{loc}(\mR_+; L^\infty(\mR^d))\cap C(\mR_+\times D)$ a classical solution of Dirichlet problem \eqref{IBP0} 
if it satisfies the following integral equation in the pointwise sense:
$$
u(t,x)=\varphi(x)+\int^t_0\big(\sL^{(\alpha)}_{\kappa,b} u+f\big)(s,x)\dif s\mbox{ in $\mR_+\times D$, $u|_{\mR_+\times D^c}=0$},
$$
which of course requires that $u$ is at least $C^1$-differentiable in $D$
and in the domain of $\sL^{(\alpha)}_\kappa$.
\ed
Our first aim is to show the following result.
\bt\label{MAIN}
Let $D$ be a bounded $C^2$-domain and $\alpha\in(0,2)$, $\beta\in(0,1)$.
Suppose  \eqref{Con1} and 
$$
b\in C^\beta \mbox{ if $\alpha\in[1,2)$, and $b=0$ if $\alpha\in(0,1)$}.
$$  
Then there exists $\theta_0\in(0,\frac{\alpha}{2})$ such that for any $\theta\in (0,\theta_0]$ and $\gamma\in(0,\beta]$ with $\alpha+\gamma\notin\mN$, 
if one of the following two conditions holds:
$$
{\it (i) }\ \theta\geq\gamma;\quad  {\it (ii)}\ \theta<\gamma\ \mbox{ and }\ |\kappa(x,z)-\kappa(x,z')|\leq \kappa_1|z-z'|^\gamma,
$$
then for all $f\in \mB^{(\alpha-\theta)}_{\gamma}(D)$ and $\varphi\in \cC^{(-\theta)}_{\alpha+\gamma}(D)$,
there is a unique classical solution $u\in \mB^{(-\theta)}_{\alpha+\gamma}(D)$ to equation \eqref{IBP0},
and there is a constant $c=c(\kappa_0,\kappa_1,\gamma,\theta,\alpha,\beta, d,\|b\|_{C^\beta})>0$ such that for all $T>0$,
$$
\|u\|_{\mB^{(-\theta)}_{\alpha+\gamma;T}(D)}+\|\p_tu\|_{\mB^{(\alpha-\theta)}_{\gamma;T}(D)}\leq c\|\varphi\|^{(-\theta)}_{\alpha+\gamma}
+c(1+T)\|f\|_{\mB^{(\alpha-\theta)}_{\gamma;T}(D)}.
$$
Moreover, the unique solution $u$ has the following probabilistic representation:
\begin{align}\label{EP0}
u(t,x)=\mE_x \Big(\varphi(X_{t}){\bf 1}_{\{\tau_{D}>t\}}\Big)
+\mE_x\left(\int^{t\wedge\tau_{D}}_0f(t-s,X_s)\dif s\right),
\end{align}
where $(X,\mP_x; x\in\mR^d)$ is the Markov process associated with $\sL^{(\alpha)}_{\kappa, b}$ and 
$\tau_D:=\inf\{t\geq 0: X_t\notin D\}$ is the first  exit time of $X$ from $D$. We also have the following estimate:
\begin{align}\label{ES2}
\mE_x\left(\int^{t\wedge\tau_{D}}_0|f(t-s,X_s)|\dif s\right)\leq cd^\theta_x\|f\|_{L^\infty_t(\cC^{(\alpha-\theta)}_0(D))}.
\end{align}
\et
\br\rm
Notice that in the estimate \eqref{ES2}, $f$ is allowed to be explosive near the boundary.
\er
Next we consider the supercritical case and show the following results.
\bt\label{MAIN1}
Let $\alpha, \beta\in(0,1)$ with $\alpha+\beta>1$ and $\gamma\in(1-\alpha,\beta]$. Suppose \eqref{Con1}, $b\in C^\beta$ and
$$
|\kappa(x,z)-\kappa(x,z')|\leq \kappa_1|z-z'|^\gamma.
$$
Let $D$ be a bounded $C^2$-domain, $\varphi\in \cC^{(0)}_{\alpha+\gamma}(D)$, $f\in \mB^{(0)}_\gamma(D)$.
We have the following conclusions:
\begin{enumerate}[{\bf (A)}]
\item Suppose that $b(z_0)\cdot\vec{n}(z_0)<0$  for each $z_0\in\p D$.  Equation \eqref{IBP0} admits a unique solution 
$$
u\in L^\infty_{loc}(\mR_+; C^{\alpha+\gamma}_{loc}(D)\cap L^\infty(\mR^d))\cap  C([0,\infty)\times D).
$$
\item Suppose that $b(z_0)\cdot\vec{n}(z_0)=0$ for each $z_0\in\p D$.  Equation \eqref{IBP0} admits a unique solution 
$$
u\in L^\infty_{loc}(\mR_+; C^{\alpha+\gamma}_{loc}(D)\cap L^\infty(\mR^d))\cap  C([0,\infty)\times D)\cap C((0,\infty)\times\bar D).
$$
Moreover, we also have the following boundary decay estimate: for some $\theta\in(0,1)$,
$$
|u(t,x)|\leq c \Big(\|f\|_\infty  +\|\varphi\|_\infty /t\Big)d_x^{\theta},\ \ t>0,\  x\in D.
$$
\item Suppose that $b(z_0)\cdot\vec{n}(z_0)>0$ for each $z_0\in\p D$.  Equation \eqref{IBP0} admits a unique solution 
$$
u\in L^\infty_{loc}(\mR_+; C^{\alpha+\gamma}_{loc}(D)\cap L^\infty(\mR^d))\cap  C([0,\infty)\times D)\cap C((0,\infty)\times\bar D).
$$
Moreover, we also have the following boundary decay estimate:
$$
|u(t,x)|\leq c \Big(\|f\|_\infty+\|\varphi\|_\infty/t\Big)d_x,\ \ t>0,\  x\in D.
$$
\end{enumerate}
In all cases, the unique solution $u$ still has the probabilistic representation \eqref{EP0}.
\et

We would like to make some comments about the above results. As mentioned above, in the supercritical case, the classical perturbation method
does not work. We shall use the viscosity approximation argument to show the existence. While, the uniqueness will be a consequence of  the probabilistic representation.
To reach this aim, we need to show that in case {\bf (A)}, the process does not touch the boundary when it exits from the domain $D$, and in cases {\bf (B)} and {\bf (C)}
the mean time of the process exiting from the domain $D$ has some decay estimates when the starting point approaches to the boundary.
Here a quite natural question is that whether we can consider the mixed case, that is, the general drift $b$. For this purpose, we  
define
\begin{align*}
\Gamma_>&:=\{z\in\p D: b(z)\cdot\vec{n}(z)>0\},\\
\Gamma_=&:=\{z\in\p D: b(z)\cdot\vec{n}(z)=0\},\\
\Gamma_<&:=\{z\in\p D: b(z)\cdot\vec{n}(z)<0\}.
\end{align*}
When $\alpha\in[\frac{1}{2},1)$, we have the following partial affirmative result.
\bt\label{MAIN2}
Let $\alpha\in[\frac{1}{2},1)$, $\beta\in[2(1-\alpha),1]$ and $\gamma\in(1-\alpha,\beta]$. Suppose \eqref{Con1}, $b\in C^\beta$ and
$$
|\kappa(x,z)-\kappa(x,z')|\leq \kappa_1|z-z'|^\gamma.
$$
Let $D$ be a bounded $C^2$-domain.
For any $\varphi\in \cC^{(0)}_{\alpha+\gamma}(D)$ and $f\in \mB^{(0)}_\gamma(D)$, there is a unique solution $u$ to \eqref{IBP0} in the class that
$$
u\in L^\infty_{loc}(\mR_+; C^{\alpha+\gamma}_{loc}(D)\cap L^\infty(\mR^d))\cap  C([0,\infty)\times D)\cap C((0,\infty)\times(D\cup\Gamma_=\cup\Gamma_>)),
$$
and which is given by the probabilistic representation \eqref{EP0}. Moreover, we have

\noindent (i) For each $z\in\Gamma_>$, there are $\delta, c>0$ such that
$$
\sup_{x\in D\cap B_\delta(z)}(d^{-1}_x|u(t,x)|)\leq c \Big(\|f\|_\infty+\|\varphi\|_\infty/t\Big),\ \ t>0.
$$
(ii) For each $z\in\Gamma^o_=$ (the interior of $\Gamma_=$), there are $\theta,\delta, c>0$ such that
$$
\sup_{x\in D\cap B_\delta(z)}(d^{-\theta}_x|u(t,x)|)\leq c \Big(\|f\|_\infty+\|\varphi\|_\infty/t\Big),\ \ t>0.
$$
(iii) For each $x\in D$, it holds that
$$
\mP_x(X_{\tau_D}\in\Gamma_<)=0,
$$
where $(X,\mP_x; x\in\mR^d)$ is the Markov process associated with $\sL^{(\alpha)}_{\kappa, b}$ and 
$\tau_D:=\inf\{t\geq 0: X_t\notin D\}$.
\et

Let us explain why we need to assume $\beta\geq 2(1-\alpha)$ in the above result which leads to $\alpha\geq \frac{1}{2}$.
Since $\Gamma_<$, $\Gamma_>$ and $\Gamma_=^o$ are
relatively open subsets of $\p D$, it is relatively easy to show that 
$\Gamma_<$ is inaccessible for the process $(\P_x, X_t)$ (see Lemma \ref{Le71} below),
and the points in $\Gamma_>\cup\Gamma^o_=$ are $t$-regular in the sense of \cite[page 206]{Fr} (see Lemma \ref{Le72} below).
However, for any boundary point $z_0\in\Gamma_=\setminus\Gamma_=^o$, in order to use the information $b(z_0)\cdot\vec{n}(z_0)=0$ and $b\in C^\beta$ to
show that $z_0$  is $t$-regular, we need to choose the exterior tangent ball $B$ with $B\cap \bar D=\{z_0\}$ so that $\nabla d_B(z_0)=\vec{n}(z_0)$.
For the exterior tangent ball, the fact $|x-z_0|^2\leq c\cdot{\rm dist}(x, B)$ for $x\in D$ 
leads to $\beta\geq 2(1-\alpha)$ (see \eqref{Dist} below).
Thus, dropping the condition $\beta\geq 2(1-\alpha)$ is left as an open problem.
 
\noindent1.3. {\bf  Example.}
Let $d=1$ and $D=(0,1)$. Let $Z_t$ be an one-dimensional symmetric $\alpha$-stable process with $\alpha\in(0,1)$. For each $t>0$ and $x\in D$, define
$$
X^x_t:=x\e^t+Z_t+\frac{1}{2}(\e^t-1)+\int^t_0Z_s\e^{t-s}\dif s,
$$
where $\tau^x_D:=\inf\{t\geq 0: X^x_t\notin D\}$.
Notice that $X^x_t$ solves the following SDE:
$$
\dif X^x_t=(X^x_t-\tfrac{1}{2})\dif t+\dif Z_t,\ \ X^x_0=x.
$$
Let $\varphi(x):=\1_{\{x\in D\}}\cdot \sin(3\pi x+\pi /2)$ and $u(t,x):=\bE (\varphi(X^x_{t}){\bf 1}_{\{\tau^x_{D}>t\}})$. 
Clearly, $|u(t,x)|\leq 1$ and $\varphi(x)\in\cC^{(0)}_2(D)$. 
By Theorem \ref{MAIN1}, for some $\eta>1$ and any $0<a<b<1$, we have
$$
u\in L^\infty_{loc}(\mR_+; C^{\eta}((a,b)))\cap C((0,\infty)\times [0,1]),
$$
and for any $t>0$ and $x\in D$, it holds that
$$
u(t,x)=\sin(3\pi x+\pi/2)+\int^t_0\Big(\Delta^{\frac{\alpha}{2}}u(s,x)+(x-\tfrac{1}{2})\p_x u(s,x)\Big)\dif s,
$$
and
$$
|u(t,x)|\leq c(x\wedge(1-x))_+/t, \ x\in D, t>0.
$$

\noindent1.4. {\bf  Plans and notations.}
This paper is organized as follows: In Section 2, we introduce some preliminaries about nonlocal operators. In particular, we prove 
a new Bernstein type inequality by heat kernel estimates, which allows us to establish the Schauder theory in the whole space 
in Section 3 for supercritical PDEs with rough kernels by using Littlewood-Paley theory.
In Section 4, we prove the Schauder interior estimate in weighted H\"older spaces, which is an analogue in the elliptic case as in \cite{Gi-Tr}.
In Section 5, we prove the probabilistic representation for general Dirichlet problem and also give some basic estimates of the first
exit time of the associated Markov process from a bounded domain. In Sections 6 and 7, we give the proofs of Theorems \ref{MAIN}, \ref{MAIN1} and \ref{MAIN2}, respectively.
Finally, we prove some supplementary facts   and give some numerical simulations for Example 1.3 in Appendix.
Before concluding this section, we introduce some notations used throughout this paper.
\begin{itemize}
\item $\mR_+:=[0,\infty)$ and $\mN_0:=\mN\cup\{0\}$. For a real number $a\in\mR$, we write $a_+:=\max(a,0)$.
\item For $R>0$ and $x \in\mR^d$, $B_R(x):=\{y\in\mR^d:|y-x|<R\}$ and in particular, $B_R:=B_R(0)$; 
$B_R^+:= \{x=(x_1,\cdots, x_d)\in B_R: x_1>0\}$. 
\item For $x,y\in\mR^d$, we use $x\cdot y$ or $\<x,y\>$ to denote the inner product in $\mR^d$.
\item For a set $D\subset\mR^d$, $D^c:=\mR^d\setminus D$, $D_0\Subset D$ means that ${\rm dist}(D_0, D^c)>0$,
and for $x\in\mR^d$,
$$
d_x:=d_{D^c}(x):={\rm dist}(x, D^c),\ \ \lambda_D:={\rm diam}(D)=\sup_{x,y\in D}|x-y|.
$$
\item Let $\sA$ and $\sB$ be two abstract operators acting on functions. The commutator between $\sA$ and $\sB$ is defined by
$$
[\sA,\sB] f:=\sA \sB f-\sB\sA f.
$$
\item For $T>0$ and a Banach space $\mB$, we denote $L^\infty_T(\mB):=L^\infty([0,T];\mB).$
\item 
Let $\chi:\mR^d\to[0,1]$ be a smooth function with $\chi|_{B_1}=1$  and $\chi |_{B^c_2}=0$. Define
\begin{align}\label{Cut}
\chi_R(x):=\chi(x/R),\ \ \chi^{x_0}_R(x):=\chi_R(x-x_0),\ \ R>0,\ \ x_0\in\mR^d.
\end{align}
\item The letter $c$ with or without subscripts denotes an unimportant constant.
\item We use $A\lesssim B$ to denote $A\leq c B$ for some unimportant constant $c>0$.
\end{itemize}

\section{Preliminaries}
Let $\sS$ be the Schwartz space of all rapidly decreasing functions, and $\sS'$ the dual space of $\sS$ 
called Schwartz generalized function (or tempered distribution) space. Given $f\in\sS$,
let $\cF f=\hat f$  be the Fourier transform defined by
$$
\hat f(\xi):=\int_{\mR^d}\e^{-\mathrm{i}\xi\cdot x} f(x)\dif x,\ \xi\in\mR^d.
$$
Let $\phi_0:\mR^d\to[0,1]$ be a smooth radial function with 
$$
\phi_0(\xi)=1,\ |\xi|\leq 1,\ \phi_0(\xi)=0,\ |\xi|\geq 3/2.
$$
Define
$$
\phi_1(\xi):=\phi_0(\xi)-\phi_0(2\xi).
$$
It is easy to see that $\phi_1\geq 0$, ${\rm supp}(\phi_1)\subset B_{3/2}\setminus B_{1/2}$ and
\begin{align}\label{EE1}
\phi_0(2\xi)+\sum_{j=0}^k\phi_1(2^{-j}\xi)=\phi_0(2^{-k}\xi)\stackrel{k\to\infty}{\to} 1.
\end{align}
In particular, if $|j-j'|\geq 2$, then
$$
\mathrm{supp}\phi_1(2^{-j}\cdot)\cap\mathrm{supp}\phi_1(2^{-j'}\cdot)=\emptyset.
$$
From now on we shall fix such $\phi_0$ and $\phi_1$. 
We introduce the following definitions. 
\bd
The block operator $\Delta_j$ are defined on $\sS'$ by
$$
\Delta_j f:=
\left\{
\begin{array}{ll}
\cF^{-1}(\phi_0(2\cdot) \cF f), & j=-1, \\
\cF^{-1}(\phi_1(2^{-j}\cdot) \cF f),& j\geq 0.
\end{array}
\right.
$$
For $s\in\mR$ and $p,q\in[1,\infty]$, the Besov space $B^s_{p,q}$ is defined as the set of all $f\in\sS'$ such that
$$
\|f\|_{B^s_{p,q}}:=\left(\sum_{j\geq -1}2^{jsq}\|\Delta_j f\|_p^q\right)^{1/q}<\infty,\ q\in[1,\infty)
$$
with usual modification for $q=\infty$, where $\|\cdot\|_p$ stands for the usual $L^p$-norm. 
\ed

\br\rm
It is well known that for any $0<s\notin\mN$ (cf. \cite[Theorem 2.36]{Ba-Ch-Da}),
\begin{align}\label{Ch}
c_s \|f\|_{B^s_{\infty,\infty}}\leq \|f\|_{C^s}\leq c'_s \|f\|_{B^s_{\infty,\infty}}. 
\end{align}
\er

We first recall the following Bernstein's inequality (cf. \cite[Lemma 2.1]{Ba-Ch-Da}). 
\bl\label{Bernstein}
(Bernstein's inequality) 
For any $k=0,1,2,\cdots$, there is a constant $c=c(k,d)>0$ such that for all $1\leq p\leq q\leq\infty$ and $j\geq -1$,
\begin{align}\label{B1}
\|\nabla^k\Delta_j f\|_q\leq c2^{(k+d(\frac{1}{p}-\frac{1}{q}))j}\|\Delta_jf\|_p.
\end{align}
\el

In the following we consider operator:
\begin{align}\label{LL0}
\sL^{(\alpha)}_0 u(x):= \int_{\mR^d}\Xi^{(\alpha)} u(x,z)\kappa(z)|z|^{-d-\a}\dif z,
\end{align}
where $\alpha\in(0,2)$ and $\Xi^{(\alpha)}u(x,z)$ is defined by \eqref{XI}, $\kappa(z)$ satisfies that for some $\kappa_0>0$,
\begin{align}\label{Kappa0}
\kappa_0^{-1}\leq \kappa(z)\leq\kappa_0,\ \1_{\a=1}\int_{r<|z|<R}z\cdot\kappa(z)\dif z=0,\ 0<r<R<\infty.
\end{align}
Let $(L_t)_{t\geq 0}$ be the L\'evy process with L\'evy measure $\nu(\dif z)=\kappa(z)|z|^{-d-\alpha}\dif z$.
It is well known that  under \eqref{Kappa0},
$L_t$ admits a smooth density $p_t(x)$, which enjoys the following two-sided estimates (see \cite[Theorem 2.1]{Ch-Zh11}):
for some $c_1\geq 1$,
\begin{align}\label{TS}
c^{-1}_1 \frac{t}{(t^{1/\alpha}+|x|)^{d+\alpha}}\leq p_t(x)\leq c_1 \frac{t}{(t^{1/\alpha}+|x|)^{d+\alpha}};
\end{align}
and if we define
$$
P_tf(x):=\mE f(L_t+x)=\int_{\mR^d}f(x+y)p_t(y)\dif y,
$$
then for any $f\in C^2_b(\mR^d)$,
\begin{align}\label{Se}
\p_t P_t f=\sL^{(\alpha)}_0P_t f=P_t\sL^{(\alpha)}_0 f.
\end{align}

Now we aim to prove the following Bernstein's type inequality. The crucial point is that the constant $c$ does not depend on the integrability index $p$, 
which allows us to derive
the Schauder estimate for supercritical nonlocal operators in Lemma \ref{Le12} below.
\bl\label{Le34}
Under \eqref{Kappa0}, there is a constant $c_0=c_0(\kappa_0,\alpha,d)>0$ such that for all $p\in[2,\infty)$ and  $f\in C_c(\mR^d)$,
\begin{align}\label{Ber}
\int_{\mR^d}|\Delta_jf|^{p-2}\Delta_jf\cdot\sL^{(\alpha)}_0 \Delta_jf\dif x\leq -c_02^{\alpha j}\|\Delta_jf\|_p^p,\ \ j=0,1,2,\cdots,
\end{align}
and for $j=-1$,
\begin{align}\label{Ber0}
\int_{\mR^d}|\Delta_{-1}f|^{p-2}\Delta_{-1}f\cdot\sL^{(\alpha)}_0 \Delta_{-1}f\dif x\leq 0.
\end{align}
\el
\begin{proof}
(i) Let $h=\cF^{-1} \phi_0$ be the inverse Fourier transform of $\phi_0$. Define
$$
h_{-1}(x):=\cF^{-1} \phi_0(2\cdot)(x)=2^dh(2^{-1}x),
$$
and for $j\geq 0$,
\begin{align}\label{EE7}
h_j(x):=\cF^{-1}\phi_1(2^{-j}\cdot)(x)=2^{jd}h(2^jx)-2^{(j-1)d}h(2^{j-1}x).
\end{align}
By definition it is easy to see that 
$$
\Delta_j f(x)=(h_j*f)(x)=\int_{\mR^d}h_j(x-y)f(y)\dif y,\ \ j\geq -1.
$$
By scaling, it suffices to prove \eqref{Ber} for $j=0$. Below, for simplicity we let $g =\Delta_0 f$.
\\
\\
(ii) We have the following claim: there is a constant $c_0=c_0(\kappa_0,\alpha,d)>0$ such that for all $p\in[1,\infty]$ and $t\in(0,1)$,
\begin{equation}\label{Decay}
\|P_tg\|_p\leq \e^{-c_0t}\|g\|_p.
\end{equation}
Let $\psi(\xi)$ be the L\'evy exponent of L\'evy process $L_t$. It is well known that
$\hat{p}_t(\xi)=\e^{-\psi(\xi)t}$.
Let $\varphi$ be a nonnegative smooth function with support in $B_{1/2}$ and $\varphi(0)=1.$ Define
 $$
 q^\delta_t:=\cF^{-1}\big(\e^{-\psi(\xi)t-\delta\varphi(\xi)t}\big)=p_t*\cF^{-1}\big(\e^{-\delta\varphi(\xi) t}\big).
 $$
Since by definition ${\rm supp}(\hat{g})={\rm supp}
(\phi_1\hat{f})\subseteq B_{3/2}\setminus B_{1/2}$,  we have 
$$ 
\hat{p_t}\cdot\hat{g}=\hat{q}_t^\delta\cdot \hat{g}\Rightarrow p_t*g=q^\delta_t*g.
$$
Hence,  by Young's inequality for convolutions, we get for all $p\in[1,\infty]$,
$$
\|P_tg\|_p=\|p_t*g\|_p=\|q^\delta_t*g\|_p\leq \|q_t^\delta\|_1 \|g\|_p. 
$$
If we can show that for some $\delta>0$, $q_t^\delta$ is {\em nonnegative} , then it follows that
$$
\|q_t^\delta\|_1=\hat q_t^\delta(0)=\e^{-\psi(0)t-\delta\varphi(0)t}=\e^{-\delta t},
$$
and the desired estimate \eqref{Decay} follows. 
\\
\\
(iii) To show the positivity of $q_t^\delta$ for some $\delta>0$, notice that
$$
\hat{q}^\delta_t(\xi)=\hat{p}_t(\xi)\cdot\e^{-\delta\varphi(\xi) t}=\hat{p}_t(\xi)(1+\hat r^\delta_t(\xi)),
$$
 where 
$$
r_t^\delta(x)=\cF^{-1}(\e^{-\delta\phi_1(\cdot) t}-1)(x)=(2\pi)^{-d}\int_{\mR^d} (\e^{-\delta \phi_1(\xi)t}-1)\e^{{\rm i} x\cdot\xi}\dif \xi. 
$$
Since $\varphi$ has support in $B_{1/2}$, it is easy to see that for any $m\in \mN_0$, there is a $c_1=c_1(m,\phi_1)>0$ such that for all $\delta\in(0,1)$ and $t\in(0,1)$,
$$
 \|\,|\cdot|^{2m}r^\delta_t(\cdot)\|_\infty\leq \|\D^m (\e^{-\delta \varphi (\cdot)t}-1)\|_{1}\leq c_1\delta t. 
$$
Therefore, by \eqref{TS} we get
$$
|r^\delta_t(x)|\leq c_1\delta t(1\wedge|x|^{-d-2})\leq c_2\delta p_t(x)
$$
and
$$
q^\delta_t=p_t+p_t*r^\delta_t\geq p_t-c_2\delta\cdot p_t*p_t=p_t-c_2\delta p_{2t}\geq p_t-2c_2\delta p_{t},
$$
which yields $q^\delta_t\geq p_t/2>0$ by choosing $\delta=1/(4c_2)$. So, the claim \eqref{Decay} is proven.
\\
\\
(iv) Since $g=\Delta_0 f\in\sS(\mR^d)$, there is a constant $c>0$ such that
$$
|g (x)|\leq c(1+|x|)^{-d-\alpha}.
$$ 
Hence, by \eqref{TS},  for $|x|\geq 1$ and $t\in(0,1]$,
\begin{equation}\label{Decay2}
\begin{split} 
|P_tg (x)|\leq& \int_{|y|\leq \frac{|x|}{2}} p_t(x-y)\,|g (y)|\dif y+\int_{|y|>\frac{|x|}{2}} p_t(x-y)\,|g (y)|\dif y\\
\leq &c|x|^{-d-\a} \|g \|_{1}+c|x|^{-d-\a}\int_{|y|>\frac{|x|}{2}} p_t(x-y)\dif y\leq c|x|^{-d-\a}. 
\end{split}
\end{equation}
Moreover, 
\begin{align}\label{bo}
\sup_{t\in(0,1)} \|P_tg  \|_\infty\leq \|g \|_\infty,\ \ \sup_{t\in(0,1)}\|\sL^{(\alpha)}_0 P_tg \|_\infty\leq \|\sL^{(\alpha)}_0 g \|_\infty<\infty,
\end{align}
and  for $p\in[2,\infty)$, by the chain rule and \eqref{Se}, 
$$
|P_t g(x) |^p-|g(x) |^p=\int^t_0\p_s|P_sg(x) |^p\dif s=p\int_0^t |P_sg(x) |^{p-2}P_sg(x) \cdot P_s\sL^{(\alpha)}_0 g(x)  ~~\dif s.
$$
Thus, by \eqref{Decay2} and \eqref{bo},
\begin{align*}
\sup_{t\in(0,1)}\frac{|\,|P_t g (x)|^p-|g (x)|^p|}{t}\leq c(1+ |x|)^{-d-\a}\in L^1(\mR^d).
\end{align*}
By the dominated convergence theorem, we obtain
\begin{align}\label{UY1}
\frac{\d^+ }{\dif t} \|P_t g \|^p_p\Big|_{t=0}
=\lim_{t\downarrow 0}\int_{\mR^d}\frac{|P_t g |^p-|g |^p}{t} \dif x
= p\int_{\mR^d}|g |^{p-2}g \cdot\sL^{(\alpha)}_0 g\, \dif x.
\end{align}
On the other hand, by \eqref{Decay},
$$
\frac{\dif^+}{\dif t}\|P_t g \|^p_p\Big|_{t=0}=\lim_{t\downarrow0} \frac{\|P_t g \|_p^p -\| g \|_p^p}{t}
 \leq \lim_{t\to0}\frac{\e^{-c_0p t}-1}{t}\| g \|_p^p = -c_0p\|g \|_p^p. 
$$
Combining the above two estimates and recalling $g=\Delta_0f$, we obtain \eqref{Ber} for $j=0$.
\\
\\
(v) Finally, for $j=-1$, by \eqref{UY1} with $g=\Delta_{-1} f$ and $\|P_tg \|_p\leq\|g\|_p$, we have \eqref{Ber0}.
\end{proof}
\br\rm
Estimate \eqref{Ber} with constant $c_0$ depending on $p$ was proved in \cite{Ch-Zh-Zh} by using Bernstein's inequality established in \cite{Ch-Mi-Zh}. 
One may ask whether \eqref{Ber} holds for $\alpha=2$, that is, for some $c_0=c_0(d)>0$, all $p\in[2,\infty)$ and $f\in C_c(\mR^d)$,
\begin{align}\label{ES1}
\int_{\mR^d}|\Delta_jf|^{p-2}\Delta_jf\cdot\Delta \Delta_jf\dif x\leq -c_02^{2 j}\|\Delta_jf\|_p^p,\ \ j=0,1,2,\cdots.
\end{align}
Let $g=\Delta_0f$. Notice that by the integration by parts,
\begin{align*}
\int_{\mR^d}|g|^{p-2}g\cdot\Delta g\dif x
=-(p-1)\int_{\mR^d}|g|^{p-2}|\nabla g|^2\dif x
=-\frac{4(p-1)}{p^2}\int_{\mR^d}|\nabla|g|^{p/2}|^2\dif x.
\end{align*}
Thus if \eqref{ES1} holds, then we would have
$$
\|\nabla|g|^{p/2}\|^2_2\geq \frac{c_0 p^2}{4(p-1)}\|g\|_p^p,\ \ p\in[2,\infty).
$$
However, by \cite[p.58, Lemma 2.8]{Ba-Ch-Da}, there is a $c>0$ independent of $p\geq 2$ such that
$$
\|\nabla|g|^{p/2}\|^2_2\geq c\|g\|^p_p.
$$
Therefore, we conjecture that it is not possible to find a constant $c_0>0$ independent of $p\geq 2$ so that \eqref{ES1} holds.
\er

We also need the following H\"older estimate of nonlocal operators.
\bl\label{Le26}
Let $\alpha\in(0,2)$ and $\sL^{(\alpha)}_0$ be defined by \eqref{LL0} with $\kappa(z):\mR^d\to\mR$ satisfying
$$
|\kappa(z)|\leq\kappa_0,\ \ {\bf 1}_{\alpha=1}\int_{r<|{\alpha=1}z|<R}z\cdot\kappa(z)\dif z=0,\ 0<r<R<\infty.
$$
Then there is a constant $c=c(\alpha,d)>0$ such that for all $\gamma\in\mR$ and $f\in B^{\a+\gamma}_{\infty,\infty}$,
$$
\|\sL^{(\alpha)}_0 f\|_{B^\gamma_{\infty,\infty}}\leq c\|\kappa\|_\infty\|f\|_{B^{\a+\gamma}_{\infty,\infty}}. 
$$
\el
\begin{proof}
Let $\tilde{\phi}_1$ be another smooth function supported in $B_2\setminus B_{1/4}$ with $\tilde{\phi}_1=1$ on $B_{3/2}\setminus B_{1/2}$. 
Let $\tilde h:=\cF^{-1}(\tilde\phi_1)$. Since $\tilde h\in\sS$, it is easy to see that for some $c=c(\alpha,d)>0$,
$$
\|\sL^{(\alpha)}_0 \tilde h\|_1\leq c\kappa_0<\infty.
$$
Let $\tilde{h}_j:=\cF^{-1} (\tilde{\phi_1}(2^{-j}\cdot))$ for $j=0,1,2,\cdots$. By scaling, we have 
$$
\|\sL^{(\alpha)}_0 \tilde{h}_j\|_1\leq c\kappa_0 2^{\a j},\ \ j=0,1,2,\cdots.
$$ 
Since $\widehat{\Delta_j f}=\phi_1(2^{-j}\cdot) \hat f=\tilde\phi_1(2^{-j}\cdot)\phi_1(2^{-j}\cdot)\hat f$, we have $\Delta_j f=\tilde h_j*\Delta_j f$ and
$$
\|\Delta_j \sL^{(\alpha)}_0  f\|_\infty= \|\sL^{(\alpha)}_0 (\tilde h_j* (\Delta_j f))\|_\infty
\leq\|\sL^{(\alpha)}_0 \tilde h_j\|_1\|\Delta_j f\|_\infty \leq c\kappa_02^{\a j} \|\Delta_j f\|_\infty.
$$
Similarly, one can show
$$
\|\Delta_{-1} \sL^{(\alpha)}_0  f\|_\infty\leq c\kappa_0\|\Delta_{-1}f\|_\infty.
$$
 Hence,
\begin{align*}
\|\sL^{(\alpha)}_0 f\|_{B^\gamma_{\infty,\infty}}= \sup_{j\geq -1}2^{\gamma j}\|\Delta_j \sL^{(\alpha)}_0  f\|_\infty
\leq c\kappa_0 \sup_{j\geq -1} 2^{\gamma j} 2^{\a j} \|\Delta_j f\|_\infty=c\kappa_0 \|f\|_{B^{\alpha+\gamma}_{\infty,\infty}}.
\end{align*}
The proof is complete.
\end{proof}

\section{Schauder's estimates of nonlocal parabolic equations}

In this section we establish the global Schauder estimate for the following nonlocal equation
\begin{align}\label{Cau}
\p_t u=\sL^{(\alpha)}_{\kappa} u+b\cdot\nabla u+f,\ \ u(0)=\varphi,
\end{align}
where $\alpha\in(0,2)$ and $\sL^{(\alpha)}_{\kappa}$ is defined by \eqref{Oper}.
The following commutator estimate will be used several times below.  
\bl\label{N7}
Let $\alpha\in(0,2)$ and
$\kappa(x,z):\mR^d\times\mR^d\to\mR$ be a bounded measurable function 
and satisfy that for some $\beta\in(0,1)$, $\kappa_1>0$, and all $x,x',z\in\mR^d$,
$$
|\kappa(x,z)-\kappa(x',z)|\leq\kappa_1|x-x'|^\beta.
$$
Let $\eta\in((\alpha-1)\vee 0,\alpha\wedge 1)$ and $\gamma\in(0,\beta]$. For any $R>0$, there exists 
a constant $c_R=c_R(\|\kappa\|_\infty,\kappa_1,\eta,\alpha,\gamma,d)>0$ 
such that for all $u\in C^{\eta+\gamma}$,
$$
\left\|[\chi_R,\sL^{(\alpha)}_\kappa] u\right\|_{C^\gamma}\leq c_R\|u\|_{C^{\gamma+\eta}}, 
$$
where $\chi_R$ is defined by \eqref{Cut}, and $c_R\to0$ as $R\to \infty$. 
\el
\begin{proof}
By definition \eqref{Oper}, we can write
\begin{align}\label{w-eps2}
[\chi_R,\sL^{(\alpha)}_\kappa] u(x)=:w^{(1)}_R(x)+w^{(2)}_R(x),
\end{align}
where
\begin{align*}
w^{(1)}_R(x)&:=\int_{|z|>R}\Big((\chi_R(x)-\chi_R(x+z))u(x+z)+(z^{(\alpha)}\cdot\nabla\chi_R(x)) u(x)\Big)\frac{\kappa(x, z)}{|z|^{d+\a}} \dif z,\\
w^{(2)}_R(x)&:=\int_{|z|\leq R}\Big((\chi_R(x)-\chi_R(x+z))u(x+z)+(z^{(\alpha)}\cdot\nabla\chi_R(x)) u(x)\Big)\frac{\kappa(x, z)}{|z|^{d+\a}} \dif z.
\end{align*}
For $w^{(1)}_R$, it is easy to see that
$$
|w^{(1)}_R(x)|\leq c \|\kappa\|_\infty\|u\|_\infty\left(\int_{|z|>R}(1+|z^{(\alpha)}|)|z|^{-d-\alpha}\dif z\right)\leq c_R \|u\|_\infty, 
$$
where $c_R\to0$ as $R\to \infty$. Similarly, we have
$$
[w^{(1)}_R]_\gamma\leq c_R \|u\|_{C^\gamma} \ \mbox{ with } \lim_{R\to\infty}c_R=0 . 
$$
For $w^{(2)}_R$, we treat it in two cases.
\\
\\
(Case $\alpha\in(0,1)$): Noticing that by definition,
\begin{align*}
w^{(2)}_R(x)=\int_{|z|\leq R}(\chi_R(x)-\chi_R(x+z))u(x+z)\frac{\kappa(x, z)}{|z|^{d+\a}} \dif z,
\end{align*}
we have
\begin{align*}
\|w^{(2)}_R\|_\infty\leq c\|\nabla\chi_R\|_\infty\|u\|_\infty\int_{|z|\leq  R}|z|^{1-d-\alpha}\dif z\leq c   R^{-\alpha}\|u\|_\infty
\end{align*} 
and
\begin{align*}
[w^{(2)}_R]_{\gamma}&\leq c\Big([\nabla\chi_R]_\gamma\|u\|_\infty
+\|\nabla\chi_R\|_\infty\|u\|_{C^\gamma}\Big)\int_{|z|\leq R}|z|^{1-d-\alpha}\dif z\leq c R^{-\a}\|u\|_{C^\gamma}.
\end{align*}
(Case $\alpha\in[1,2)$): By definition we can write
\begin{align*}
w^{(2)}_R(x)&=\int_{|z|\leq R}\Big((\chi_R(x)-\chi_R(x+z))u(x+z)+z\cdot(\nabla\chi u)(x)\Big)\frac{\kappa(x,z)}{|z|^{d+\alpha}}\dif z\\
&=\int_{|z|\leq  R}z\cdot\int^{\bf 1}_0\Big((\nabla\chi_R u)(x)-\nabla\chi_R(x+sz)u(x+z)\Big)\dif s\ \frac{\kappa(x,z)}{|z|^{d+\alpha}}\dif z.
\end{align*}
Fix $\eta\in(\alpha-1,1)$. We have
\begin{align*}
\|w^{(2)}_R\|_\infty\leq ([\nabla\chi_R]_\eta\|u\|_\infty+\|\nabla\chi_R\|_\infty\|u\|_{\eta})
\int_{|z|\leq R}|z|^{1+\eta-d-\alpha}\dif z\leq c R^{\eta-\a}\|u\|_{C^\eta},
\end{align*} 
and\begin{align*}
[w^{(2)}_R]_{\gamma}\leq&
([\nabla\chi_R]_{\gamma+\eta}\|u\|_\infty+[\nabla\chi_R]_{\eta}\|u\|_{C^\gamma})  \int_{|z|\leq  R}|z|^{1+\eta-d-\alpha}\dif z\\
&+([\nabla\chi_R]_{\gamma}[u]_{\eta}+\|\nabla\chi_R\|_\infty\|u\|_{C^{\gamma+\eta}})
\int_{|z|\leq  R}|z|^{1+\eta-d-\alpha}\dif z
\leq  c  R^{\eta-\a}\|u\|_{C^{\gamma+\eta}}.
\end{align*}
Combining the above two cases, we get the desired estimate.
\end{proof}

For $\gamma\in(0,1)$, $\alpha\in(0,2)$ and $T>0$, write
\begin{align}\label{SPAC11}
\bB^{\gamma}_T:=L^\infty_T(C^\gamma),\ \ 
\bA^{\alpha,\gamma}_T:=\Big\{u\in \bB^{\alpha+\gamma}_T, 
\p_t u\in \bB^{\gamma}_T\Big\}.
\end{align}
We first establish  Schauder's estimate for 
kernel $\kappa(x,z)=\kappa(z)$ by using Lemmas \ref{Le34} and \ref{N7}. 
\bl\label{Le12}
Let $\sL^{(\alpha)}_0$ be defined by \eqref{LL0} with $\kappa(z)$ satisfying \eqref{Kappa0}. Suppose that
$b\in C^\beta$ for some $\beta\in((1-\alpha)\vee 0,1)$. 
For any $\gamma\in (0,\beta]$ with $\alpha+\gamma\notin\mN$, there are constants $c=c(\kappa_0,\alpha,\beta,\gamma,d)>0$
and $m=m(\alpha,\beta,\gamma)>0$ such that for any $T>0$ and $u\in \bA^{\alpha,\gamma}_T$,
\begin{align}\label{N8}
\|u\|_{\bB^{\alpha+\gamma}_T}\leq 
c\Big(\|u(0)\|_{C^{\alpha+\gamma}}+\|f\|_{\bB^\gamma_T}
+(1+\|b\|^m_{C^\beta})\|u\|_{\bB^0_T}\Big),
\end{align}
where $f:=\p_t u-\sL^{(\alpha)}_0 u-b\cdot\nabla u$. 
\el
\br\rm
Notice that the above $b\cdot \nabla u$ is well defined in the distributional sense since $\beta+\a>1$ (see \cite{Ba-Ch-Da}). 
\er
\begin{proof}[Proof of Lemma \ref{Le12}]
(i) We first assume that $u$ has compact support. 
Using operator $\Delta_j$ act on both sides of $\p_t u=\sL^{(\alpha)}_0 u+b\cdot\nabla u+f$, we get
\begin{align*}
\p_t \Delta_ju=\sL^{(\alpha)}_0\Delta_j u+(b\cdot\nabla) \Delta_j u+[\Delta_j, b\cdot\nabla ]u+\Delta_jf.
\end{align*}
For $p>2$, multiplying both sides by $|\Delta_j u|^{p-2}\Delta_j u$ and then integrating in $x$, we obtain
\begin{align}\label{N1}
\begin{split}
\frac{\p_t\|\Delta_j u\|^p_p}{p}& =\int_{\mR^d}|\Delta_j u|^{p-2}\Delta_j u \sL^{(\alpha)}_0\Delta_j u\dif x
+\int_{\mR^d}|\Delta_j u|^{p-2} (\Delta_j u) \,(b\cdot\nabla) \Delta_j u\dif x\\
&+\int_{\mR^d}|\Delta_j u|^{p-2} (\Delta_j u) \,[\Delta_j,b\cdot\nabla]u\dif x
+\int_{\mR^d}|\Delta_j u|^{p-2} (\Delta_j u)\Delta_jf\dif x.
\end{split}
\end{align}
For the first term denoted by $\sI_j$, by Lemma \ref{Le34}, there is a constant $c_0>0$ such that
\begin{align}\label{N2}
\sI_{-1}\leq 0,\ \ \sI_j\leq-c_02^{\alpha j}\|\Delta_j u\|^p_p,\ \ j=0,1,2,\cdots,\ \ p\geq 2.
\end{align}
For the second term denoted by $\sZ_j$,  let $S_j:=\sum_{i={-1}}^{j-1}\Delta_j$ and make the following decomposition:
\begin{align*}
\sZ_j&=\int_{\mR^d}((b-S_jb)\cdot\nabla) \Delta_j u\,|\Delta_j u|^{p-2}\Delta_j u\dif x 
  +\int_{\mR^d}(S_jb\cdot\nabla) \Delta_j u\,|\Delta_j u|^{p-2}\Delta_j u\dif x =:\sZ^{(1)}_j+\sZ^{(2)}_j.
\end{align*}
For $\sZ^{(1)}_j$, by Bernstein's inequality \eqref{B1}, we have
\begin{align}\label{N3}
\begin{split}
\sZ^{(1)}_j&\leq\sum_{k\geq j}\|(\Delta_kb\cdot\nabla) \Delta_j u\|_p\|\Delta_j u\|^{p-1}_p
\leq\sum_{k\geq j}\|\Delta_kb\|_{\infty}\|\nabla\Delta_j u\|_p\|\Delta_j u\|^{p-1}_p\\
&\leq 2^j\sum_{k\geq j}\|\Delta_kb\|_{\infty}\|\Delta_j u\|^p_p
\lesssim2^{(1-\beta)j}\|b\|_{B^\beta_{\infty,\infty}}\|\Delta_j u\|^{p}_p.
\end{split}
\end{align}
Here and below, the constant contained in $\lesssim$ is independent of $p$.
For $\sZ^{(2)}_j$,   by the divergence theorem and \eqref{B1} again, we have
\begin{align}\label{N4}
\begin{split}
\sZ^{(2)}_j&=\frac{1}{p}\int_{\mR^d}(S_jb\cdot\nabla) |\Delta_j u|^p\dif x=-\frac{1}{p}\int_{\mR^d}\div S_j  b\, |\Delta_j u|^p\dif x\\
&\leq \frac{1}{p}\|\div S_j  b\|_{\infty}\|\Delta_j u\|^p_{p}\leq \frac{1}{p}\sum_{k\leq j-1}\|\div  \Delta_k b\|_{\infty}\|\Delta_j u\|^p_{p}\\
&\lesssim \sum_{k\leq j-1}2^{k}\|\Delta_kb\|_{\infty}\|\Delta_j u\|^p_{p}\lesssim 2^{(1-\beta)j}\|b\|_{B^\beta_{\infty,\infty}}\|\Delta_j u\|^p_{p}.
\end{split}
\end{align}
Combining \eqref{N1}-\eqref{N4} and by H\"older's inequality, we obtain
\begin{align*}
\p_t\|\Delta_j u\|^p_p/p&\leq-c_02^{\alpha j}{\bf 1}_{\{j\geq 0\}}\|\Delta_j u\|^p_p
+c2^{(1-\beta)j}\|b\|_{B^\beta_{\infty,\infty}}\|\Delta_j u\|^p_{p}\\
&\quad+\|[\Delta_j,b\cdot\nabla]u\|_p\|\Delta_j u\|_p^{p-1}+\|\Delta_j u\|^{p-1}_{p}\|\Delta_j f\|_p.
\end{align*}
Dividing both sides by $\|\Delta_j u\|_p^{p-1}$ and by Young's inequality 
for products (due to $\beta+\alpha>1$),
we get for some $c_1,c_2>0$ and all $j\geq -1$ and $p>2$,
\begin{align*}
\p_t\|\Delta_j u(t)\|_p&\leq-c_12^{\alpha j}\|\Delta_j u(t)\|_p+g^j_p(t),
\end{align*}
where
$$
g^j_p(t):=\left(c_2\|b\|_{C^\beta}^{\frac{\alpha}{\alpha+\beta-1}}+c_12^{-\alpha}\right)\|\Delta_j u(t)\|_p+\|[\Delta_j,b\cdot\nabla]u(t)\|_p+\|\Delta_j f(t)\|_p.
$$
By Gronwall's inequality we have
$$
\|\Delta_j u(t)\|_p\leq \e^{-c_12^{\alpha j}t}\|\Delta_j u(0)\|_p+ \int^t_0 \e^{-c_1 2^{\alpha j}(t-s)}g^j_p(s)\dif s.
$$
Letting $p$ go to infinity, we obtain
\begin{align*}
\|\Delta_j u(t)\|_\infty&\leq \e^{-c_12^{\alpha j}t}\|\Delta_j u(0)\|_\infty
+ \int^t_0 \e^{-c_1 2^{\alpha j}(t-s)}g^j_\infty(s)\dif s\\
&\leq \e^{-c_12^{\alpha j}t}\|\Delta_j u(0)\|_\infty
+ \sup_{s\in[0,t]}g^j_\infty(s)/(c_12^{\alpha j}).
\end{align*}
Hence,
\begin{align}\label{JK91}
c_1 2^{\alpha j}\|\Delta_j u(t)\|_\infty\leq c_1 2^{\alpha j}\|\Delta_j u(0)\|_\infty+\sup_{s\in[0,t]}g^j_\infty(s).
\end{align}
Noticing that by the commutator estimate proved in \cite[Lemma 2.1]{Ch-Zh-Zh},
$$
\|[\Delta_j,b\cdot\nabla]u\|_\infty\lesssim 2^{-\gamma j}
\|b\|_{C^{\beta}}\|u\|_{C^{1-\beta+\gamma}},
$$
we have
\begin{align*}
g^j_\infty(t)\lesssim \left(\|b\|_{C^\beta}^{\frac{\alpha}{\alpha+\beta-1}}+1\right)\|\Delta_j u(t)\|_\infty+ 2^{-\gamma j}\|b\|_{C^\beta}\|u\|_{C^{1-\beta+\gamma}}+
\|\Delta_j f(t)\|_\infty.
\end{align*}
Since $\|u\|_{C^s}\asymp\|u\|_{B^s_{\infty,\infty}}=\sup_{j\geq -1} \Big(2^{sj}\|\Delta_j u(t)\|_\infty\Big)$ for $0\leq s\notin\mN$,
by \eqref{JK91} and  \eqref{Interpo2}, we get
\begin{align*}
\|u(t)\|_{C^{\alpha+\gamma}}
&\leq c\left(\|u(0)\|_{C^{\alpha+\gamma}}+\|f\|_{\bB^{\gamma}_t}
+\left(\|b\|_{C^\beta}^{\frac{\alpha}{\alpha+\beta-1}}+1\right)\|u\|_{\bB^{\gamma}_t}
+\|b\|_{C^\beta}\|u\|_{\bB^{1-\beta+\gamma}_t}\right)\\
&\leq \frac{1}{2}\|u\|_{\bB^{\alpha+\gamma}_t}+
c\left(\|u(0)\|_{C^{\alpha+\gamma}}+\|f\|_{\bB^{\gamma}_t}
+\Big(1+\|b\|_{C^\beta}^m\Big)\|u\|_{\bB^0_t}\right),
\end{align*}
which yields the desired estimate.

(ii)  For general $u\in \bA^{\a,\gamma}_T$, let $u_R=u\cdot\chi_R$, where $\chi_R$ is defined in \eqref{Cut}. We have
$$
\p_t u_R =\sL^{(\a)}_0 u_R +b\cdot \nabla u_R +f\cdot\chi_R + u b\cdot \nabla \chi_R + [\chi_R, \sL^{(\a)}_0] u. 
$$
By what we have proved in step (i) and Lemma \ref{N7}, we have 
\begin{align*}
\|u_R\|_{\bB^{\alpha+\gamma}_T}\leq 
c\Big(\|u_R(0)\|_{C^{\alpha+\gamma}}+\|f_R\|_{\bB^\gamma_T}
+(1+\|b\|^m_{C^\beta})\|u_R\|_{\bB^0_T}\Big)+c R^{-1} \|ub\|_{\bB^\gamma_T}+c_R\|u\|_{\bB^{\a+\gamma}_T},
\end{align*}
where $\lim_{R\to\infty}c_R=0$. Letting $R\to \infty$ for both sides, we complete the proof. 
\end{proof}

To extend Lemma \ref{Le26} to the variable coefficient case, we need the following lemma. 
\bl\label{Var}
Let $\alpha\in(0,2)$ and
$\kappa(x,z):\mR^d\times\mR^d\to\mR$ be a bounded measurable function 
and satisfy that for some $\beta\in(0,1)$ with $\alpha+\beta\notin\mN$ and $\kappa_1>0$, all $x,x',z\in\mR^d$,
\begin{align}\label{CoN11}
|\kappa(x,z)-\kappa(x',z)|\leq\kappa_1|x-x'|^\beta,\ {\bf 1}_{\alpha=1}\int_{r<|z|<R}z\cdot\kappa(x,z)\dif z=0,\ r<R.
\end{align}
For any $\eps\in (0,1)$, there are constants $c=c(\alpha,\beta,d)>0$ and $c_\eps=c_\eps(\kappa_1,\alpha,\beta,d)>0$ such that
for all $f\in C^{\alpha+\beta}$,
$$
\|\sL^{(\alpha)}_\kappa f\|_{C^\beta}\leq (c\|\kappa\|_\infty+\eps)\|f\|_{C^{\alpha+\beta}}+c_\eps\|f\|_\infty.
$$
\el
\begin{proof}
Fix $y\in\mR^d$ and define 
$$
\sL^{(\alpha)}_yf(x):= \int_{\mR^d}\Xi^{(\alpha)} f(x,z)\kappa(y,z)|z|^{-d-\a}\dif z. 
$$
We have
$$
|\sL^{(\alpha)}_\kappa f(x)-\sL^{(\alpha)}_\kappa f(y) |\leq |\sL^{(\alpha)}_\kappa f(x)-\sL^{(\alpha)}_yf(x) |+ |\sL^{(\alpha)}_y f(x)-\sL^{(\alpha)}_y f(y)|.
$$
For the first term denoted by $\sI$, by definition  and  \eqref{Interpo2}, we have for $\theta\in(0,\beta)$,
\begin{align*}
\sI\leq& \left|\int_{\mR^d}\Xi^{(\alpha)} f(x,z)(\kappa(x,z)-\kappa(y,z))|z|^{-d-\a}\dif z\right|\\
\leq &\kappa_1|x-y|^\beta \int_{\mR^d} |\Xi^{(\alpha)} f(x,z)|\,|z|^{-d-\a}\dif z\\
\leq & c_\theta|x-y|^\beta\|f\|_{C^{\a+\theta}}\leq |x-y|^\beta\left(\eps\|f\|_{C^{\alpha+\beta}}+c_\eps\|f\|_\infty\right),
\end{align*} 
where we have used that $|\Xi^{(\alpha)} f(x,z)|\leq \|f\|_{C^{\a+\theta}}(|z|^{\alpha+\theta}\wedge 1)$.
For the second term, by Lemma \ref{Le26} and \eqref{Ch}, we have
$$
|\sL^{(\alpha)}_y f(x)-\sL^{(\alpha)}_y f(y) |\leq |x-y|^\beta \|\sL^{(\alpha)}_yf\|_{C^\beta}\leq c\|\kappa\|_\infty\|f\|_{C^{\a+\beta}} |x-y|^\beta. 
$$
Combining the above inequalities, we complete the proof. 
\end{proof}

Now we can show the following variable coefficients estimate.
\bt\label{Main2}
Let $\beta\in(0,1)$ with $\alpha+\beta>1$ and $\gamma\in (0,\beta]$ with $\alpha+\gamma\notin\mN$. 
Under \eqref{Con1} and $b\in C^\beta$,
there are constants $c=c(\kappa_0,\alpha,\beta,\gamma,d)>0$ and $m=m(\alpha,\beta,\gamma)>0$ such that
for all  $T>0$ and $u\in \bA^{\alpha,\gamma}_T$,
\begin{align}\label{N88}
\|u\|_{\bB^{\alpha+\gamma}_T}\leq 
c\Big(\|u(0)\|_{C^{\alpha+\gamma}}+(1+\|b\|^m_{C^\beta})\|u\|_{\bB^0_T}+\|f\|_{\bB^{\gamma}_T}\Big),
\end{align}
where $f:=\p_t u-\sL^{(\alpha)}_\kappa u-b\cdot\nabla u$.
\et
\begin{proof}
We use the freezing coefficients argument to show \eqref{N88}.
Fix $x_0\in\mR^d$ and $\eps\in(0,1)$. Let $\chi^{x_0}_\eps$ be defined as in \eqref{Cut} with $R=\eps$ and define
$$
\kappa_\eps(x,z):=[\kappa(x,z)-\kappa(x_0,z)]\chi^{x_0}_\eps(x),\ \
\sL^{(\alpha)}_0 u(x):= \int_{\mR^d}\Xi^{(\alpha)} u(x,z)\kappa(x_0,z)|z|^{-d-\a}\dif z.
$$
Let $u\in \bA^{\alpha,\gamma}_T$. Define 
$$
f:=\p_t u-\sL^{(\alpha)}_\kappa u-b\cdot\nabla u,\ \ u_\eps:=u\chi^{x_0}_\eps.
$$
It is easy to see that
\be\label{loc}
\begin{split}
\p_t u_\eps=&\sL^{(\alpha)}_0 u_\eps+b\cdot \nabla u_\eps+\chi^{x_0}_\eps f-ub\cdot \nabla \chi^{x_0}_\eps
+\chi^{x_0}_\eps \big(\sL^{(\alpha)}_\kappa u-\sL^{(\alpha)}_0 u\big)+[\chi^{x_0}_\eps,\sL^{(\alpha)}_0]u.
\end{split}
\ee
Since $\gamma\in(0,\beta]$, we obviously have
\be\label{loc-1}
\|\chi^{x_0}_\eps f-ub\cdot \nabla \chi^{x_0}_\eps\|_{C^\gg} \leq c_\eps (\|f\|_{C^\gamma}+\|u\|_{C^\gamma}).
\ee
Noticing that
$$
\chi^{x_0}_\eps(x) (\sL^{(\alpha)}_\kappa u-\sL^{(\alpha)}_0 u)(x)=\int_{\mR^d}\Xi^{(\alpha)} u(x,z) \kappa_\eps(x,z)|z|^{-d-\a}\dif z=:\sL^{(\a)}_{\kappa_\eps}u(x),
$$
by \eqref{Con1} and Lemma \ref{Var}, 
we have
\begin{align}\label{N6}
\|\sL^{(\a)}_{\kappa_\eps}u\|_{C^\gamma}\leq (c\eps^\beta+\eps)  
\|u\|_{C^{\a+\gamma}}+ c_{\eps}\|u\|_{\infty}. 
\end{align}
Moreover, by Lemma \ref{N7} and \eqref{Interpo2}, we also have
\begin{align}\label{N77}
\|[\chi^{x_0}_\eps,\sL^{(\alpha)}_0]u\|_{C^\gamma}
\leq \eps  \|u\|_{C^{\a+\gamma}}+ c_{\eps}\|u\|_{\infty}. 
\end{align}
By  \eqref{N8}, \eqref{loc}, \eqref{loc-1}, \eqref{N6} and \eqref{N77}, choosing $\eps$ small enough we get
\begin{align*}
\|u\|_{L^\infty_T(C^{\alpha+\gamma} (B_{\eps/2}(x_0)))}&\leq
\|u_\eps\|_{\bB^{\alpha+\gamma}_T}
\leq\tfrac{1}{2}\|u\|_{\bB^{\a+\gamma}_T}
+ c_\eps\Big(\|u(0)\|_{C^{\alpha+\gamma}}+(1+\|b\|^m_{C^\beta})\|u\|_{\bB^0_T}+\|f\|_{\bB^\gamma_T}\Big),
\end{align*}
where $c_\eps$ is independent of $x_0$. Thus we obtain \eqref{N88} by taking supremum in $x_0\in\mR^d$.
\end{proof}
\br\label{RE35}
\rm
\begin{enumerate}[(i)]
\item When $\alpha\in[1,2)$, 
by the a priori estimate \eqref{N88}, it is by now standard to 
show the existence of a solution $u\in \bA^{\alpha,\gamma}_T$ to Cauchy problem  \eqref{Cau} 
by continuity method (see the proof of Theorem \ref{MAIN} below).
\item When $\alpha\in(0,1)$, consider  the following viscosity approximation equation:
\begin{align}\label{VIS}
\p_t u=\nu\Delta^{1/2}u+\sL^{(\alpha)}_\kappa u+b\cdot\nabla u+f,\ \ \nu>0.
\end{align}
Under the assumptions of Theorem \ref{Main2}, 
from the proof of Theorem \ref{Main2}, it is easy to see that the following uniform estimate holds
\begin{align}\label{N888}
\|u\|_{\bB^{\alpha+\gamma}_T}\leq 
c\Big(\|\varphi\|_{C^{\alpha+\gamma}}+\|f\|_{\bB^\gamma_T}+(1+\|b\|^m_{C^\beta})\|u\|_{\bB^0_T}\Big),
\end{align}
where the constants $c,m>0$ are independent of $T,\nu>0$.
From this viscosity approximation and uniform estimate \eqref{N888}, we can also show 
the existence of a solution $u\in \bB^{\alpha+\gamma}_T$ to supercritical equation \eqref{Cau} by a standard compact argument (see the proof of Theorem\ref{MAIN1} below).
\end{enumerate}
\er

\section{Schauder's interior estimates}

The following simple lemma is quite useful, which provides a way of treating the weighted H\"older norm by the usual H\"older's norm.
\bl\label{Le21}
Let $D$ be a bounded domain and  $\mu\in(0,1)$. For $x_0\in D$ and $u\in C(D)$, define
$$
R:=\mu\,{\rm dist}(x_0,\p D),\ \ \ u^{x_0}_R(x):=u(Rx+x_0).
$$
For any $\theta\in\mR$, $\gamma\geq 0$ and $m\in\mN$ with $m\mu<1$, there is a constant $c=c(\theta,\gamma,\mu,m,d)\geq 1$ such that
for all $u\in\cC^{(\theta)}_\gamma(D)$,
\begin{align}\label{HL6}
c^{-1}\| u\|^{(\theta)}_{\gamma; D}\leq\sup_{x_0\in D}\left(R^{\theta}\|u^{x_0}_R\|_{\gamma; B_m}\right)\leq c\| u\|^{(\theta)}_{\gamma; D}.
\end{align}
\el
\begin{proof}
Write $\cN:=\sup_{x_0\in D}R^{\theta}\|u^{x_0}_R\|_{\gamma; B_m}$. For any $\delta\geq 0$, noticing that 
$$
[u^{x_0}_R]_{\delta;B_m}=R^\delta[u]_{\delta;B_{mR}(x_0)},
$$
we have
\begin{align}\label{GF2}
R^{\theta}[u^{x_0}_R]_{\delta; B_m}\leq (1-m\mu)^{-|\theta+\delta|} [u]^{(\theta)}_{\delta; D},
\end{align}
which in turn implies that $\cN\leq c\| u\|^{(\theta)}_{\gamma; D}$ for some $c>1$.
Next, for $k=0,\cdots,[\gamma]$, we have
$$
d^{\theta+k}_{x_0}|\nabla^k u(x_0)|=R^{\theta}/\mu^{\theta+k}|\nabla^k u^{x_0}_R(0)|\leq \cN/\mu^{\theta+k}.
$$
For $x_0,y_0\in D$ with $d_{x_0}\leq d_{y_0}$, if $|x_0-y_0|\leq R$, then 
\begin{align*}
d^{\theta+\gamma}_{x_0,y_0}\frac{|\nabla^{[\gamma]} u(x_0)-\nabla^{[\gamma]} u(y_0)|}{|x_0-y_0|^{\gamma-[\gamma]}}
\leq R^{\theta}/\mu^{\theta+\gamma} [u^{x_0}_R]_{\gamma; B_m}\leq \cN/\mu^{\theta+\gamma}. 
\end{align*}
If $|x_0-y_0|>R$, then 
\begin{align*}
d^{\theta+\gamma}_{x_0,y_0}\frac{|\nabla^{[\gamma]} u(x_0)-\nabla^{[\gamma]} u(y_0)|}{|x_0-y_0|^{\gamma-[\gamma]}}
\leq R^{\theta+[\gamma]}/\mu^{\theta+\gamma} \Big(|\nabla^{[\gamma]} u(x_0)|+|\nabla^{[\gamma]} u(y_0)|\Big)\leq \cN/\mu^{\theta+\gamma}.
\end{align*}
Thus we obtain $\| u\|^{(\theta)}_{\gamma; D}\leq c\cN$ by taking supremum with respect to $x_0,y_0\in D$.
\end{proof}

As a corollary we have the following interpolation result.
\bl\label{Le42}
Let $0\leq\gamma_0<\gamma_1<\gamma_2$ with $\gamma_1\notin\mN$ and 
$r:=(\gamma_2-\gamma_1)/(\gamma_2-\gamma_0)$. 
For any $T>0$ and $\theta\in\mR$, there is a constant $c>0$ such that for all $u\in \mB^{(\theta-\gamma_2)}_{\gamma_2; T}(D)$ with
$\p_t u\in \mB^{(\theta-\gamma_0)}_{\gamma_0; T}(D)$,
\begin{align}\label{EX1}
\|u(t_1)-u(t_0)\|_{\cC^{(\theta-\gamma_1)}_{\gamma_1}(D)}\leq c (t_1-t_0)^r 
\|\p_t u\|^r_{\mB^{(\theta-\gamma_0)}_{\gamma_0; T}(D)}\|u\|^{1-r}_{\mB^{(\theta-\gamma_2)}_{\gamma_2;T}(D)}.
\end{align}
In particular, if $(u_n)_{n\in\mN}$ is a bounded sequence in $\mB^{(\theta-\gamma_2)}_{\gamma_2; T}(D)$ with
$(\p_t u_n)_{n\in\mN}$ bounded in $\mB^{(\theta-\gamma_0)}_{\gamma_0; T}(D)$, then there are $u\in\mB^{(\theta-\gamma_2)}_{\gamma_2; T}(D)$
and a subsequence $u_{n_k}$ such that for any $\eps\in(0,\gamma_2-\gamma_0)$ and $D_0\Subset D$,
\begin{align}\label{EX2}
\lim_{k\to\infty}\|u_{n_k}-u\|_{L^\infty_T(C^{\gamma_2-\eps}(\bar D_0))}=0.
\end{align}
\el
\begin{proof}
First of all, by \eqref{HL6} and the usual interpolation inequality \eqref{Interpo2}, we have 
\begin{align}\label{GF44}
\|u\|^{(\theta-\gamma_1)}_{\gamma_1; D}\leq c\left(\|u\|^{(\theta-\gamma_0)}_{\gamma_0;D}\right)^r\left(\|u\|^{(\theta-\gamma_2)}_{\gamma_2;D}\right)^{1-r}.
\end{align}
For any $\beta\in(0,r)$ and $q>1/\beta$, by Garsia-Rademich-Rumsey's inequality (see \cite{St-Va}), we have
\begin{align*}
&\|u(t_1)-u(t_0)\|^q_{\cC^{(\theta-\gamma_1)}_{\gamma_1}(D)}\lesssim(t_1-t_0)^{\beta q-1}\int^{t_1}_{t_0}\!\!\int^{t}_{t_0}
\|u(t)-u(s)\|^q_{\cC^{(\theta-\gamma_1)}_{\gamma_1}(D)}(t-s)^{-1-\beta q}\dif s\dif t\\
&\quad\lesssim(t_1-t_0)^{\beta q-1}\int^{t_1}_{t_0}\!\!\int^{t}_{t_0}
\|u(t)-u(s)\|^{r q}_{\cC^{(\theta-\gamma_0)}_{\gamma_0}(D)}
\|u(t)-u(s)\|^{(1-r)q}_{\cC^{(\theta-\gamma_2)}_{\gamma_2}}(t-s)^{-1-\beta q}\dif s\dif t\\
&\quad\lesssim(t_1-t_0)^{\beta q-1}\left(\int^{t_1}_{t_0}\!\!\int^{t}_{t_0}\frac{|t-s|^{r q}}{(t-s)^{1+\beta q}}
\dif s\dif t\right)\|\p_t u\|^{rq}_{\mB^{(\theta-\gamma_0)}_{\gamma_0;T}(D)}
\|u\|^{(1-r)q}_{\mB^{(\theta-\gamma_2)}_{\gamma_2;T}(D)},
\end{align*}
which gives \eqref{EX1}. As for \eqref{EX2}, it is a  direct consequence of  \eqref{EX1} and Ascolli-Arzela's lemma.
\end{proof}

We prepare the following crucial lemma for later use.
\bl
Let $D$ be a bounded domain and $\alpha\in(0,2)$. Let
$\kappa(x,z):\mR^d\times\mR^d\to\mR$ be bounded by $\kappa_0$ and satisfy that for some $\beta\in(0,1)$ and $\kappa_1>0$,
$$
|\kappa(x,z)-\kappa(x',z)|\leq\kappa_1|x-x'|^\beta.
$$
Let $\gamma\in[0,\beta]$ and $\theta\in[0,\alpha\wedge 1)$. Suppose that one of the following two conditions holds:
$$
{\it (i) }\ \theta\geq\gamma;\quad  {\it (ii)}\ \theta<\gamma\ \mbox{ and }\ |\kappa(x,z)-\kappa(x,z')|\leq\kappa_2|z-z'|^\gamma.
$$
Then there is a constant $c=c(\alpha,\beta,\gamma,\theta,d,\kappa_0,\kappa_1,\kappa_2,\lambda_D)>0$ such that for all  $x_0\in D$,
\begin{align}\label{GF4}
\Big\|\Big(\sL^{(\alpha)}_\kappa\big((1-\chi^{x_0}_{2R}) u\big)\Big)^{x_0}_R\Big\|_{\gamma; B_1}
\leq cR^{\theta-\alpha}[u]_{\theta; D},\ \ u\in \cC^{(-\theta)}_\theta(D),
\end{align}
where $R:={\rm dist}(x_0,\p D)/8$ and $f^{x_0}_R(x):=f(Rx+x_0)$ for a function $f$.
\el
\begin{proof}
(i) Assume $\theta\geq\gamma$. Let $\bar u:=(1-\chi^{x_0}_{2R}) u$. Noticing that for any $x\in B_R(x_0)$ 
 and $z\in \bar B_R$, $\bar u(x+z)=0$,  
by definition \eqref{XI}, we have for $x\in B_R(x_0)$,
\begin{align}\label{GP1}
\sL^{(\alpha)}_\kappa\bar u(x)=\int_{\mR^d}\Xi^{(\alpha)}\bar u(x,z)\frac{\kappa(x,z)}{|z|^{d+\alpha}}\dif z
=\int_{|z|>R}\Xi^{(\alpha)}\bar u(x,z)\frac{\kappa(x,z)}{|z|^{d+\alpha}}\dif z.
\end{align}
Let $x^*_0\in\p D$ be such that
$$
{\rm dist}(x_0,x^*_0)={\rm dist}(x_0,\p D)=8R.
$$
Fix $u\in \cC^{(-\theta)}_\theta(D)$. Since $\cC^{(-\theta)}_{\theta}(D)\subset C^\theta(\mR^d)$, we have for any $x\in B_R(x_0)$ and $z\in\mR^d$,
\begin{align}\label{GP3}
|u(x+z)-u(x^*_0)|\leq (|x-x^*_0|+|z|)^\theta[u]_{\theta; \mR^d}=(9R+|z|)^\theta[u]_{\theta; D}.
\end{align}
For  $x\in B_R(x_0)$, since $\chi^{x_0}_{2R}(x)=1$ and $u(x^*_0)=0$, by definition \eqref{XI} we have
\begin{align}\label{GP2}
\Xi^{(\alpha)}\bar u(x,z)=(1-\chi^{x_0}_{2R}(x+z)) (u(x+z)-u(x^*_0)),\ \ z\in\mR^d.
\end{align}
Thus, by \eqref{GP3}, for any $x\in B_R(x_0)$ and $z\in\mR^d$, we have
$$
|\Xi^{(\alpha)}\bar u(x,z)|\leq|u(x+z)-u(x^*_0)|\leq(9R+|z|)^\theta[u]_{\theta; D},
$$
which yields by \eqref{GP1},
\begin{align}\label{PY1}
\big[\big(\sL^{(\alpha)}_\kappa \bar u\big)^{x_0}_R\big]_{0; B_1}
=[\sL^{(\alpha)}_\kappa\bar u]_{0; B_R(x_0)}
\lesssim [u]_{\theta; D}\int_{|z|>R} \frac{(R+|z|)^{\theta}}{|z|^{d+\alpha}}\dif z\lesssim R^{\theta-\alpha}[u]_{\theta; D}.
\end{align}
On the other hand, for any $x,x'\in B_R(x_0)$, since $\theta\geq\gamma$, by \eqref{GP2} and \eqref{GP3}, we have
\begin{align*}
&|\Xi^{(\alpha)}\bar u(x,z)-\Xi^{(\alpha)}\bar u(x',z)|\leq
|\chi^{x_0}_{2R}(x+z)-\chi^{x_0}_{2R}(x'+z)|\, |u(x+z)-u(x^*_0)|\\
&\quad+|u(x+z)-u(x'+z)|^{\frac{\gamma}{\theta}}(|u(x+z)-u(x^*_0)|+|u(x'+z)-u(x^*_0)|)^{1-\frac{\gamma}{\theta}}\\
&\quad\lesssim R^{-1} |x-x'|  (R+|z|)^\theta[u]_{\theta; D}+|x-x'|^\gamma(R+|z|)^{\theta-\gamma}[u]_{\theta; D},
\end{align*}
which yields by \eqref{GP1} and H\"older continuity of $\kappa(x,z)$ in $x$,
\begin{align}\label{PY2}
\big[\big(\sL^{(\alpha)}_\kappa \bar u\big)^{x_0}_R\big]_{\gamma; B_1}
=R^\gamma[\sL^{(\alpha)}_\kappa\bar u]_{\gamma; B_R(x_0)}
\lesssim [u]_{\theta; D}\int_{|z|>R} \frac{(R+|z|)^{\theta}}{|z|^{d+\alpha}}\dif z\lesssim R^{\theta-\alpha}[u]_{\theta; D}.
\end{align}
(ii)  Since $\bar u(z)=0$ for $z\in B_{2R}(x_0)$,  we can write for $x\in B_R(x_0)$,
$$
\sL^{(\alpha)}_\kappa\bar u(x)=\int_{\mR^d}\bar u(z)\frac{\kappa(x,x-z)}{|x-z|^{d+\alpha}}\dif z
=\int_{B^c_{2R}(x_0)}\bar u(z)\frac{\kappa(x,x-z)}{|x-z|^{d+\alpha}}\dif z.
$$
Since $|\kappa(x,z)-\kappa(x',z')|\lesssim(|x-x'|+|z-z'|)^\gamma$, we have for $x,y\in B_R(x_0)$, 
\begin{align*}
|\sL^{(\alpha)}_\kappa\bar u(x)-\sL^{(\alpha)}_\kappa\bar u(y)|
&=\left|\int_{B^c_{2R}(x_0)}\bar u(z) \left(\frac{\kappa(x, z-x)}{|z-x|^{d+\a}}-\frac{\kappa(y, z-y)}{|z-y|^{d+\a}}\right)\dif z\right|\\
&\lesssim  [u]_{0;D} \int_{B^c_{2R}(x_0)} \left|\frac{\kappa(x, z-x)}{|z-x|^{d+\a}}-\frac{\kappa(y, z-y)}{|z-y|^{d+\a}}\right| \dif z
\lesssim [u]_{0;D}|x-y|^{\gamma} R^{-\alpha},
\end{align*}
which implies that for $\theta<\gamma$,
\begin{align}\label{PY4}
\big[\big(\sL^{(\alpha)}_\kappa \bar u\big)^{x_0}_R\big]_{\gamma; B_1}
=R^\gamma[\sL^{(\alpha)}_\kappa\bar u]_{\gamma; B_R(x_0)}\lesssim R^{\gamma-\alpha}[u]_{0;D}
\leq \lambda_D^{\gamma-\theta} R^{\theta-\alpha}[u]_{0;D}.
\end{align}
Combining \eqref{PY1}, \eqref{PY2} with \eqref{PY4}, we obtain \eqref{GF4}.
\end{proof}
Now we can show the following local result of nonlocal operators in weighted H\"older spaces.
\bt\label{Th0}
Let $D$ be a bounded domain and $\alpha\in(0,2)$. Let
$\kappa(x,z):\mR^d\times\mR^d\to\mR$ be bounded by $\kappa_0$ and satisfy that for some $\beta\in(0,1)$ and $\kappa_1>0$,
\begin{align}\label{CoN1}
|\kappa(x,z)-\kappa(x',z)|\leq\kappa_1|x-x'|^\beta,\ \1_{\a=1}\int_{r<|z|<R}z\cdot\kappa(x,z)\dif z=0,\ r<R.
\end{align}
Let $\gamma\in(0,\beta]$ with $\alpha+\gamma\notin\mN$ and $\theta\in(-\infty,\alpha\wedge 1)$. 
Suppose that one of the following two conditions holds:
$$
{\it (i) }\ \theta\geq\gamma;\quad  {\it (ii)}\ \theta<\gamma\ \mbox{ and }\ |\kappa(x,z)-\kappa(x,z')|\leq\kappa_2|z-z'|^\gamma.
$$
Then there is a constant $c=c(\alpha,\beta,\gamma,\theta,d,\kappa_0,\kappa_1,\kappa_2,\lambda_D)>0$ such that for all $u\in \cC^{(-\theta)}_{\alpha+\gamma}(D)$,
\begin{align}\label{HL5}
\|\sL^{(\alpha)}_\kappa u\|^{(\alpha-\theta)}_{\gamma;D}\leq c\left(\|u\|^{(-\theta)}_{\alpha+\gamma; D}+{\bf 1}_{\{\theta<0\}}\|u\|_{0; D}\right).
\end{align}
\et
\begin{proof}
For $u\in\cC^{(-\theta)}_{\alpha+\gamma}(D)$ and $x_0\in D$, let $R:={\rm dist}(x_0,\p D)/8$ and $\chi^{x_0}_{2R}$ be defined by \eqref{Cut}.
Define
$$
\tilde u:=\chi^{x_0}_{2R} u, \quad \bar u:=(1-\chi^{x_0}_{2R}) u,\ \ \kappa_R(x,z):=\kappa(Rx+x_0,Rz).
$$
By scaling, we have
\begin{align}\label{PU1}
\left(\sL^{(\alpha)}_\kappa \tilde u\right)^{x_0}_R=R^{-\alpha}\sL^{(\alpha)}_{\kappa_R} \tilde u^{x_0}_R.
\end{align}
Hence, by Lemma \ref{Var}, there is a constant $c=c(\alpha,\gamma,d,\kappa_0,\kappa_1,\lambda_D)>0$ such that
\begin{align*}
\big\|\big(\sL^{(\alpha)}_\kappa \tilde u\big)^{x_0}_R\big\|_{\gamma; B_1}
\leq R^{-\alpha} \|\sL^{(\alpha)}_{\kappa_R} \tilde u^{x_0}_R\|_{C^\gamma}
\leq c R^{-\alpha} \|\tilde u^{x_0}_R\|_{C^{\alpha+\gamma}}.
\end{align*}
Since $\tilde u^{x_0}_R=\chi_2 u^{x_0}_R$ has support in $B_4$, we further have by \eqref{HL6},
\begin{align}\label{PY5}
\big\|\big(\sL^{(\alpha)}_\kappa \tilde u\big)^{x_0}_R\big\|_{\gamma; B_1}
\leq cR^{-\alpha} \|u^{x_0}_R\|_{\alpha+\gamma;B_4}\leq c R^{\theta-\alpha}\| u\|^{(-\theta)}_{\alpha+\gamma; D},
\end{align}
which together with \eqref{GF4} yields
$$
\big\|\big(\sL^{(\alpha)}_\kappa u\big)^{x_0}_R\big\|_{\gamma; B_1}\leq c R^{\theta-\alpha}\left( \| u\|^{(-\theta)}_{\alpha+\gamma; D}
+{\bf 1}_{\{\theta<0\}}\lambda_D^{-\theta}\|u\|_{0; D}\right).
$$
Thus we obtain the desired estimate by \eqref{HL6}.
\end{proof}

Below we show the interior estimates of Dirichlet problems in weighted H\"older spaces.

\bt\label{Main0}
Let $D$ be a bounded domain and $\alpha\in(0,2)$, $\beta\in(0,1)$. Suppose \eqref{Con1} and $b\in C^\beta$.
For given $\gamma\in(0,\beta]$ with $\alpha+\gamma\notin\mN$ and $\theta\in[0,\alpha\wedge 1)$,
let $u\in \mB^{(-\theta)}_{\alpha+\gamma}(D)$ satisfy 
\begin{align}\label{EQQ1}
\p_t u=\sL^{(\alpha)}_\kappa u+\1_{\alpha\in[1,2)}b\cdot\nabla u+f\mbox{ in $\mR_+\times D$.} 
\end{align}
If one of the following two conditions holds:
$$
{\it (i) }\ \theta\geq\gamma;
\quad  {\it (ii)}\ \mbox{  $\theta<\gamma$ and }|\kappa(x,z)-\kappa(x,z')|\leq\kappa_1|z-z'|^\gamma,
$$
then there is a constant $c=c(d,\kappa_0,\kappa_1,\alpha,\beta,\gamma,\theta, \lambda_D)>0$ such that for all $T>0$,
\begin{align}\label{INS}
\|u\|_{\mB^{(-\theta)}_{\alpha+\gamma; T}(D)}\leq 
c\left(\|u(0)\|_{\cC^{(-\theta)}_{\alpha+\gamma}(D)}+\|f\|_{\mB^{(\alpha-\theta)}_{\gamma; T}(D)}+\|u\|_{\mB^{(-\theta)}_{0; T}(D)}\right),
\end{align}
provided that the right hand side is finite.
\et
\begin{proof}
For $x_0\in D$, let $R:=d_{D^c}(x_0)/8$ and define
$$
u_R(t,x):=R^{-\alpha}u(R^\alpha t, Rx+x_0),\ \ w_R(t,x):=u_R(t,x)\chi_1(x)
$$
and
$$
\kappa_R(x,z):=\kappa(Rx+x_0,Rz),\ \ b_R(x):=b(Rx+x_0),\  f_R(t, x):=f(R^\alpha t,Rx+x_0).
$$
By definitions and scaling, it is easy to see that
\begin{align}\label{EQ1}
\p_t w_R=\sL^{(\alpha)}_{\kappa_R} w_R+{\bf 1}_{\alpha\in[1,2)} R^{\alpha-1}b_R\cdot\nabla w_R+f_R\chi_1+g_R\ \ in \ \ \mR_+\times \mR^d,
\end{align}
where
\begin{align*}
g_R&:=\chi_1\sL^{(\alpha)}_{\kappa_R}((1-\chi_2) u_R)
+[\chi_1,\sL^{(\alpha)}_{\kappa_R}](\chi_2 u_R)\\
&\quad-R^{\alpha-1}u_R{\bf 1}_{\alpha\in[1,2)}b_R\cdot\nabla\chi_1=:g^{(1)}_R+g^{(2)}_R+g^{(3)}_R.
\end{align*}
Noticing that
$$
u_R(R^{-\alpha}t,x)=R^{-\alpha} u(t,Rx+x_0)=:R^{-\alpha}u^{x_0}_R(t),
$$
by the global Schauder's estimate \eqref{N88}, we have
\begin{align}\label{JH93}
\begin{split}
&R^{-\alpha}\sup_{t\in[0,T]}\|u^{x_0}_R(t)\|_{\alpha+\gamma; B_1}
\leq\|u_R\|_{L^\infty_{TR^{-\alpha}}(C^{\alpha+\gamma}(B_1))}\leq\|w_R\|_{\bB^{\alpha+\gamma}_{TR^{-\alpha}}}\\
&\quad\lesssim \|w_R(0)\|_{C^{\alpha+\gamma}}+\|w_R\|_{\bB^0_{TR^{-\alpha}}}
+\|f_R\chi_1\|_{\bB^\gamma_{TR^{-\alpha}}}+\|g_R\|_{\bB^\gamma_{TR^{-\alpha}}}.
\end{split}
\end{align}
Fix $\theta\in[0,\alpha\wedge 1)$. Let us estimate the right hand side of \eqref{JH93}.
First of all, it is easy to see that by \eqref{HL6},
\begin{align}\label{EP1}
\|w_R(0)\|_{C^{\alpha+\gamma}}\leq \|u_R(0)\|_{\alpha+\gamma; B_2}
\lesssim R^{\theta-\alpha}\|u(0)\|^{(-\theta)}_{\alpha+\gamma; D}
\end{align}
and
\begin{align}\label{EP2}
\begin{split}
\|w_R\|_{\bB^0_{TR^{-\alpha}}}&\leq\|u_R\|_{L^\infty_{TR^{-\alpha}}(C^0(B_2))}\lesssim
R^{\theta-\alpha}\|u\|_{\mB^{(-\theta)}_{0;T}(D)},\\
\|f_R\chi_1\|_{\bB^\gamma_{TR^{-\alpha}}}&\lesssim\|f_R\|_{L^\infty_{TR^{-\alpha}}(C^\gamma(B_2))}
\lesssim R^{\theta-\alpha}\|f\|_{\mB^{(\alpha-\theta)}_{\gamma; T}(D)}.
\end{split}
\end{align}
Next we estimate $\|g_R\|_{\bB^\gamma_{TR^{-\alpha}}}\leq\|g^{(1)}_R\|_{\bB^\gamma_{TR^{-\alpha}}}
+\|g^{(2)}_R\|_{\bB^\gamma_{TR^{-\alpha}}}+\|g^{(3)}_R\|_{\bB^\gamma_{TR^{-\alpha}}}$. 
Since the time variable does not play any role in the following calculations, 
we drop it and estimate $\|g_R\|_{C^\gamma}$. 
For $g^{(1)}_R$,  noticing that
$$
\sL^{(\alpha)}_{\kappa_R}((1-\chi_2) u_R)(x)=\sL^{(\alpha)}_{\kappa}((1-\chi^{x_0}_{2R})u)(Rx+x_0),
$$
by the definition of $g^{(1)}_R$, we have for any $\eps\in(0,1)$,
\begin{align*}
\|g^{(1)}_R\|_{C^\gamma}&\lesssim \left\|\left(\sL^{(\alpha)}_{\kappa}((1-\chi^{x_0}_{2R})u)\right)^{x_0}_R\right\|_{\gamma; B_2}
\stackrel{\eqref{GF4}}{\lesssim} R^{\theta-\alpha}\|u\|^{(-\theta)}_{\theta; D}
\stackrel{\eqref{Interpo2},\eqref{HL6}}{\lesssim} R^{\theta-\alpha}\Big(\eps\|u\|^{(-\theta)}_{\alpha+\gamma; D}
+c_\eps\|u\|^{(-\theta)}_{0; D}\Big).
\end{align*}
For $g^{(2)}_R$, by Lemma  \ref{N7} with $R=1$ there and Young's inequality,
for any $\eps\in(0,1)$, there is a constant $c_\eps>0$ such that
\begin{align*}
\|g^{(2)}_R\|_{C^\gamma}&\leq \eps\|\chi_2u_R\|_{C^{\alpha+\gamma}}+c_\eps\|\chi_2u_R\|_{C^{0}}
\lesssim\eps\|u_R\|_{\alpha+\gamma; B_6}+ c_\eps\|u_R\|_{0; B_6}\\
&\leq R^{-\alpha}\left(\eps\|u^{x_0}_R\|_{\alpha+\gamma; B_6}+ c_\eps\|u^{x_0}_R\|_{0; B_6}\right)
\stackrel{\eqref{HL6}}{\lesssim} R^{\theta-\alpha}\left(\eps\|u\|^{(-\theta)}_{\alpha+\gamma; D}+c_\eps\|u\|^{(-\theta)}_{0; D}\right).
\end{align*}
For $g^{(3)}_R$, by \eqref{Interpo2} and \eqref{HL6} again, we have
\begin{align*}
\|g^{(3)}_R\|_{C^\gamma}\leq c\|u_R\|_{\gamma; B_2}
&\leq\eps\|u_R\|_{\alpha+\gamma; B_2}+c_\eps\|u_R\|_{0; B_2} 
\lesssim R^{\theta-\alpha}\left(\eps\|u\|^{(-\theta)}_{\alpha+\gamma; D}+c_\eps \|u\|^{(-\theta)}_{0; D}\right).
\end{align*}
Combining the above calculations, we get for any $\eps\in(0,1)$,
\begin{align}\label{EP3}
\begin{split}
\|g_R\|_{\bB^\gamma_{TR^{-\alpha}}}
&\leq R^{\theta-\alpha}\left(\eps\|u\|_{\mB^{(-\theta)}_{\alpha+\gamma;T}(D)}
+c_\eps\|u\|_{\mB^{(-\theta)}_{0;T}(D)}\right).
\end{split}
\end{align}
By \eqref{JH93}, \eqref{EP1}, \eqref{EP2} and \eqref{EP3}, we obtain that for any $\eps\in(0,1)$,
\begin{align*}
R^{-\theta}\sup_{t\in[0,T]}\|u^{x_0}_R(t)\|_{\alpha+\gamma; B_1}
&\leq \eps\|u\|_{\mB^{(-\theta)}_{\alpha+\gamma;T}(D)}+
c\left(\|u(0)\|^{(-\theta)}_{\alpha+\gamma; D}+\|f\|_{\mB^{(\alpha-\theta)}_{\gamma;T}(D)}\right)+c_\eps\|u\|_{\mB^{(-\theta)}_{0;T}(D)},
\end{align*}
which implies \eqref{INS} by Lemma \ref{Le21} and choosing $\eps$ small enough.
\end{proof}

When $\alpha\in(0,1)$ and the drift $b$ is non-zero, as explained in the introduction, by scaling equation \eqref{EQ1}, one sees that $R^{\alpha-1}$ will blow up
as $R\to 0$. Therefore we have to choose suitable $\theta$ to eliminate the factor $R^{\alpha-1}$ appearing in \eqref{EQ1}. 
Moreover, in order to show the existence of a solution to
the supercritical Dirichlet problem, we shall use the viscosity approximation method. We have

\bt\label{Main44}
Let $D$ be a bounded domain and $\alpha,\beta\in(0,1)$ with $\alpha+\beta>1$. Suppose \eqref{Con1} and $b\in C^\beta$.
For $\gamma\in(1-\alpha,\beta]$, $\theta\in[0,1)$ and $\nu>0$, 
let $u\in \mB^{(-\theta)}_{1+\gamma}(D)$ satisfy 
$$
\p_t u=\nu\Delta^{1/2}u+\sL^{(\alpha)}_\kappa u+b\cdot\nabla u+f\mbox{ in $\mR_+\times D$.} 
$$
If in addition for some $\kappa_1>0$, 
$$
\ |\kappa(x,z)-\kappa(x,z')|\leq\kappa_1|z-z'|^\gamma,
$$
then there is a constant $c=c(\nu,d,\kappa_0,\alpha,\gamma,\theta, \lambda_D)>0$ such that for all $T>0$,
\begin{align}\label{INSS}
\|u\|_{\mB^{(-\theta)}_{1+\gamma; T}(D)}\leq 
c\left(\|u(0)\|_{\cC^{(-\theta)}_{1+\gamma}(D)}+\|f\|_{\mB^{(1-\theta)}_{\gamma; T}(D)}+\|u\|_{\mB^{(-\theta)}_{0; T}(D)}\right),
\end{align}
and there are $\theta_0=\theta_0(\alpha,\beta,\gamma)>0$ and 
a constant $c=c(d,\kappa_0,\kappa_1,\alpha,\beta,\gamma,\theta_0, \lambda_D)>0$ such that for all $T>0$ and $\nu>0$,
\begin{align}\label{Uni}
\|u\|_{\mB^{(\theta_0)}_{\alpha+\gamma;T}(D)}\leq 
c\left(\|u(0)\|_{\cC^{(\theta_0)}_{\alpha+\gamma}(D)}+\|f\|_{\mB^{(\alpha+\theta_0)}_{\gamma;T}(D)}+\|u\|_{\mB^{(0)}_{0;T}(D)}\right).
\end{align}
\et
\begin{proof}
(i) By \eqref{INS} with $\alpha=1$ and \eqref{HL5}, we have
\begin{align*}
\|u\|_{\mB^{(-\theta)}_{1+\gamma; T}(D)}&\leq c\left(\|u(0)\|_{\cC^{(-\theta)}_{1+\gamma}(D)}+\|\sL^{(\alpha)}_\kappa u\|_{\mB^{(1-\theta)}_{\gamma; T}(D)}
+\|f\|_{\mB^{(1-\theta)}_{\gamma; T}(D)}+\|u\|_{\mB^{(-\theta)}_{0; T}(D)}\right)\\
&\leq c\left(\|u(0)\|_{\cC^{(-\theta)}_{1+\gamma}(D)}+\|u\|_{\mB^{(1-\theta-\alpha)}_{\alpha+\gamma; T}(D)}
+\|f\|_{\mB^{(1-\theta)}_{\gamma; T}(D)}+\|u\|_{\mB^{(-\theta)}_{0; T}(D)}\right).
\end{align*}
By \eqref{GF44} and Young's inequality, we further have for any $\eps\in(0,1)$,
$$
\|u\|_{\mB^{(-\theta)}_{1+\gamma; T}(D)}
\leq c\|u(0)\|_{\cC^{(-\theta)}_{1+\gamma}(D)}+\eps\|u\|_{\mB^{(-\theta)}_{1+\gamma; T}(D)}
+c_\eps\|u\|_{\mB^{(1+\gamma-\theta)}_{0; T}(D)}
+c\|f\|_{\mB^{(1-\theta)}_{\gamma; T}(D)}+\|u\|_{\mB^{(-\theta)}_{0; T}(D)},
$$
which in turn gives \eqref{INSS} by $\|u\|_{\mB^{(1+\gamma-\theta)}_{0; T}(D)}\leq\lambda_D^{1+\gamma}\|u\|_{\mB^{(-\theta)}_{0; T}(D)}$.
\\
\\
(ii) To show the uniform estimate \eqref{Uni}, we follow the proof of Theorem \ref{Main0}. For $x_0\in D$, let $R:=d_{D^c}(x_0)/8$ and 
$u_R,w_R,\kappa_R, b_R$ and $f_R$ be defined as in Theorem \ref{Main0}.
By definitions, it is easy to see that
$$
\p_t w_R=\nu R^{\alpha-1}\Delta^{1/2}w_R+\sL^{(\alpha)}_{\kappa_R} w_R+ R^{\alpha-1}b_R\cdot\nabla w_R+f_R\chi_1+g_R\ \ in \ \ \mR_+\times \mR^d,
$$
where
\begin{align*}
g_R&:=\Big\{\nu R^{\alpha-1}\chi_1\Delta^{1/2}((1-\chi_2)u_R)+\chi_1\sL^{(\alpha)}_{\kappa_R}((1-\chi_2) u_R)\Big\}\\
&\quad+\nu R^{\alpha-1}[\chi_1,\Delta^{1/2}](\chi_2u_R)+[\chi_1,\sL^{(\alpha)}_{\kappa_R}](\chi_2 u_R)\\
&\quad-R^{\alpha-1}u_Rb_R\cdot\nabla\chi_1=:g^{(1)}_R+g^{(2)}_R+g^{(3)}_R+g^{(4)}_R.
\end{align*}
Let $m$ be as in \eqref{N888}. By the global Schauder's estimate \eqref{N888}, we have
\begin{align}\label{JH33}
\begin{split}
&R^{-\alpha}\sup_{t\in[0,T]}\|u^{x_0}_R(t)\|_{\alpha+\gamma; B_1}
\leq\|u_R\|_{L^\infty_{TR^{-\alpha}}(C^{\alpha+\gamma}(B_1))}\leq\|w_R\|_{\bB^{\alpha+\gamma}_{TR^{-\alpha}}}\\
&\quad\lesssim \|w_R(0)\|_{C^{\alpha+\gamma}}
+(1+R^{m(\alpha-1)})\|w_R\|_{\bB^0_{TR^{-\alpha}}}+\|f_R\chi_1\|_{\bB^\gamma_{TR^{-\alpha}}}+\|g_R\|_{\bB^\gamma_{TR^{-\alpha}}}.
\end{split}
\end{align}
Here and below, the constant contained in $\lesssim$ is independent of $R\in(0,\lambda_D), x_0\in D$ and  $\nu,\eps\in(0,1)$.
As in the proof of Theorem \ref{Main0}, noticing that 
$$
g^{(1)}_R(x)=\chi_1\left\{\Delta^{1/2}+\sL^{(\alpha)}_\kappa\right\}\left((1-\chi^{x_0}_{2R})u\right)(Rx+x_0),
$$
by \eqref{GF4} with $\theta=0$ there, we have
$$
\|g^{(1)}_R\|_{C^\gamma}\lesssim \left\|\left(\left\{\Delta^{1/2}+\sL^{(\alpha)}_\kappa\right\}((1-\chi^{x_0}_{2R})u)\right)^{x_0}_R\right\|_{\gamma; B_2}\lesssim 
(R^{-1}+R^{-\alpha})[u]_{0;D}.
$$
Fix $\eta\in(0,\alpha)$. By Lemma \ref{N7}, \eqref{Interpo2} and Young's inequality, we have for any $\eps\in(0,1)$,
\begin{align*}
\|g^{(2)}_R\|_{C^\gamma}&\lesssim R^{\alpha-1}\|\chi_2 u_R\|_{C^{\gamma+\eta}}\leq R^{\alpha-1}\|u_R\|_{\gamma+\eta; B_4}\\
&\leq \eps\|u_R\|_{\alpha+\gamma; B_4}+c_\eps R^{(\alpha-1)(\alpha+\gamma)/(\alpha-\eta)}\|u_R\|_{0; B_4}\\
&\lesssim R^{-\theta-\alpha}\left(\eps\|u\|^{(\theta)}_{\alpha+\gamma; D}+c_\eps R^{\theta+(\alpha-1)(\alpha+\gamma)/(\alpha-\eta)}\|u\|_{0; D}\right),
\end{align*}
and 
\begin{align*}
\|g^{(3)}_R\|_{C^\gamma}&\leq \eps\|u_R\|_{\alpha+\gamma;B_4}+c_\eps\|u_R\|_{0;B_4}
\lesssim R^{-\theta-\alpha}\left(\eps\|u\|^{(\theta)}_{\alpha+\gamma; D}+c_\eps\|u\|_{0; D}\right).
\end{align*}
For $g^{(4)}_R$,  by \eqref{Interpo2} again, we have
\begin{align*}
\|g^{(4)}_R\|_{C^\gamma}\leq cR^{\alpha-1}\|u_R\|_{\gamma; B_2}
&\leq\eps\|u_R\|_{\alpha+\gamma; B_2}+c_\eps R^{(\alpha-1)(\alpha+\gamma)/\alpha}\|u_R\|_{0; B_2}\\ 
&\lesssim R^{-\theta-\alpha}\left(\eps\|u\|^{(\theta)}_{\alpha+\gamma; D}+c_\eps R^{\theta+(\alpha-1)(\alpha+\gamma)/\alpha} \|u\|_{0; D}\right).
\end{align*}
Combining the above calculations and choosing $\theta$ large enough, we get for any $\eps\in(0,1)$,
\begin{align}\label{EP33}
\begin{split}
\|g_R\|_{\bB^\gamma_{TR^{-\alpha}}}
&\leq R^{-\theta-\alpha}\left(\eps\|u\|_{\mB^{(\theta)}_{\alpha+\gamma;T}(D)}+c_\eps\|u\|_{\mB^{(0)}_{0;T}(D)}\right),
\end{split}
\end{align}
and also
\begin{align}\label{EP22}
(1+R^{m(\alpha-1)})\|w_R\|_{L^\infty_{TR^{-\alpha}}(C^0)}\lesssim
\lambda_D^{m(\alpha-1)+\theta}R^{-\theta-\alpha}\|u\|_{\mB^{(0)}_{0;T}(D)}.
\end{align}
By \eqref{EP1}, \eqref{EP2}, \eqref{JH33}, \eqref{EP33} and \eqref{EP22}, we obtain that for all $\eps\in(0,1)$,
\begin{align*}
R^{\theta}\sup_{t\in[0,T]}\|u^{x_0}_R(t)\|_{\alpha+\gamma; B_1}
&\leq \eps\|u\|_{\mB^{(\theta)}_{\alpha+\gamma;T}(D)}+
c_\eps\Big(\|u(0)\|^{(\theta)}_{\alpha+\gamma; D}+\|f\|_{\mB^{(\alpha+\theta)}_{\gamma;T}(D)}+\|u\|_{\mB^{(0)}_{0;T}(D)}\Big),
\end{align*}
which implies  the desired estimate by Lemma \ref{Le21} and choosing $\eps$ small enough.
\end{proof}

\section{Probabilistic representation for Dirichlet problem}

Let $\Omega$ be the space of  all c\`adl\`ag functions from $\mR_+$ to $\mR^d$, which is  endowed with the Skorokhod topology.
Let $X_t(\omega)=\omega_t$ be the coordinate process over $\Omega$ and
$\{\sF_t^0; t\geq 0\}$ the  natural filtration generated by $X$. 
For Borel sets $A,D\subset \mR^d$, denote by $\sigma_A$ and $\tau_D$ the hitting time of $A$ and the first exit time of $D$ respectively, i.e.,
$$
\sigma_A:=\inf\{t\geq 0: X_t\in A\},\ \ \tau_D:=\inf\{t\geq 0: X_t\notin D\}.
$$
The following relation will be used frequently in the strong Markov property: for $A\subset D$,
\begin{align}
\tau_D=\sigma_A+\tau_D\circ\theta_{\sigma_A}\ \mbox{ on }  \sigma_A<\tau_D,
\end{align}
where $\theta_t(\omega):=\omega_{t+\cdot}$ is the usual shift operator on $\Omega$. 

Below we shall present a general probabilistic representation for Dirichlet problem of nonlocal parabolic operators. 
Let $J(x,z)$ be a nonnegative measurable function on 
$\mR^d\times\mR^d$, which is  a L\'evy jump kernel and satisfies that for some $\vartheta\in(0,2]$,
$$
J^{(\vartheta)}_0(x):=\int_{\mR^d}(1\wedge|z|^\vartheta)J(x,z)\dif z<\infty,\ \forall x\in\mR^d.
$$
Let $\sL_J$ be the nonlocal L\'evy operator associated with $J$, that is, for any $f\in C^2_b(\mR^d)$,
$$
\sL_J f(x):=\int_{\mR^d}(f(x+z)-f(x)-\1_{\{|z|\leq 1\}}z\cdot\nabla f(x))J(x,z)\dif z.
$$
Throughout this section we always assume that
\begin{enumerate}[{\bf (MP)} ]
\item $J^{(\vartheta)}_0$ and $b$ are bounded measurable, and 
for each $x\in\mR^d$, there is a unique probability measure $\mP_x$ over $\sF^0_\infty$ so that for any $u\in \bB^{\vartheta\vee 1}_T$ with
$\p_t u\in \bB^{0}_T$ (see \eqref{SPAC11} for a definition of space $\bB^{\gamma}_T$),
\begin{align}\label{EE3}
u(t,X_t)-u(0,x)-\int^t_0(\p_su+\sL_J u+b\cdot\nabla u)(s,X_s)\dif s
\end{align}
is an $\sF^0_t$-martingale starting from zero under $\mP_{x}$. In particular, $\{X,\mP_{x};{x\in\mR^d}\}$  forms a family of strong Markov processes (see \cite{Et-Ku}).
We shall denote by $\sF_t$ the augmentation filtration
of $\sF^0_t$ with respect to $(\mP_x)_{x\in\mR^d}$, and $P_tf(x):=\E_x f(X_t)$. 
Moreover, we also require $P_t: C^\infty_c(\mR^d)\to C_b(\mR^d)$.
\end{enumerate}

\br\rm\label{RE-Feller}
Let $\beta\in(0,1)$ with $\alpha+\beta>1$.
Under {\bf (H$^\beta_\kappa$)} and $b\in C^\beta$, the above assumption {\bf (MP)} is satisfied for $J(x,z):=\kappa(x,z)/|z|^{d+\alpha}$.
In fact, since the coefficients are bounded continuous, the existence of martingale solutions is well-known (see \cite[p.536, Theorem 2.31]{Ja-Shi}). 
We only show the uniqueness.  
Given $T>0$ and $f\in C^\infty_c(\mR^d)$, let $u\in \bB^{\alpha+\beta}_T$ be the unique solution of 
the following nonlocal equation (see Remark \ref{RE35}),
$$
\p_t u+\sL^{(\alpha)}_{\kappa} u+b\cdot\nabla u+f=0,\ \ u(T,x)=0.
$$
By \eqref{EE3}, we get
$$
u(0,x)=\mE_x\left(\int^T_0 f(X_s)\dif s\right).
$$
Since the left hand side does not depend on $\mP_x$, the uniqueness follows by \cite[Corollary 6.2.4]{St-Va}. 
Moreover, again by Remark \ref{RE35}, we have $P_t: C_c^\infty(\mR^d)\to C_b(\mR^d)$.
\er

\subsection{Probabilistic representation}
The following L\'evy system is a crucial tool in the study of jump processes (see \cite{Ch-Ki-So3} and \cite{Ch-Zh1} for a proof).
\bt
Let $f$ be a nonnegative measurable function on $\mR_+\times\mR^d\times\mR^d$
vanishing on $\{(s,x,y): x=y\}$. 
For any  $x\in\mR^d$ and stopping time $\tau$, it holds that
\begin{align}\label{Levy}
\mE_x\left(\sum_{s\leq \tau}f(s,X_{s-},X_s)\right)=\mE_x\left(\int^\tau_0\!\!\int_{\mR^d}f(s,X_s,z)J(X_s,z-X_s)\dif z\dif s\right).
\end{align}
\et
The L\'evy system will be used in many situations as follows, which exhibits the main feature of jump processes.
\bl\label{Le52}
Let $D\subset\mR^d$ be an open subset and $A\subset D^c$ a measurable subset with ${\rm dist}(A,D)>0$. For any $x\in D$, we have
\begin{align}\label{Le}
\mP_x(X_{\tau_D}\in A)=\mE_x\left(\int^{\tau_D}_0\!\!\int_A J(X_s,z-X_s)\dif z\dif s\right).
\end{align}
In particular,  if the Lebesgue measure of $A$ is zero, then
$$
\mP_x(X_{\tau_D}\in A)=0,\ \ x\in D.
$$
\el
\begin{proof}
Since ${\rm dis}(A,D)>0$, we have
$$
{\bf 1}_{\{0<\tau_D\}}{\bf 1}_{A}(X_{\tau_D})=\sum_{0<s\leq\tau_D}{\bf 1}_{\bar D}(X_{s-}){\bf 1}_{A}(X_s).
$$
Since $D$ is open, $\mP_x(\tau_D>0)=1$ for $x\in D$.
By the L\'evy system \eqref{Levy} with $f(s,x,y)={\bf 1}_{\bar D}(x){\bf 1}_A(y)$, we have
$$
\mP_x(X_{\tau_D}\in A)=\mE_x\left(\sum_{0<s\leq\tau_D}{\bf 1}_{\bar D}(X_{s-}){\bf 1}_{A}(X_s)\right)=
\mE_x\left(\int^{\tau_D}_0\!\!\int_A{\bf 1}_{\bar D}(X_s)J(X_s,z-X_s)\dif z\dif s\right),
$$
which gives \eqref{Le} because $X_s\in D$ for $s<\tau_D$.
\end{proof}

The following result states the quasi-left continuity of $X$,
which is essentially contained in \cite[page 70, Theorem 4]{Chung}. We sketch it's proof.
\bl\label{Le51}
Let $D\subset\mR^d$ be a bounded domain and $D_n\uparrow D$. 
For each $x\in\mR^d$ and  $\mP_x$-almost all $\omega\in\Omega$, it holds that
$$
\tau_{D_n}(\omega)\uparrow\tau_D(\omega),\ \ X_{\tau_{D_n}(\omega)}(\omega)\to X_{\tau_D(\omega)}(\omega).
$$
\el
\begin{proof}
Let $\tau_\infty:=\sup_n\tau_n$. Obviously, $\tau_\infty\leq \tau_D$. 
Moreover, we also have $X_{\tau_{D_n}}\to X_{\tau_\infty}$ a.s., which follows by the same argument as in the proof of \cite[page 70, Theorem 4]{Chung}.
Now since $X_{\tau_n}\in D^c_n$ for each $n\in\mN$
and $D_n\uparrow D$, we must have $X_{\tau_\infty}\in D^c$, which implies that
$\tau_\infty=\tau_D$ a.s.
\end{proof}
To present the probabilistic representation and a maximum principle of nonlocal Dirichlet problem, we introduce the following class of functions pair:
for $\gamma>0$,
\begin{align}\label{AL5}
\bH^\gamma(D):=
\left\{
(u,f)\,\Bigg|\,
\begin{aligned}
& u\in L^\infty_{loc}(\mR_+; C^\gamma_{loc}(D)\cap L^\infty(\mR^d))\cap C(\mR_+\times(\p D)^c)\\
& \quad\p_t u\in L^\infty_{loc}(\mR_+; C_{loc}(D)),\quad f\in L^\infty_{loc}(\mR_+\times D)
\end{aligned}
\right\}.
\end{align}
\bt\label{Th44}
Let $D$ be a bounded domain, and $b\in L^\infty_{loc}(D)$ and $(u,f)\in \bH^{\vartheta\vee 1}(D)$ satisfy
\begin{align}\label{NH2}
\p_t u=\sL_J u+b\cdot\nabla u+f\ \mbox{ on }\ \mR_+\times D.
\end{align}
Suppose that $\p D$ has Lebesgue zero measure, and  one of the following two conditions holds:
\begin{enumerate}[{\rm (i)}]
\item $f\geq 0$ or $f\in L^\infty_{loc}(\mR_+; C_b(D))$, and for all  $x\in D$, $\mP_x(X_{\tau_D}\in\p D)=0$.
\item $f\in L^\infty_{loc}(\mR_+; C_b(D))$, and $\p D=\Gamma_0\cup\Gamma_1$, where $\Gamma_0$ and $\Gamma_1$ are two disjoint measurable sets,
and for all $x\in D$, $\mP_x(X_{\tau_D}\in\Gamma_0)=0$ and $u\in C((0,\infty)\times (D\cup\Gamma_1))$ with $u|_{(0,\infty)\times\Gamma_1}=0$.
\end{enumerate}
Then for all $x\in\mR^d$ and $t>0$, it holds that
\begin{align}\label{Rep}
\begin{split}
u(t,x)&=\mE_x \Big(u(0, X_t); t<\tau_D\Big)+\mE_x\left(\int^{t\wedge\tau_{D}}_0f(t-s,X_s)\dif s\right)\\
&\quad+\mE_x \Big(u(t-\tau_D, X_{\tau_{D}}); \tau_D\leq t, X_{\tau_D}\notin\p D\Big).
\end{split}
\end{align}
In particular, we have the following maximum principle:
\begin{align}\label{Max}
\|u\|_{L^\infty_T(C^0(D))}\leq \|u(0)\|_{C^0(D)}+T\|f\|_{L^\infty_T(C^0(D))}+\|u\|_{L^\infty_T(C^0(\bar D^c))}.
\end{align}
\et
\begin{proof}
For $x\notin D$,  there is nothing to prove since  $\mP_x(\tau_D=0)=1$ by the definition of $\tau_D$.
We assume $x\in D$. Let $\rho_\eps$ be a family of mollifiers with support in $B_\eps$. 
Define $u_\eps:=u*\rho_\eps$. Let $D_n\uparrow\!\uparrow D$. Fix $t>0$. Applying \eqref{EE3} to
function $(s,y)\mapsto u_\eps(t-s,y)$, we have
$$
\mE_x u_\eps(t-t\wedge\tau_{D_n},X_{t\wedge\tau_{D_n}})=u_\eps(t,x)-\mE_x\int^{t\wedge\tau_{D_n}}_0
(\p_su_\eps-\sL_J u_\eps-b\cdot\nabla u_\eps)(t-s,X_s)\dif s.
$$
Fix $y\in D_n$. For $\eps<{\rm dist}(D_n, D^c)/2=:\delta_n$, we drop the time variable and write
\begin{align*}
\sL_J u_\eps(y)&=\int_{|z|<\delta_n}(u_\eps(y+z)-u_\eps(y)-\1_{\{|z|\leq 1\}}z\cdot\nabla u_\eps(y))J(y,z)\dif z\\
&+\int_{|z|\geq\delta_n}(u_\eps(y+z)-u_\eps(y)-\1_{\{|z|\leq 1\}}z\cdot\nabla u_\eps(y))J(y,z)\dif z=:I^{(1)}_\eps(y)+I^{(2)}_\eps(y).
\end{align*}
Since $u\in C^{\vartheta\vee 1}_{loc}(D)$, by $\|J^{(\vartheta)}_0\|_\infty<\infty$ and the dominated convergence theorem, we have
$$
\lim_{\eps\to 0}I^{(1)}_\eps(y)=\int_{|z|<\delta_n}(u(y+z)-u(y)-\1_{\{|z|\leq 1\}}z\cdot\nabla u(y))J(y,z)\dif z,
$$
and by $u\in L^\infty(\mR^d)\cap C(\mR^d\setminus\p D)$,
$$
\lim_{\eps\to 0}I^{(2)}_\eps(y)=\int_{|z|\geq\delta_n}(u(y+z)-u(y)-\1_{\{|z|\leq 1\}}z\cdot\nabla u(y))J(y,z)\dif z.
$$
Moreover, it is easy to see that
$$
\sup_{\eps\in(0,\delta_n)}\|(\p_s u_\eps-\sL_J u_\eps-b\cdot\nabla u_\eps\|_{L^\infty([0,t]\times D_n)}<\infty.
$$
Hence, by the dominated convergence theorem and \eqref{NH2}, we have
$$
\mE_x\int^{t\wedge\tau_{D_n}}_0
(\p_su_\eps-\sL_J u_\eps-b\cdot\nabla u_\eps)(t-s,X_s)\dif s
\stackrel{\eps\to 0}{\to} \mE_x\int^{t\wedge\tau_{D_n}}_0f(t-s,X_s)\dif s.
$$
Notice that for fixed $n$, by Lemma \ref{Le52} and $|\p D|=0$ we have
$$
\mP_x(X_{\tau_{D_n}}\in \p D)=0.
$$
Since $u\in C(\mR_+\times(\mR^d\setminus\p D))$ is bounded, by the dominated convergence theorem again, we have
$$
\lim_{\eps\to 0}\mE_x u_\eps(t-t\wedge\tau_{D_n},X_{t\wedge\tau_{D_n}})=\mE_x u(t-t\wedge\tau_{D_n},X_{t\wedge\tau_{D_n}}).
$$
Combining the above limits, we get
\begin{align}\label{Rep0}
\mE_x u(t-t\wedge\tau_{D_n},X_{t\wedge\tau_{D_n}})=u(t,x)-\mE_x\int^{t\wedge\tau_{D_n}}_0f(t-s,X_s)\dif s.
\end{align}
Write
\begin{align}\label{HS1}
\mE_x u(t-t\wedge\tau_{D_n},X_{t\wedge\tau_{D_n}})
=\mE_x \Big(u(t-\tau_{D_n},X_{\tau_{D_n}});\tau_{D_n}\leq t\Big)
+\mE_x \Big(u(0,X_t);\tau_{D_n}>t\Big).
\end{align}
Since $\tau_{D_n}\uparrow\tau_D$, we have
\begin{align}\label{HS2}
\lim_{n\to\infty}\mE_x \Big(u(0,X_t);\tau_{D_n}>t\Big)=\mE_x \Big(u(0,X_t);\tau_{D}>t\Big).
\end{align}
Let 
$$
\Omega_0=\cup_{n}\{\omega: X_{\tau_{D_n}(\omega)}(\omega)\in  \bar D^c\},\ \ 
\Omega^c_0=\cap_{n}\{\omega: X_{\tau_{D_n}(\omega)}(\omega)\in  D\setminus \bar D^c_n\}.
$$ 
Notice that for any $\omega\in\Omega_0$, there is a $n_0$ such that for all $n\geq n_0$,
$X_{\tau_{D_{n}}(\omega)}(\omega)=X_{\tau_{D}(\omega)}(\omega)\in \bar D^c$. Hence,
\begin{align}\label{HS3}
\lim_{n\to\infty}\mE_x \Big(u(t-\tau_{D_n},X_{\tau_{D_n}}); \tau_{D_n}\leq t, \Omega_0\Big)=
\mE_x\Big( u(t-\tau_{D}, X_{\tau_D});\tau_{D}\leq t,X_{\tau_D}\notin \p D\Big).
\end{align}
Moreover, by Lemma \ref{Le51}, we have
\begin{align}\label{HS6}
\mP_x(\omega: X_{\tau_{D_n}(\omega)}(\omega)\to X_{\tau_D(\omega)}(\omega))=1\Rightarrow \mP_x(\Omega^c_0)=\mP_x(X_{\tau_D}\in\p D).
\end{align}
\\
(i) If $\mP_x(X_{\tau_D}\in\p D)=0$, then 
$\mP_x(\Omega^c_0)=0$. Hence, combining this with \eqref{HS1}-\eqref{HS3}, and taking limits $n\to\infty$ for \eqref{Rep0}, 
by the monotone convergence theorem or the dominated convergence theorem, we obtain \eqref{Rep}.
\\
\\
(ii) Write $\Omega^c_0=\Omega_1+\Omega_2$, where
$$
\Omega_1:=\Omega^c_0\cap\{X_{\tau_D}\in\Gamma_0\},\ \Omega_2:=\Omega^c_0\cap\{X_{\tau_D}\in\Gamma_1\}.
$$
Fix $x\in D$. By the assumption we have
\begin{align}\label{WP1}
\lim_{n\to\infty}\mE_x \Big(u(t-\tau_{D_n},X_{\tau_{D_n}});\tau_{D_n}\leq t, \Omega_1\Big)=0.
\end{align}
To treat $\mE_x \Big(u(t-\tau_{D_n},X_{\tau_{D_n}});\tau_{D_n}\leq t, \Omega_2\Big)$, notice that there is a countable set $\sT_x\subset\mR_+$ so that
\begin{align}\label{HS4}
\mP_x(\tau_D=t)=0,\ t\notin\sT_x.
\end{align}
Since $u\in C((0,\infty)\times (D\times\Gamma_1))$ and $u|_{(0,\infty)\times\Gamma_1}=0$, for $t\notin\sT_x$ we have
\begin{align*}
\lim_{n\to\infty}\mE_x \Big(u(t-\tau_{D_n},X_{\tau_{D_n}});\tau_{D_n}\leq t, \Omega_2\Big)
&=\lim_{n\to\infty}\mE_x \Big(u(t-\tau_{D_n},X_{\tau_{D_n}});\tau_{D_n}\leq\tau_D< t, \Omega_2\Big)=0,
\end{align*}
which together  with \eqref{Rep0}-\eqref{WP1} yields \eqref{Rep} for $t\notin\sT_x$.
Next we assume $t\in\sT_x$, and choose $t_n\notin\sT_x$ so that $t_n\downarrow t$. By what we have proved and \eqref{HS4}, it holds that
\begin{align*}
u(t_n,x)&=\mE_x\Big(u(t_n-\tau_D, X_{\tau_{D}}); t_n>\tau_D, X_{\tau_D}\notin\p D\Big)+\mE_x\Big(u(0, X_{t_n}); t_n<\tau_D\Big)\\
&\quad+\mE_x\left(\int^{t_n}_0f(s,X_{t_n-s})\1_{\{t_n-s<\tau_D\}}\dif s\right)=:I^{(1)}_n+I^{(2)}_n+I^{(3)}_n.
\end{align*}
For $I^{(1)}_n$, by $u\in C([0,\infty)\times\bar D^c)$, we have
\begin{align*}
\lim_{n\to\infty}I^{(1)}_n=\mE_x\Big(u(t-\tau_D, X_{\tau_{D}}); t\geq \tau_D, X_{\tau_D}\notin\p D\Big).
\end{align*}
For $I^{(2)}_n$ and $I^{(3)}_n$, since $t\mapsto X_t$ is right continuous and $u,f$ are bounded continuous in $D$, 
by the dominated convergence theorem, we have
$$
\lim_{n\to 0}I^{(2)}_n=\mE_x\Big(u(0, X_{t}); t<\tau_D\Big),\quad 
\lim_{n\to 0}I^{(3)}_n=\mE_x\left(\int^{t}_0f(s,X_{t-s})\1_{\{t-s<\tau_D\}}\dif s\right).
$$
Combining the above limits, we obtain \eqref{Rep} for any $t>0$.
\end{proof}
\br\rm
The above case (ii) will be used in the supercritical case. Notice that the condition $u|_{(0,\infty)\times\Gamma_1}=0$ 
can be replaced with that for each $z\in\p D$, the limit $\lim_{t\to 0} u(t,z)$ exists and is denoted by $u^0_{\Gamma_1}(z)$. If so, we need an extra term 
$\mE_x(u^0_{\Gamma_1}(X_{\tau_D});\tau_D=t, X_{\tau_D}\in\Gamma_1)$ in \eqref{Rep}.
\er
We also need the following simple estimate.  
\bl
Let $D$ be a bounded domain and $f\in C_b(\mR^d)\cap C^{2}_{loc}(D)$ and $g\in C_b(D)$ satisfy
$$
\sL_J f(x)+b(x)\cdot\nabla f(x)\leq g(x),\ \ x\in D.
$$
Then for any $x\in D$, it holds that
\begin{align}\label{Dyn}
\mE_x f(X_{\tau_D})\leq f(x)+\mE_x\int^{\tau_D}_0g(X_s)\dif s.
\end{align}
\el
\begin{proof}
Let $D_n\uparrow\!\uparrow D$. As in the proof of Theorem \ref{Th44} and by the assumption, we have
$$
\mE_x f(X_{\tau_{D_n}})=f(x)+\mE_x\int^{\tau_{D_n}}_0(\sL_J f+b\cdot\nabla f)(X_s)\dif s
\leq f(x)+\mE_x\int^{\tau_{D_n}}_0g(X_s)\dif s.
$$
By taking limits $n\to\infty$, we obtain the desired estimate.
\end{proof}

\subsection{Estimates of exit times}

The following lemma is well-known (see \cite{So-Vo}, \cite{Ch-Ki-So3} and \cite{Ch-Zh1}). For the reader's convenience, we provide the proof here. 
\bl\label{Le501}
Suppose that for some $\alpha\in(0,2)$ and $\kappa_1,\kappa_2>0$,
$$
\1_{\{|z|\leq 1\}}J(x,z)\leq\kappa_1|z|^{-d-\alpha}, \ \ \int_{|z|\geq 1}J(x,z)\dif z\leq\kappa_2.
$$
Then there is a constant $c_0>0$ such that for all $x_0\in\mR^d$ and $t>0$, $\eps\in(0,1)$,
$$
\mP_{x_0}(\tau_{B_\eps(x_0)}<t)\leq c_0\left(t/\eps^{\alpha}+\|b\|_\infty t/\eps\right).
$$
\el
\begin{proof}
Given $f\in C^2_b(\mR^d)$ with $f(0)=0$ and $f(x)=1$ for $|x|\geq 1$,    set
$$
f_\eps(x):=f((x-x_0)/\eps),\ \ \eps>0.
$$
By \eqref{EE3} and the optional stopping theorem,
\begin{align}
\mP_{x_0}\Big(\tau_{B_\eps(x_0)}<t\Big)\leq \mE_{x_0}f_\eps\Big(X_{\tau_{B_\eps(x_0)}\wedge t}\Big)
=\mE_{x_0}\int^{\tau_{B_\eps(x_0)}\wedge t}_0(\sL_J f_\eps+b\cdot\nabla f_{\eps})(X_{s})\dif s.\label{ERY2}
\end{align}
By the assumption, we have for all $\eps\in(0,1)$,
\begin{align*}
|\sL_J f_\eps(x)|&= \left|\int_{\mR^d}(f_\eps(x+z)-f_\eps(x)-{\bf 1}_{\{|z|\leq 1\}}z\cdot\nabla f_\eps(x))J(x,z)\dif z\right|\\
&\lesssim\int_{|z|\leq \eps}\frac{\|\nabla^2 f_\eps\|_\infty}{|z|^{d+\alpha-2}}\dif z+\int_{1\geq |z|\geq \eps}\frac{\|\nabla f_\eps\|_\infty}{|z|^{d+\alpha-1}}\dif z
+\|f_\eps\|_\infty\kappa_2\\
&\lesssim  (\|\nabla f^2\|_\infty+\|\nabla f\|_\infty)\eps^{-\alpha}+\|f\|_\infty\kappa_2\lesssim\eps^{-\alpha},
\end{align*}
and
$$
|b(x)\cdot\nabla f_\eps(x)|\leq \|b\|_\infty\|\nabla f\|_\infty/\eps.
$$
Substituting them into \eqref{ERY2}, we obtain the desired estimate.
\end{proof}
We need the following moment estimate of the exit time $\tau_D$ from a bounded domain $D$.
\bl
Let $D\subset \mR^d$ be a bounded domain. Suppose that for some $\alpha\in(0,2)$ and $\kappa_0>0$,
$$
J(x,z)\geq\kappa_0|z|^{-d-\alpha}.
$$
Then for any $n\in \N$, there is a constant $c_n$ depending only on $n,\alpha,d,\kappa_0$ such that for all $x\in\mR^d$,
\begin{align}\label{Stop-bounded}
\E_x \tau_D^n\leq c_{n}\lambda_D^{n\alpha}. 
\end{align}
\el
\begin{proof}
Recall that $\lambda_D={\rm diam}(D)$. Let $T:=\inf\{s>0:|\Delta X_s|>\lambda_D\}$.  
For any $t>0$, by the L\'evy system \eqref{Levy} and the assumption, we have
\begin{align*}
\P_x(T\leq t)&=\E_x\left( \sum_{s\leq T\wedge t} {\bf 1}_{\{|\Delta X_s|>\lambda_D\}}\right)
=\E_x\left( \int_0^{T\wedge t}\dif s\int_{|z|>\lambda_D}J(X_s,z)\dif z\right)\\
&\geq c_1\lambda_D^{-\alpha}\E_x (T\wedge t) \geq c_1\lambda_D^{-\alpha} t \P_x(T> t), 
\end{align*}
where $c_1=c_1(\alpha,\kappa_0,d)>0$. This together with $\{T\leq  t\}\subset\{\tau_D\leq  t\}$ implies that for all $x\in D$,
$$
\P_x (\tau_D> t)\leq  \P_x(T> t)\leq 1/(1+c_1\lambda_D^{-\alpha} t). 
$$
Moreover, for any $n\in\mN$, by the Markov property,
\begin{align*}
\mP_x(\tau_D>n t)&=\mP_x\Big(\tau_D>n t|\sF_{(n-1)t}\Big)=\mE_x\Big(\mP_{X_{(n-1)t}}(\tau_D>t); \tau_D>(n-1)t\Big)\\
&\leq 1/(1+c_1\lambda_D^{-\alpha} t)\mP_x(\tau_D>(n-1) t)\leq\cdots\leq 1/(1+c_1\lambda_D^{-\alpha} t)^n.
\end{align*}
Therefore, for any $n\in\mN$, by the change of variable, we have
\begin{align*}
\E_x\tau^n_D&=n\int^\infty_0 t ^{n-1}\mP_x(\tau_D> t)\dif t=n(n+1)^n\int^\infty_0 t^{n-1}\mP_x(\tau_D>(n+1) t)\dif  t\\
&\leq n(n+1)^n\int^\infty_0\frac{ t^{n-1}}{(1+c_1\lambda_D^{-\alpha} t)^{n+1}}\dif  t
=n(n+1)^n(c_1\lambda_D^{-\alpha})^{-n}\int^\infty_0\frac{ t^{n-1}}{(1+ t)^{n+1}}\dif  t,
\end{align*}
which yields the desired estimate.
\end{proof}
\br\rm
From the above two lemmas, if
$$
\1_{\{|z|\leq 1\}}J(x,z)\leq\kappa_1|z|^{-d-\alpha},\ \ J(x,z)\geq\kappa_0|z|^{-d-\alpha},
$$ 
then for any $n\in\mN$, there is a constant $c_n\geq 1$ such that
for all $x\in\mR^d$ and $\eps\in(0,1)$,
\begin{align}\label{EK0}
c^{-1}_n(\1_{b\equiv 0}\eps^\alpha+\1_{b\not=0}\eps^{\alpha\vee 1})^n\leq \mE_x\tau^n_{B_\eps(x)}\leq c_n\eps^{n\alpha}.
\end{align}
Indeed, the upper bound follows by \eqref{Stop-bounded}.  For the lower bound, by Lemma \ref{Le501} we have
$$
\mE_x\tau^n_{B_\eps(x)}\geq t^n\mP_x(\tau_{B_\eps(x)}\geq t)\geq t^n(1-c_0(t/\eps^{\alpha}+\|b\|_\infty t/\eps)).
$$
Taking $t=(\1_{b\equiv 0}\eps^\alpha+\1_{b\not=0}\eps^{\alpha\vee 1})/(4c_0)$, we obtain the lower bound.
\er
The following  two lemmas about the estimate of the first exit time are useful.
\bl\label{Le59}
Let $D\subset\mR^d$ be an open subset. Suppose that there is a constant $c_0>0$ such that for each $x\in D$, there is a neighborhood $Q_x\Subset D$ of $x$ such that
\begin{align}\label{CN1}
\mP_x (X_{\tau_{Q_x}}\in D^c)\geq c_0.
\end{align}
Then for any  $A\subset D^c$,
\begin{align}\label{CN2}
\sup_{x\in D}\mP_x (X_{\tau_D}\in A)\leq \sup_{x\in D}\mP_x(X_{\tau_{Q_x}}\in A)/c_0.
\end{align}
In particular, if $A$ has Lebesgue zero measure, then $\mP_x (X_{\tau_D}\in A)=0$ for each $x\in D$.
\el
\begin{proof}
Since $X_{\tau_D}=X_{\tau_D}\circ\theta_{\tau_{Q_x}}$ on $\tau_{Q_x}<\tau_D$, by the strong Markov property we have
\begin{align*}
\mP_x(X_{\tau_D}\in A)&=\mP_x(X_{\tau_D}\in A; X_{\tau_{Q_x}}\in D^c)+\mP_x(X_{\tau_D}\in A;X_{\tau_{Q_x}}\in D)\\
&=\mP_x(X_{\tau_{Q_x}}\in A)+\mE_x\Big(\mP_{X_{\tau_{Q_x}}}(X_{\tau_D}\in A); X_{\tau_{Q_x}}\in D\Big)\\
&\leq\mP_x(X_{\tau_{Q_x}}\in A)+\sup_{x\in D}\mP_x(X_{\tau_D}\in A)\,\mP_x(X_{\tau_{Q_x}}\in D).
\end{align*}
Taking supremum with respect to  $x\in D$ and by the assumption, we obtain
$$
\sup_{x\in D}\mP_x(X_{\tau_D}\in A)\leq \sup_{x\in D}\mP_x(X_{\tau_{Q_x}}\in A)+\sup_{x\in D}\mP_x(X_{\tau_D}\in A) (1-c_0),
$$
which implies the desired estimate \eqref{CN2}. If $A$ has Lebesgue zero measure, since ${\rm dist}(A,Q_x)>0$,
by Lemma \ref{Le52}, we have $\mP_x(X_{\tau_{Q_x}}\in A)=0$ for each $x\in D$.
So $\mP_x (X_{\tau_D}\in A)=0$.
\end{proof}

\bl\label{Le55}
Let $V, U, D$ be three bounded domains. Suppose that $V\Subset U$, $V\cap D\not=\emptyset$ and for some $\alpha\in(0,2)$ and $\kappa_0>0$,
$$
J(x,z)\geq\kappa_0|z|^{-d-\alpha}.
$$
Then there is a constant $c_0=c_0(U,V,D)>0$ such that for all $x\in V\cap D$,
\begin{align}\label{EG7}
\mE_x\tau_D\leq c_0\mE_x\tau_{U\cap D}.
\end{align}
\el
\begin{proof}
For $x\in V\cap D$, since $\tau_D=\tau_{U\cap D}+\tau_D\circ\theta_{\tau_{U\cap D}}$, by the strong Markov property, we have
\begin{align}\label{ET5}
\begin{split}
\mE_x\tau_D&\leq\mE_x\tau_{U\cap D}+\mE_x(\tau_D; \tau_{U\cap D}<\tau_D)\\
&\leq 2\mE_x\tau_{U\cap D}+\mE_x \Big(\mE_{X_{\tau_{U\cap D}}}\tau_D; \tau_{U\cap D}<\tau_D\Big)\\
&\leq 2\mE_x\tau_{U\cap D}+\mP_x(\tau_{U\cap D}<\tau_D)\sup_{x\in D}\mE_x\tau_D.
\end{split}
\end{align}
Let $f$ be a $C^2$-function with
$$
f(x)=0,\ x\in V, \ f(x)=1,\ x\in U^c.
$$
It is easy to see that
$$
c_1:=\|\sL_J f+b\cdot\nabla f\|_\infty<\infty.
$$
Since $\tau_{U\cap D}<\tau_D$ implies $X_{\tau_{U\cap D}}\in D\setminus U$, by \eqref{Dyn},  we have for $x\in V\cap D$,
$$
\mP_x(\tau_{U\cap D}<\tau_D)\leq \mE_x f(X_{\tau_{U\cap D}})\leq c_1\mE_x\tau_{U\cap D}.
$$
Substituting this into \eqref{ET5} and by \eqref{Stop-bounded}, we obtain \eqref{EG7}.
\end{proof}
The following result states that the process always jumps out a bounded domain without touching the boundary.
\bp\label{Pr512}
Let $D$ be a bounded domain satisfying the uniformly exterior cone condition. 
Let $\alpha\in(0,2)$. Suppose that $b\equiv 0$ if $\alpha\in(0,1)$,  and for some $\kappa_0,\kappa_1,\kappa_2>0$,
\begin{align}\label{HJ4}
J(x,z)\geq \kappa_0|z|^{-d-\alpha},\ \ J(x,z)\1_{\{|z|\leq 1\}}\leq \kappa_1|z|^{-d-\alpha},\ \ \int_{|z|\geq 1}J(x,z)\dif z\leq\kappa_2.
\end{align}
Then for each $x\in D$, it holds that
\begin{align}\label{EE2}
\mP_x(X_{\tau_D}\in \p D)=0.
\end{align}
Moreover, if $\mP_x\circ X^{-1}_s(\dif y)\ll\dif y$ for each $s>0$, then 
\begin{align}\label{EE22}
\mP_x(X_{\tau_D-}\in \p D, X_{\tau_D}\notin \bar D^c)=0.
\end{align}
\ep
\begin{proof}
(i) We first show \eqref{EE2}.
Fix a point $x\in D$. Let $z_x\in\p D$ be such that 
${\rm dist}(x,\p D)=|x-z_x|.$ Set $R:=(|x-z_x|\wedge 1)/3$ and $Q_x:=B_{R}(x)$.
Since $\p D$ satisfies uniformly exterior cone condition, there is a cone $C_\theta$ with vertex $z_x$ and angle $\theta>0$ not depending on the point $x$ 
such that $C_\theta \cap B_R(z_x)\subset D^c$ and $|C_\theta\cap B_R(z_x)|\geq c_0R^{-d}$, where $c_0=c_0(\theta,d)$.  
Noticing that for some $c_1=c_1(\lambda_D)>0$,
$${\rm diam}(B_R(x)\cup B_R(z_x)) \leq c_1  R,
$$
by formula \eqref{Le} and the assumption, we have
\begin{align*}
\mP_x(X_{\tau_{Q_x}}\in D^c)&=\mE_x\int^{\tau_{Q_x}}_0\!\!\!\dif s\!\int_{D^c}J(X_{s}, z-X_{s})\dif z
\geq\kappa_0\mE_x\int^{\tau_{Q_x}}_0\!\!\!\dif s\!\int_{C_\theta \cap B_R(z_x)}\frac{\dif z}{|z-X_s|^{d+\alpha}}\\
&\geq\kappa_0\mE_x\tau_{Q_x}\int_{C_\theta \cap B_R(z_x)}\frac{\dif z}{(c_1R)^{d+\alpha}}
\geq c R^{-\alpha}\mE_x\tau_{Q_x}.
\end{align*}
Since $b\equiv 0$ for $\alpha\in(0,1)$, by \eqref{EK0} with $\eps=R$ and $n=1$, we get
$$
\mP_x(X_{\tau_{Q_x}}\in D^c)\geq c_0>0,
$$
where $c_0$ is independent of $x\in D$.
Thus by Lemma \ref{Le59}, we obtain \eqref{EE2}.
\\
\\
(ii) For \eqref{EE22}, it suffices to show that for any $\eps,T>0$ (see \cite{Ar-Bi-Ca}),
$$
\mP_x(X_{T\wedge\tau_D-}\in\p D, X_{T\wedge\tau_D}\in \mR^d\backslash D^\eps)=0,
$$
where $D^\eps:=\{x\in\mR^d: {\rm dist}(x, D)\leq\eps\}$.
By the L\'evy system \eqref{Levy}, we have
\begin{align*}
\mP_x(X_{T\wedge\tau_D-}\in\p D, X_{\tau_D}\in
\mR^d\backslash D^\eps) 
&=\mE_x\int^{T\wedge\tau_D}_0
{\bf 1}_{\{X_{s}\in\p D\}}\int_{D^c_\eps}J(X_{s}, z-X_{s})\dif z\dif s\\
&\leq c_\eps\mE_x\int^{T}_0{\bf 1}_{\{X_{s}\in\p D\}}\dif s=0.
\end{align*}
The proof is complete.
\end{proof}

\br\rm\label{Re13}
Notice that $J(x,z):=\nu|z|^{-d-1}+\kappa(x,z)|z|^{-d-\alpha'}$ with $\nu>0$ and $0\leq \kappa(x,z)\leq\kappa_1$, $\alpha'\in(0,1)$, satisfies \eqref{HJ4} for $\alpha=1$.
\er

\section{Subcritical and critical cases: Proof of Theorem \ref{MAIN}}
\subsection{Distance functions}
Let $U\subset\mR^d$ be any open subset.  
The aim of this subsection is to prove some basic estimates for $\sL^{(\alpha)}_\kappa d^\theta_{U^c}$ when $U$ is a half space, a ball's complement
and a bounded $C^2$-domain, respectively. 
Fix $\alpha\in(0,2)$.  Throughout this subsection we assume
\begin{align}\label{Con2}
\kappa^{-1}_0\leq\kappa(x,z)\leq\kappa_0,\ \ {\bf 1}_{\alpha=1}\int_{r<|z|<R}z\,\kappa(x,z)\dif z=0,\ 0<r<R<\infty.
\end{align}
\bl\label{Est-g}
Let $Q:=\{x=(x_1,\cdots,x_d)\in\mR^d: x_1>0\}$ be the half space. Under \eqref{Con2}, 
there are constants $\theta_0\in(0,\tfrac{\alpha}{2})$ 
and  $c_0>0$ only depending on $\kappa_0,d,\alpha$ such that for any $\theta\in(0, \theta_0]$,
\begin{align}\label{ES10}
\left(\sL^{(\alpha)}_\kappa d_{Q^c}^{\theta}\right)(x)\leq -c_0d^{\theta-\a}_{Q^c}(x)=-c_0x_1^{\theta-\a},\ \ x\in Q.
\end{align}
\el
\begin{proof}
Let ${\bf e}_1:=(1,0,\cdots,0)$ and for $z=(z_1,\cdots, z_d)$ and $x=(x_1,\cdots, x_d)$,
$$
z_1^*:=(z_2,\cdots, z_d)\in\mR^{d-1},\ \ \tilde\kappa_x(z):=\kappa(x,(x_1z_1,z^*_1)). 
$$
Notice that $d_{Q^c}(x)=(x_1)_+=x_1\vee 0$.  For $x\in Q$, by scaling we have
$$
\left(\sL^{(\alpha)}_\kappa d_{Q^c}^{\theta}\right)(x)=x_1^{\alpha-\theta}\left(\sL^{(\alpha)}_{\tilde \kappa_x} d_{Q^c}^{\theta}\right)({\bf e}_1).
$$
Hence, it suffices to prove \eqref{ES10} for $x={\bf e}_1$ and $\kappa(x,z)=\kappa(z)$.
Noticing that 
$$
\nabla d^\theta_{Q^c}({\bf e}_1)=(\theta,0,\cdots,0),
$$
by definition we have
\begin{align*}
\left(\sL^{(\alpha)}_\kappa d^\theta_{Q^c}\right)({\bf e}_1)
&=\int_{\mR^d} \left((1+z_1)_+^\theta-1-\theta \big(\1_{\a=1}\1_{\{|z|\leq 1\}}+\1_{\alpha\in(1,2)}\big)z_1   \right)\frac{\kappa(z)}{|z|^{d+\a}} \dif z\\
&=\int_{\{z_1>-1\}} \left((1+z_1)^\theta-1-\theta \big(\1_{\a=1}\1_{\{|z|\leq 1\}}+\1_{\alpha\in(1,2)}\big) z_1  \right)\frac{\kappa(z)}{|z|^{d+\a}} \dif z\\
&-\theta  \int_{\{z_1\leq -1\}} \big(\1_{\a=1}\1_{\{|z|\leq 1\}}+\1_{\alpha\in(1,2)}\big)\frac{z_1\kappa(z)}{|z|^{d+\a}} \dif z
-\int_{\{z_1\leq -1\}}\frac{\kappa(z)}{|z|^{d+\a}} \dif z=:\sum_{i=1}^3I_i(\theta).
\end{align*}
For $I_1(\theta)$, by the change of variables ($z_1=s, z_1^*=sy$ with $y\in\mR^{d-1}$), we have
\begin{align*}
|I_1(\theta)|
&\leq\kappa_0\int^\infty_{-1}\!\!\!\int_{\mR^{d-1}} \left|(1+z_1)^\theta-1-\theta z_1 \big(\1_{\a=1}\1_{\{|z|\leq 1\}}+\1_{\alpha\in(1,2)}\big)  \right|
\frac{\dif z_1\dif z^*_1}{(|z_1|^2+|z^*_1|^2)^{(d+\a)/2}}\\
&=\kappa_0\int^\infty_{-1}\!\!\!\int_{\mR^{d-1}} \left|(1+s)^\theta-1-\theta s \Big(\1_{\a=1}{\bf 1}_{|s|^2(1+|y|^2)\leq 1}+\1_{\alpha\in(1,2)}\Big)  \right|
\frac{s^{-1-\alpha}\dif s\dif y}{(1+|y|^2)^{(d+\a)/2}}.
\end{align*}
By elementary calculations, one sees that
$$
\lim_{\theta\to 0}|I_1(\theta)|=0.
$$
For $I_2(\theta)$, as above we have
\begin{align*}
|I_2(\theta)|&\leq \kappa_0\theta\int^{-1}_{-\infty}\!\int_{\mR^{d-1}} \big(\1_{\a=1}\1_{\{|z|\leq 1\}}+\1_{\alpha\in(1,2)}\big)
\frac{|z_1|\dif z_1\dif z^*_1}{(|z_1|^2+|z^*_1|^2)^{(d+\a)/2}}\\
&=\kappa_0\theta\int^\infty_1\!\int_{\mR^{d-1}} \big(\1_{\a=1}{\bf 1}_{s^2(1+|y|^2)\leq 1}+\1_{\alpha\in(1,2)}\big)
\frac{s^{-\alpha}\dif s\dif y}{(1+|y|^2)^{(d+\a)/2}}\leq c\kappa_0\theta.
\end{align*}
For $I_3(\theta)$, we have
\begin{align*}
I_3(\theta)&\leq -\kappa_0^{-1}\int^{-1}_{-\infty}\!\int_{\mR^{d-1}}\frac{\dif z_1\dif z_1^*}{(|z_1|^2+|z^*_1|^2)^{(d+\a)/2}}
=-\kappa_0^{-1}\int^\infty_1\frac{\dif s}{s^{1+\alpha}}\int_{\mR^{d-1}} 
\frac{\dif y}{(1+|y|^2)^{(d+\a)/2}}.
\end{align*}
Combining the above calculations, we obtain
$$
\left(\sL^{(\alpha)}_\kappa d^\theta_{Q^c}\right)({\bf e}_1) \leq (-c_1\kappa_0^{-1}+c_\theta\kappa_0),
$$
where $c_\theta\to 0$ as $\theta\to 0$ and $c_1=c_1(\alpha,d)$.
Thus we obtain the desired estimate by letting $\theta$ small enough.
\end{proof}
\br\rm
When $\kappa\equiv1$, by more careful calculations, we have (for example, see \cite{Ro-Se1})
\begin{align}\label{ES34}
\left(\Delta^{\frac{\alpha}{2}}d^\theta_{Q^c}\right)({\bf e}_1)=
\left\{
\begin{array}{ll}
-c_0,&\ \theta\in(0,\tfrac{\alpha}{2}),\\
0,&\theta=\tfrac{\alpha}{2},\\
c_1,&\theta\in(\tfrac{\alpha}{2},\alpha),
\end{array}
\right.
\end{align}
where $c_0,c_1>0$ only depends on $d,\alpha,\theta$.
\er
Next we show the same estimate for any ball's complement.
\bl\label{Est-F1}
Under \eqref{Con2}, there exists a $\theta_0\in(0,\frac{\alpha}{2})$ such that for all  $\theta\in(0,\theta_0]$,  
there are $\delta, c>0$ only depending on $\theta,\kappa_0,d,\alpha$ such that for any $x_0\in\mR^d$ and $R>0$,
$$
\left(\sL^{(\alpha)}_\kappa d^\theta_{B_R(x_0)}\right)(x)\leq -c\cdot d_{B_R(x_0)}^{\theta-\alpha}(x),\ \ x\in B_{R(1+\delta)}(x_0)\setminus B_R(x_0).
$$
\el
\begin{proof}
Noticing that $d_{B_R(x_0)}(x)=(|x-x_0|-R)_+$, we have
$$
d_{B_R(x_0)}(x)=R\, d_{B_1}((x-x_0)/R),
$$
and by scaling,
$$
\left(\sL^{(\alpha)}_\kappa d^\theta_{B_R(x_0)}\right)(x)=R^{\theta-\alpha}\left(\sL^{(\alpha)}_{\kappa_R} d^\theta_{B_1}\right)((x-x_0)/R),
$$
where $\kappa_R(x,z):=\kappa(x,Rz)$.
Hence, without loss of generality, we may assume $x_0=0$ and $R=1$. 
For simplicity we write
$$
h(x):=d^\theta_{B_1}(x)=(|x|-1)_+^\theta, \ g(x)=(x_1-1)_+^\theta.
$$
Let $r\geq 0$ and $x_r=(1+r, 0,\cdots,0)\in \mR^d$. 
Noticing that
$$
h(x_r)=g(x_r),\ \ \nabla h(x_r)=\theta r^{\theta-1} {\mathbf {e}_1}=\nabla g(x_r),
$$ 
we have
\begin{align*}
 \sL^{(\alpha)}_\kappa (h-g)(x_r)&=  \int_{\mR^d} (h(x_r+z)-g(x_r+z)) \frac{\kappa(x_r, z)}{|z|^{d+\alpha}}\dif z.
\end{align*}
Notice the following elementary inequality: for $\theta\in(0,1)$ and $a,b\geq 0$,
$$
((1+a)^2+b^2)^{1/2}\leq 1+a+b^2,\ (a+b)^\theta-a^\theta\leq(a^{\theta-1}b)\wedge b^\theta.
$$
For any $z=(z_1, z^*_1)\in \mR^d$, if $r+z_1\geq 0$, then
\begin{align*}
h(x_r+z)-g(x_r+z)&=\left\{\big((1+r+z_1)^2+|z_1^*|^2\big)^{1/2}-1\right\}^{\theta}-|r+z_1|^\theta\\
&\leq (|r+z_1|+|z_1^*|^2)^\theta-|r+z_1|^\theta\leq(|r+z_1|^{\theta-1}|z_1^*|^2)\wedge |z_1^*|^{2\theta};
\end{align*}
if $r+z_1<0$, then $g(x_r+z)=0$ and
\begin{align*}
h(x_r+z)&\leq\Big(\big((1+r+z_1)^2+|z_1^*|^2\big)^{1/2}-1\Big)_+^{\theta}
\leq\Big(\big((1-|r+z_1|)^2+|z_1^*|^2\big)^{1/2}-1\Big)_+^{\theta}\\
&\leq\Big(\big(1+|r+z_1|^2+|z_1^*|^2\big)^{1/2}-1\Big)^{\theta}\leq |r+z_1|^{2\theta}+|z_1^*|^{2\theta}.
\end{align*}
Hence,
\begin{align}
\sL^{(\alpha)}_\kappa (h-g)(x_r)&\leq\kappa_0  \int_{\{r+z_1\geq 0\}}(|r+z_1|^{\theta-1}|z_1^*|^2)\wedge |z_1^*|^{2\theta}\frac{\dif z}{|z|^{d+\alpha}}\no\\
&\quad+\kappa_0  \int_{\{r+z_1<0\}}(|r+z_1|^{2\theta}+|z_1^*|^{2\theta})\frac{\dif z}{|z|^{d+\alpha}}\no\\
&\leq \kappa_0
\left(\Big(\tfrac{r}{2}\Big)^{\theta-1}\int_{|z|<r/2}\frac{|z|^2\dif z}{|z|^{d+\alpha}}
+\int_{|z|>r/2}\frac{|z|^{2\theta}\dif z}{|z|^{d+\alpha}}\right)\no\\
&\quad+\kappa_0 r^{2\theta} \int_{\{r+z_1<0\}}\frac{\dif z}{|z|^{d+\alpha}}
+\kappa_0\int_{\{r+z_1<0\}}\frac{\dif z}{|z|^{d+\alpha-2\theta}}\no\\
&\leq c_d\kappa_0
\left(\tfrac{r^{1+\theta-\alpha}}{2-\alpha}
+\tfrac{r^{2\theta-\alpha}}{\alpha-2\theta}\right).\label{Rg}
\end{align}
On the other hand, by Lemma \ref{Est-g}, there are $\theta_0\in(0,\frac{\alpha}{2})$ and $c_0>0$ such that for all $\theta\in(0,\theta_0]$,
$$
\sL^{(\alpha)}_{\kappa} g(x_r)\leq -c_0 r^{\theta-\alpha},
$$
which together with \eqref{Rg} yields
$$
\sL^{(\alpha)}_{\kappa} h(x_r)\leq -c_0 r^{\theta-\alpha}+ c_1r^{2\theta-\alpha}\leq-c_0 r^{\theta-\alpha}/2,
$$
provided $r\in(0,(c_0/(2c_1))^{1/\theta}\wedge 1)$. By rotational invariance, we obtain the desired estimate.
\end{proof}
\br\rm
By \eqref{ES34}, for any $\theta\in(0,\tfrac{\alpha}{2})$, there are $c,\delta>0$ such that for all $x_0\in\mR^d$ and $R>0$,
\begin{align}\label{ES44}
\left(\Delta^{\tfrac{\alpha}{2}}d^\theta_{B_R(x_0)}\right)(x)\leq -c\cdot d_{B_R(x_0)}^{\theta-\alpha}(x),\ \ x\in B_{R(1+\delta)}(x_0)\setminus B_R(x_0).
\end{align}
\er

Now we extend the above lemma to general bounded $C^2$-domain.
\bl\label{Le66}
Let $D$ be a bounded $C^2$-domain and $b$ a bounded vector field. Under \eqref{Con2},
there exists a $\theta_0\in(0,\frac{\alpha}{2})$ such that for all $\theta\in(0,\theta_0]$, 
there are $\eps,c_0\in(0,1)$ only depending on $\theta,\|b\|_\infty,\kappa_0,d,\alpha, D$ such that
$$
\left(\sL^{(\alpha)}_\kappa d^\theta_{D^c}+{\bf 1}_{\alpha\in[1,2)}b\cdot\nabla d^\theta_{D^c}\right)(x)
\leq -c_0d^{\theta-\alpha}_{D^c}(x),\ \ x\in D\setminus D_\eps,
$$
where $D_\eps:=\{x\in D: d_{D^c}(x)>\eps\}$.
\el
\begin{proof}
Since $D$ is a bounded $C^2$-domain, for each boundary point $z_0\in\p D$, there is an exterior tangent ball 
$B_r(y_0)\subset D^c$ which touches $D$ at $z_0$, where the radius $r$ does not depend on $z_0$.
Moreover, it is well known that $d_{D^c}$ is a $C^2$-function on $D\setminus D_\eps$ provided $\eps$ small enough
(see \cite[Lemma 14.16]{Gi-Tr}). Fix $x_0\in D\setminus D_\eps$. Let $z_0\in\p D$ be the unique boundary point such that
$$
d_{D^c}(x_0)=|x_0-z_0|.
$$
Let $B:=B_r(y_0)$ be the exterior tangent ball of $D$ at point $z_0$ so that $d_{D^c}(x_0)=d_B(x_0)$.
Let $\theta,\delta$ be as in Lemma \ref{Est-F1}.
Without loss of generality we assume $\delta=\eps=\delta\wedge\eps$. By Lemma \ref{Est-F1}, we have
\begin{align}\label{ES7}
\left(\sL^{(\alpha)}_\kappa d^\theta_B\right)(x_0)\leq -c_0d^{\theta-\alpha}_B(x_0)=-c_0d^{\theta-\alpha}_{D^c}(x_0).
\end{align}
Since $B$ lies in the outside of $D$, it is easy to see that
\begin{align}\label{PU2}
d_{D^c}(x)\leq d_B(x),\ \forall x\in\mR^d,
\end{align}
which together with $d_{D^c}(x_0)=d_B(x_0)$ and the maximum principle yields
\begin{align}\label{PU3}
\nabla d_{D^c}(x_0)=\nabla d_B(x_0).
\end{align}
Hence, by $d_{D^c}(x_0)=d_B(x_0)$ and \eqref{PU2}, \eqref{PU3},
$$
\left(\sL^{(\alpha)}_\kappa d^\theta_{D^c}\right)(x_0)\leq \left(\sL^{(\alpha)}_\kappa d^\theta_B\right)(x_0)\stackrel{\eqref{ES7}}{\leq}-c_0d^{\theta-\alpha}_{D^c}(x_0).
$$
On the other hand, it is easy to see that
$$
|b\cdot\nabla d^\theta_{D^c}|(x_0)\leq\theta\|b\|_\infty d^{\theta-1}_{D^c}(x_0).
$$
Therefore,
$$
\left(\sL^{(\alpha)}_\kappa d^\theta_{D^c}\right)(x_0)\leq (-c_0+\1_{\alpha\in[1,2)}\theta\|b\|_\infty\eps^{\alpha-1})d^{\theta-\alpha}_{D^c}(x_0),\ \ x_0\in D\setminus D_\eps.
$$
Choosing $\theta$ small enough we get the desired estimate.
\end{proof}
\br\rm
By \eqref{ES44}, for any $\theta\in(0,\frac{\alpha}{2})$, there are $\eps,c_0>0$ such that
\begin{align}\label{ES54}
\left(\Delta^{\frac{\alpha}{2}} d^\theta_{D^c}\right)(x)\leq -c_0d^{\theta-\alpha}_{D^c}(x),\ \ x\in D\setminus D_\eps.
\end{align}
\er
\subsection{A maximum principle in weighted H\"older spaces}
In this subsection we show a maximum principle in weighted H\"older spaces by using the barrier function in Lemma \ref{Le66}.
\bl\label{Le77}
Let $D$ be a bounded $C^2$-domain, and $b$ a bounded measurable vector field and $\gamma\in(0,1)$. Under \eqref{Con2} and {\bf (MP)},
there exists a $\theta_0\in(0,\frac{\alpha}{2})$ such that for all $\theta\in(0,\theta_0]$,
there is a constant $c=c(d,\alpha,\theta,\kappa_0,D)>0$ such that for any pair of $(u,f)\in \bH^{\alpha+\gamma}(D)$ satisfying 
$$
\p_t u=\sL^{(\alpha)}_{\kappa} u+1_{\alpha\in[1,2)}b\cdot\nabla u+f\ \mbox{ on }\ \mR_+\times D,\ \ u|_{\mR_+\times D^c}=0,\ u(0)=0,
$$
it holds that for all $T>0$, 
\begin{align}\label{ES3}
\|u\|_{\mB^{(-\theta)}_{0;T}(D)}\leq c(1+T)\|f\|_{\mB^{(\alpha-\theta)}_{0;T}(D)}.
\end{align}
\el
\begin{proof}
Let $\theta,\eps,c_0$ be as in Lemma \ref{Le66}. Let $d_x=d_{D^c}(x)$ and $D_\eps:=\{x\in D: d_x>\eps\}$.  
Fix $T>0$ and define
$$
\cN:=\|u\|_{L^\infty_T(C^0(D_\eps))}+\|f\|_{L^\infty_T(\cC^{(\alpha-\theta)}_0(D))},\ \ w(x):=\cN d^\theta_x/(c_0\eps^\theta).
$$
Then by the definition of $\mB^{(\alpha-\theta)}_{0;T}(D)$ (see \eqref{SPAC1}) and Lemma \ref{Le66}, we have
$$
|f(t,x)|\leq \cN d^{\theta-\alpha}_x\leq -\left(\sL^{(\alpha)}_\kappa+{\bf 1}_{\alpha\in[1,2)}b\cdot\nabla\right) w(x)\ 
\mbox{ in $[0,T]\times D\setminus D_\eps$},
$$
and so,
$$
\p_t (u-w)=\p_t u\leq \sL^{(\alpha)}_\kappa(u-w)+{\bf 1}_{\alpha\in[1,2)}b\cdot\nabla (u-w) \ \mbox{ in $[0,T]\times D\setminus D_\eps$},
$$
with
$$
u-w\leq 0\ \mbox{ in $[0,T]\times (D^c\cup D_\eps)$}.
$$
Since (i) in Theorem \ref{Th44} is satisfied (see Proposition \ref{Pr512}), by  \eqref{Rep}, it is easy to see that
\begin{align*}
u(t,x)\leq w(x)=\cN d^\theta_x/(c_0\eps^\theta)\ \mbox{ in $[0,T]\times D\setminus D_\eps$}.
\end{align*}
Hence, by the definition of $\cN$,
\begin{align}\label{PU4}
\sI_\eps:=\sup_{t\in[0,T]}\sup_{x\in D\setminus D_\eps}d^{-\theta}_x|u(t,x)|
&\leq \Big(\|u\|_{L^\infty_T(C^0(D_\eps))}+\|f\|_{L^\infty_T(\cC^{(\alpha-\theta)}_0(D))}\Big)/(c_0\eps^\theta).
\end{align}
On the other hand, by the maximum principle \eqref{Max}, we also have for any $\delta\in(0,\eps]$,
\begin{align}\label{PU5}
\|u\|_{L^\infty_T(C^0(D_\delta))}\leq T\|f\|_{L^\infty_T(C^0(D_\delta))}+\|u\|_{L^\infty_T(C^0(D\setminus D_\delta))}.
\end{align}
Since $\|u\|_{L^\infty_T(C^0(D_\eps))}\leq\|u\|_{L^\infty_T(C^0(D_\delta))}+\|u\|_{L^\infty_T(C^0(D\setminus D_\delta))}$
and $\sI_\delta\leq\sI_\eps$, by \eqref{PU4} and \eqref{PU5} we further have
\begin{align*}
\sI_\eps&\leq \left(2\|u\|_{L^\infty_T(C^0(D\setminus D_\delta))}+T\|f\|_{L^\infty_T(C^0(D_\delta))}+\|f\|_{L^\infty_T(\cC^{(\alpha-\theta)}_0(D))}\right)/(c_0\eps^\theta)\\
&\leq \left(2\delta^\theta\sI_\eps+(T\delta^{\theta-\alpha}+1)\|f\|_{L^\infty_T(\cC^{(\alpha-\theta)}_0(D))}\right)/(c_0\eps^\theta).
\end{align*}
By letting $\delta$ be small enough, we get
$$
\sI_\eps\leq c_\eps(1+T)\|f\|_{L^\infty_T(\cC^{(\alpha-\theta)}_0(D))},
$$
which together with \eqref{PU5} yields the desired estimate.
\end{proof}
\br\rm\label{Re68}
Let $\gamma\in (0,1)$ and $(u,f)\in \bH^{\alpha+\gamma}(D)$  solve the following Dirichlet problem:
$$
\p_t u=\Delta^{\frac{\alpha}{2}}u+f\mbox{ in $\mR_+\times D$},\ u|_{\mR_+\times D^c}=0.
$$
As above, by \eqref{ES54}, one can show that for any $\theta\in(0,\frac{\alpha}{2})$, 
there exists a constant $c=c(d,\alpha,D,\theta)>0$ such that
for any $T>0$,
$$
\|u\|_{\mB^{(-\theta)}_{0;T}(D)}\leq c(1+T)\|f\|_{\mB^{(\alpha-\theta)}_{0;T}(D)}.
$$
\er
\subsection{Proof of Theorem \ref{MAIN}}
We need the following solvability of fractional Dirichlet problem, whose proof is given in the appendix. 
The main novelty here is that $f$ is not necessarily bounded near the boundary. 
\bt\label{Th52}
Let $D$ be a bounded $C^2$-domain and $\alpha\in(0,2)$. 
For any $\theta\in(0,\tfrac{\alpha}{2})$, $\gamma\in(0,1)$ and
$f\in \mB^{(\alpha-\theta)}_\gamma(D)$, there is a unique $u\in \mB^{(-\theta)}_{\alpha+\gamma}(D)$
so that
\begin{align}\label{Lap0}
u(t,x)=\int^t_0\left(\Delta^{\frac{\alpha}{2}} u(s,x)+f(s,x)\right)\dif s,\ \ t\geq 0,\ x\in D,
\end{align}
or simply,
\begin{align}\label{Lap}
\p_t u=\Delta^{\frac{\alpha}{2}} u+f \mbox{ in } \mR_+\times D,\ \ u|_{\mR_+\times D^c}=0,\ \ u(0,\cdot)=0.
\end{align}
\et

Now we can use the continuity method to prove Theorem \ref{MAIN}.
\begin{proof}[Proof of Theorem \ref{MAIN}] 
By considering $\tilde u=u-\varphi$, without loss of generality, we may assume $\varphi=0$.
Fix $T>0$. Let $\theta$ be as in Lemma \ref{Le77}. Define a Banach space
$$
\mA^{(\theta)}_{\alpha,\gamma;T}(D):=\Big\{u\in \mB^{(-\theta)}_{\alpha+\gamma;T}(D), 
\p_t u\in \mB^{(\alpha-\theta)}_{\gamma;T}(D)\Big\}.
$$
Let $\sL^{(\alpha)}_{\kappa,b}:=\sL^{(\alpha)}_{\kappa}+{\bf 1}_{\alpha\in[1,2)}b\cdot\nabla$. For $\alpha\in[1,2)$, by \eqref{HL6}, we have
$$
\|b\cdot\nabla u\|^{(\alpha-\theta)}_{\gamma;D}\leq 
c\|b\|^{(\alpha-1)}_{\gamma; D}\|\nabla u\|^{(1-\theta)}_{\gamma;D}\lesssim\|u\|^{(-\theta)}_{1+\gamma;D}
\lesssim\|u\|^{(-\theta)}_{\alpha+\gamma;D},
$$
which together with Theorem \ref{Th0}, yields that
$$
\sL^{(\alpha)}_{\kappa,b}: \mB^{(-\theta)}_{\alpha+\gamma;T}(D)\to \mB^{(\alpha-\theta)}_{\gamma;T}(D).
$$
For $\tau\in[0,1]$, consider the operator $\sT_\tau: \mA^{(\theta)}_{\alpha,\gamma;T}(D)\to \mB^{(\alpha-\theta)}_{\gamma;T}(D)$:
$$
\sT_\tau u:=\p_tu-(1-\tau)\Delta^{\frac{\alpha}{2}}u-\tau\sL^{(\alpha)}_{\kappa,b}u.
$$
By Theorem \ref{Main0} and Lemma \ref{Le77}, there is a constant $c>0$ independent of $\tau\in[0,1]$ such that 
$$
\|u\|_{\mB^{(-\theta)}_{\alpha+\gamma;T}(D)}\leq c\|\sT_\tau u\|_{\mB^{(\alpha-\theta)}_{\gamma;T}(D)},
$$
and by \eqref{HL5},
$$
\|\p_tu\|_{\mB^{(\alpha-\theta)}_{\gamma;T}(D)}
\leq \|\sT_\tau u\|_{\mB^{(\alpha-\theta)}_{\gamma;T}(D)}+ 
\|\Delta^{\frac{\alpha}{2}}u\|_{\mB^{(\alpha-\theta)}_{\gamma;T}(D)}+\|\sL^{(\alpha)}_{\kappa,b}u\|_{\mB^{(\alpha-\theta)}_{\gamma;T}(D)}
\leq c\|\sT_\tau u\|_{\mB^{(\alpha-\theta)}_{\gamma;T}(D)}.
$$
Since $\sT_0$ is an onto mapping from $\mA^{(\theta)}_{\alpha,\gamma;T}(D)$ to $\mB^{(\alpha-\theta)}_{\gamma;T}(D)$ by Theorem \ref{Th52}, 
by \cite[Theorem 5.2]{Gi-Tr}, $\sT_1$ is also an onto mapping from $\mA^{(\theta)}_{\alpha,\gamma;T}(D)$ to $\mB^{(\alpha-\theta)}_{\gamma;T}(D)$. 
Thus we get the existence. As for the uniqueness, it follows by Lemma \ref{Le77}.

To show the probabilistic representation, for $f\in \mB^{(\alpha-\theta)}_{\gamma}(D)$, we let
$$
f_+:=f\vee 0,\ \ f_-:=(-f)\vee 0.
$$
It is easy to see that $f_+, f_-\in \mB^{(\alpha-\theta)}_{\gamma;T}(D)$.
Let $u_+$ and $u_-$ be the solutions of \eqref{IBP0} corresponding  to $(f_+,\varphi)$ and $(f_-, 0)$, respectively.
By the uniqueness, one has
$$
u=u_+-u_-.
$$
Moreover, by Remark \ref{RE-Feller} and Theorem \ref{Th44}  we have
$$
u_+(t,x)=\mE_x \Big(\varphi(X_{t}){\bf 1}_{\tau_{D}>t}\Big)
+\mE_x\left(\int^{t\wedge\tau_{D}}_0f_+(t-s,X_s)\dif s\right)
$$
and
$$
u_-(t,x)=\mE_x\left(\int^{t\wedge\tau_{D}}_0f_-(t-s,X_s)\dif s\right).
$$
Thus we get \eqref{EP0}. As for \eqref{ES2}, it follows by \eqref{ES3} and \eqref{EP0}. 
\end{proof}

\section{Supercritical case: Proof of Theorems \ref{MAIN1} and \ref{MAIN2}}

In the following we always assume $\alpha\in(0,1)$ and $b\in C^\beta$ with $\beta\in(1-\alpha,1)$.
\subsection{Boundary probabilistic estimates}
\begin{figure}[h!]\small
\centering
\includegraphics
[height=0.31\textwidth,width=1.0\textwidth]{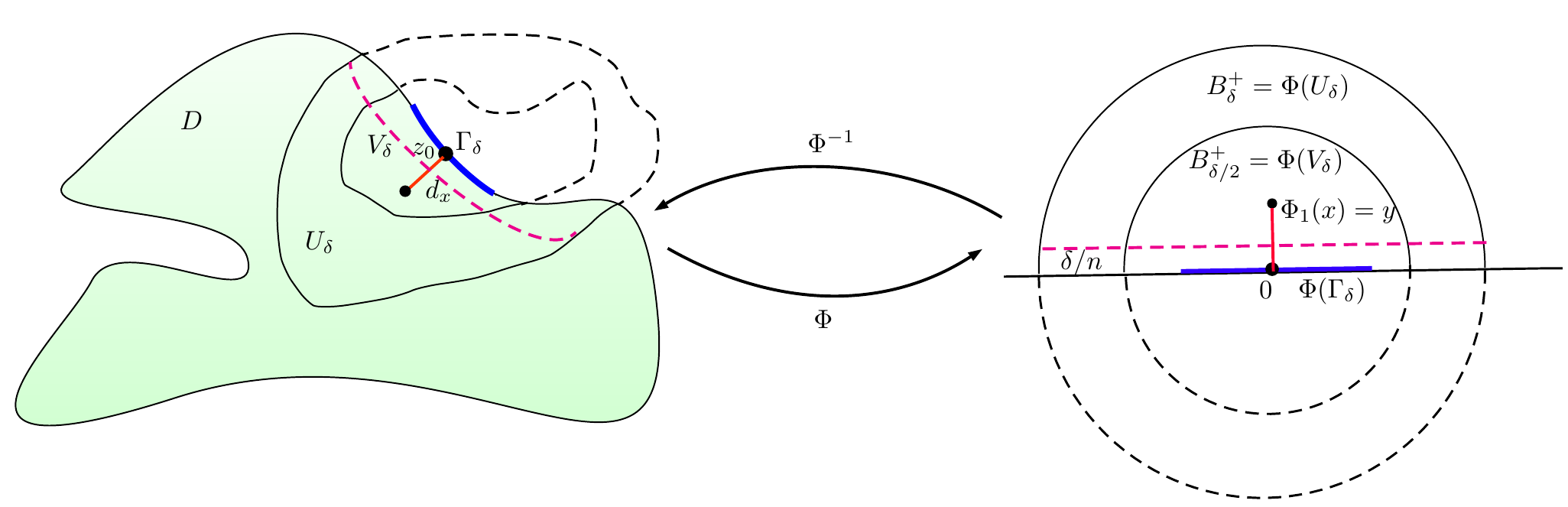}
\caption{Flatten out the boundary}\label{fig2}
\end{figure}

In this subsection we first show some estimates about the first
exit time and the exit position of the process $X_t$ from a domain in the supercritical case.
Let $D$ be a bounded domain with  $C^2$-boundary. More precisely, for any 
$z_0\in \p D$, there are a neighborhood $W\subset\mR^d$ of $z_0$ and a $C^2$-function $\phi:\mR^{d-1}\to\mR$ such that 
(upon relabeling and reorienting the coordinates axes if necessary)
$$
W\cap D=\{x\in W: x_1>\phi(x_2,\cdots, x_d)\}.
$$
Define two maps $\Phi, \Psi: \mR^d\to \mR^d$
\begin{align}\label{SRT}
\left\{
\begin{aligned}
&\Phi_1(x):=x_1-\phi(x_2,\cdots,x_d),\\
&\Phi_j(x):=x_j,\ \ j=2,\cdots, d;
\end{aligned}
\right.
\quad
\left\{
\begin{aligned}
&\Psi_1(y):=y_1+\phi(y_2,\cdots,y_d),\\
&\Psi_j(y):=y_j,\ \ j=2,\cdots, d.
\end{aligned}
\right.
\end{align}
One sees that $\Phi$ is a $C^2$-diffemorphism with $\Phi^{-1}=\Psi$.
We shall say that $\Phi$ straightens out the boundary at $z_0$. 
Below, by translation and dilation, without loss of generality, we assume $z_0=0$ so that
\begin{align}\label{Phi}
\Phi(z_0)=0\ \mbox{ and }\ \Phi(W)=B_1,\ \ \Phi(W\cap D)=B^+_1:=\{y\in B_1: y_1>0\}.
\end{align}
For $\delta\in(0,1]$, define
\begin{align}\label{HJ8}
\Gamma_\delta:=\p D\cap \Phi^{-1}(B_{\delta/4}),\ U_\delta:=\Phi^{-1}(B^+_\delta),\ \ V_\delta:=\Phi^{-1}(B^+_{\delta/2}).
\end{align}
From the above construction, it is easy to see that
\begin{align}\label{DG1}
\vec{n}(z_0)=-\nabla \Phi_1(z_0)/|\nabla \Phi_1(z_0)|,
\end{align}
where $\vec{n}(z_0)$ is the unit outward normal at point $z_0\in\p D$, and
\begin{align}\label{Equiv}
d_{D^c}(x)\leq \Phi_1(x)\leq c_0d_{D^c}(x),\ \ x\in U_1.
\end{align}
See Figure \ref{fig2} for flattening out the boundary.

The following lemma shows that if the vector field $b$ along the boundary is towards the interior of $D$, then 
the exit position $X_{\tau_D}$ would not touch the boundary.
\bl\label{Le71}
Let $z_0\in\p D$. If $b(z_0)\cdot \vec{n}(z_0)<0$, then there is a neighborhood $\Gamma\subset\p D$ of $z_0$ such that
for each $x\in D$,
$$
\mP_x(X_{\tau_D}\in\Gamma)=0.
$$
\el
\begin{proof}
Fix $z_0\in \p D$. Let $U_1$ and $\Phi$ be as above.
Since $x\mapsto b(x)\cdot\nabla \Phi_1(x)$ is continuous in $\bar U_1$, 
by \eqref{DG1} and $b(z_0)\cdot \vec{n}(z_0)<0$, without loss of generality, we may assume
\begin{align}\label{GF0}
c_0:=\inf_{x\in U_1}b(x)\cdot\nabla \Phi_1(x)>0.
\end{align}
(i) Let $\Gamma_\delta, U_\delta$ and $V_\delta$ be as in \eqref{HJ8} (see Figure \ref{fig2}). We first prove that for $\delta$ small enough,
\begin{align}\label{FD2}
\mP_x(X_{\tau_{U_\delta}}\in\Gamma_\delta)=0,\ \ x\in V_\delta.
\end{align}
Fix $\gamma\in(0,1)$. For any $n\in\mN$, let $\ell_n:\mR\to[0,\infty)$ be given by
$$
\ell_n(s):=
\left\{
\begin{array}{ll}
(n/\delta)^{\gamma+1}s,& s\in[0,\delta/n],\\
s^{-\gamma},& s\in[\delta/n,\delta],\\
\tfrac{(1-s)\delta^{-\gamma}}{1-\delta},& s\in[\delta,1],\\
0,&s\notin[0,1].
\end{array}
\right.
$$
Define
$$
f_n(x):=\ell_n(\Phi_1(x)),\ x\in U_1,\ f_n(x):=0,\ x\notin U_1
$$
and
$$
U^n_\delta:=\Phi^{-1}(\{y\in B^+_\delta, y_1>\delta/n\}).
$$
By definition and the change of variables, we have for $x=\Phi^{-1}(y)\in U^n_\delta$,
\begin{align*}
\sL^{(\alpha)}_\kappa f_n(x)&=\int_{U_1}(f_n(z)-f_n(x))\frac{\kappa(x,z-x)}{|z-x|^{d+\alpha}}\dif z
-\int_{U^c_1}f_n(x)\frac{\kappa(x,z-x)}{|z-x|^{d+\alpha}}\dif z\no\\
&\leq\|\kappa\|_\infty\int_{U_1}\frac{(f_n(z)-f_n(x))_+}{|z-x|^{d+\alpha}}\dif z
=\|\kappa\|_\infty\int_{B^+_1}\frac{(\ell_n(z_1)-y_1^{-\gamma})_+}{|\Phi^{-1}(z)-x|^{d+\alpha}}\det(\nabla\Phi^{-1}(z))\dif z\\
&\leq c\int_{\{0<z_1<y_1\}}\frac{z_1^{-\gamma}-y_1^{-\gamma}}{|z-y|^{d+\alpha}}\dif z
=cy_1^{-\gamma-\alpha}\int^{\bf 1}_0\frac{s^{-\gamma}-1}{(1-s)^{1+\alpha}}\dif s\int_{\mR^{d-1}}\frac{\dif w}{(1+|w|^2)^{(d+\alpha)/2}},
\end{align*}
where the constant $c$ does not depend on $n$ and $\delta$.
On the other hand, for $x\in U^n_\delta$, by \eqref{GF0},
\begin{align*}
b(x)\cdot\nabla f_n(x)=-\gamma \Phi_1(x)^{-\gamma-1} b(x)\cdot\nabla\Phi_1(x)\leq -\gamma c_0 \Phi_1(x)^{-\gamma-1}.
\end{align*}
Therefore, since $\alpha\in(0,1)$, if we let $\delta$ be small enough, then for all $x\in U^n_\delta$,
$$
\sL^{(\alpha)}_\kappa f_n(x)+b(x)\cdot\nabla f_n(x)\leq c\Phi_1(x)^{-\gamma-\alpha}-\gamma c_0\Phi_1(x)^{-\gamma-1}\leq 0.
$$
Thus, for $m>n$, since $f_m\in C_b(\mR^d)\cap C^2_{loc}(U^n_\delta)$,  we have for $x\in U^n_\delta$,
$$
\mP_x\Big(X_{\tau_{U^n_\delta}}\in U^m_\delta\setminus U^n_\delta \Big)(n/\delta)^\gamma\leq\mE_x f_m\Big(X_{\tau_{U^n_\delta}}\Big)\stackrel{\eqref{Dyn}}{\leq}f_m(x).
$$
Letting $m\to\infty$, we get for any $x\in U^n_\delta$,
\begin{align}\label{FD1}
\mP_x\Big(X_{\tau_{U^n_\delta}}\in (U_\delta\setminus U^n_\delta)\cup\Gamma_\delta\Big)=
\mP_x\Big(X_{\tau_{U^n_\delta}}\in U_\delta\setminus U^n_\delta \Big)\leq \Phi_1(x)^{-\gamma}(\delta/n)^\gamma,
\end{align}
where the first equality is due to Lemma \ref{Le52}.
Since $X_{\tau_U{^n_\delta}}\to X_{\tau_{U_\delta}}$ by Lemma \ref{Le51}, we have
$$
\Big\{X_{\tau_{U_\delta}}\in\Gamma_\delta\Big\}\subset\liminf_{n\to\infty}\Big\{X_{\tau_{U^n_\delta}}\in (U_\delta\setminus U^n_\delta)\cup\Gamma_\delta\Big\},
$$
which together with \eqref{FD1}, yields \eqref{FD2}.
\\
\\
(ii) Next we show that for any $x\in D$,
$$
h(x):=\mP_x(X_{\tau_D}\in\Gamma_\delta)=0.
$$
Notice that by the L\'evy system \eqref{Levy}, for $x\in D\setminus U_\delta$,
\begin{align}\label{HJ1}
\mP_x\left(X_{\tau_D\wedge\sigma_{V_\delta}}\in\Gamma_\delta\right)
=\mE_x\int^{\tau_D\wedge\sigma_{V_\delta}}_0\int_{\Gamma_\delta}\frac{\kappa(X_s,z-X_s)}{|z-X_s|^{d+\alpha}}\dif z=0,
\end{align}
where we have used that  $|z-X_s|>\delta$ for $s<\tau_D\wedge\sigma_{V_\delta}$ and $z\in\Gamma_\delta$.
By the strong Markov property,  we have for $x\in D\setminus U_\delta$, 
\begin{align*}
h(x)&=\mP_x(X_{\tau_D}\in\Gamma_\delta;\sigma_{V_\delta}<\tau_D)
+\mP_x(X_{\tau_D}\in\Gamma_\delta;\sigma_{V_\delta}\geq \tau_D)\\
&=\mE_x\Big(\mP_{X_{\sigma_{V_\delta}}}(X_{\tau_D}\in\Gamma_\delta);\sigma_{V_\delta}<\tau_D\Big)
+\mP_x\Big(X_{\tau_D\wedge\sigma_{V_\delta}}\in\Gamma_\delta\Big),
\end{align*}
and for $x\in V_\delta$,
\begin{align*}
h(x)&=\mP_x(X_{\tau_D}\in\Gamma_\delta;\tau_{U_\delta}<\tau_D)
+\mP_x(X_{\tau_D}\in\Gamma_\delta;\tau_{U_\delta}\geq \tau_D)\\
&=\mE_x\Big(\mP_{X_{\tau_{U_\delta}}}(X_{\tau_D}\in\Gamma_\delta);\tau_{U_\delta}<\tau_D\Big)
+\mP_x\Big(X_{\tau_{U_\delta}}\in\Gamma_\delta\Big).
\end{align*}
Thus, by \eqref{HJ1} and \eqref{FD2}, we get
\begin{align}
\sup_{x\in D\setminus U_\delta}h(x)&\leq\sup_{x\in V_\delta}h(x)\cdot \sup_{x\in D\setminus U_\delta}\mP_x(\sigma_{V_\delta}<\tau_D),\label{HJ2}\\
\sup_{x\in V_\delta}h(x)&\leq\sup_{x\in D\setminus U_\delta}h(x)\cdot \sup_{x\in V_\delta}\mP_x(\tau_{U_\delta}<\tau_D).\label{HJ3}
\end{align}
For $x\in D\setminus U_\delta$ and $\eps<{\rm dist}(D\setminus U_\delta, V_\delta)/2$, notice that
if $|\Delta X_{\tau_{B_\eps(x)}}|>2\rm {diam}(D)$, 
then $\sigma_{V_\delta}\geq \tau_D$. This means that the process has jumped out from $D$ before it enters into $V_\delta$. Hence,
$$
{\bf 1}_{\{\sigma_{V_\delta}\geq \tau_D\}}\geq {\bf 1}_{\{|\Delta X_{\tau_{B_\eps(x)}}|>2\lambda_D\}}
=\sum_{0<s\leq \tau_{B_\eps(x)}}{\bf 1}_{\{|\Delta X_s|>2\lambda_D\}}.
$$
By the L\'evy system \eqref{Levy}, we have
$$
\mP_x(\sigma_{V_\delta}\geq \tau_D)\geq\mE_x\int^{\tau_{B_\eps(x)}}_0\!\!\!\int_{|z|>2\lambda_D}\frac{\kappa(X_s, z)}{|z|^{d+\alpha}}\dif z\dif s\geq c_0\mE_x\tau_{B_\eps(x)},
$$
which together with \eqref{EK0} yields
$$
\inf_{x\in D\setminus U_\delta}\mP_x(\sigma_{V_\delta}\geq \tau_D)\geq\inf_{x\in D\setminus U_\delta}\mE_x\tau_{B_\eps(x)}>0.
$$
That is, $\sup_{x\in D\setminus U_\delta}\mP_x(\sigma_{V_\delta}<\tau_D)<1$. So, by \eqref{HJ2} and \eqref{HJ3},
$$
\sup_{x\in D\setminus U_\delta}h(x)=0,\ \ \sup_{x\in V_\delta}h(x)=0.
$$
The proof is thus complete.
\end{proof}
In the next lemma, we consider the following viscosity approximation operator
$$
\sA_\nu u:=\nu\Delta^{1/2} u+\sL^{(\alpha)}_\kappa u+b\cdot\nabla u,\ \ \nu\in[0,1].
$$
It is well known that the martingale problem associated with $\sA_\nu$ is well-posed (see \cite{Ch-Zh2}). 
The associated Markov process is denoted by $(X,\mP^\nu_x)$ and the expectation with respect to $\mP^\nu_x$ is denoted by $\mE^\nu_x$ .
\bl\label{Le72}
Let $z_0\in\p D$ and $\delta_0>0$.
Suppose that one of the following two conditions holds
$$
(i)\ b(z_0)\cdot \vec {n}(z_0)>0;\quad (ii)\  b(z)\cdot \vec {n}(z)=0\mbox{ for each $z\in\p D\cap B_{\delta_0}(z_0)$}.
$$
Then there are $\theta\in(0,\tfrac{\alpha}{2})$ and $\delta>0$ such that
\begin{align}\label{WP2}
\sup_{x\in D\cap B_\delta(z_0)}\sup_{\nu\in(0,1)}\left(d_{D^c}^{-\theta}(x)\mE^\nu_x\tau_D\right)<\infty.
\end{align}
Moreover, in the case (i), for $\nu=0$, we further have
\begin{align}\label{ES8}
\sup_{x\in D\cap B_\delta(z_0)}\left(d_{D^c}^{-1}(x)\mE^0_x\tau_D\right)<\infty.
\end{align}
\el
\begin{proof}
Let $U_\delta:=D\cap B_\delta(z_0) $. By Lemma \ref{Le66}, if we choose $\delta$ small enough, then
there are $\theta\in(0,\alpha/2)$ and $c_0>0$ such that for all $\nu\in(0,1)$,
\begin{align}\label{AS1}
\left(\nu\Delta^{1/2}d^\theta_{D^c}+\sL^{(\alpha)}_\kappa d^\theta_{D^c}\right)(x)\leq-c_0\nu d^{\theta-1}_{D^c}(x)-c_0  d^{\theta-\alpha}_{D^c}(x),\ \ x\in U_\delta.
\end{align}
(i) In the first case, since $x\mapsto b(x)\cdot\nabla d_{D^c}(x)$ is continuous in $\bar U_\delta$, 
by \eqref{DG1}, without loss of generality, we may assume
$$
c_1:=-\sup_{x\in U_\delta}b(x)\cdot\nabla d_{D^c}(x)>0,
$$
which implies that
\begin{align}\label{AS2}
b(x)\cdot\nabla d^\theta_{D^c}(x)\leq -\theta c_1d^{\theta-1}_{D^c}(x),\ \ x\in U_\delta.
\end{align}
(ii)  In the second case, for $x\in U_\delta$, let $z\in\p D\cap B_\delta(z_0)$ be such that $d_{D^c}(x)=|x-z|$. Since
$b(z)\cdot\nabla d_{D^c}(z)=0$, by the H\"older continuity of $b$, we have
$$
|b(x)\cdot\nabla d_{D^c}(x)|=|b(x)\cdot\nabla d_{D^c}(x)-b(z)\cdot\nabla d_{D^c}(z)|\lesssim |x-z|^\beta=d^\beta_{D^c}(x),
$$
and so,
\begin{align}\label{AS3}
|b(x)\cdot\nabla d^\theta_{D^c}(x)|=\theta d^{\theta-1}_{D^c}(x) |b(x)\cdot\nabla d_{D^c}(x)|\leq \theta c_1 d_{D^c}^{\theta+\beta-1}(x).
\end{align}
Combining \eqref{AS1}, \eqref{AS2} and \eqref{AS3}, by choosing $\delta$ small enough, we always have
$$
\left(\sA_\nu d^\theta_{D^c}\right)(x)\leq -c_2  d^{\theta-\alpha}_{D^c}(x)\leq -c_2\delta^{\theta-\alpha},  \ x\in U_\delta,
$$
where $c_2$ is independent of $\nu$ and $x$.
Thus by \eqref{Dyn}, we have for all $x\in U_\delta$,
$$
\mE^\nu_x \left(d^\theta_{D^c}(X_{\tau_{U_\delta}})\right)\leq d^\theta_{D^c}(x)-c_2\delta^{\theta-\alpha}\mE^\nu_x\tau_{U_\delta}.
$$
Therefore, for all $x\in U_{\delta/2}$,
$$
\mE^\nu_x\tau_{U_\delta}\leq \delta^{\alpha-\theta}d^\theta_{D^c}(x)/c_2,
$$
which together with Lemma \ref{Le55} (choosing $U=B_\delta(z_0)$ and $V=B_{\delta/2}(z_0)$ there) yields the desired estimate.
\\
\\
(iii) Finally, in the first case, for $\delta>0$ define $f(x):=d_{D^c}(x)\wedge\delta$. We have
\begin{align*}
|\sL^{(\alpha)}_\kappa f(x)|&\lesssim\int_{|z|\leq\delta}|f(x+z)-f(x)|\frac{\dif z}{|z|^{d+\alpha}}+\int_{|z|>\delta}|f(x+z)-f(x)|\frac{\dif z}{|z|^{d+\alpha}}\\
&\leq\int_{|z|\leq \delta}|z|\dif z/|z|^{d+\alpha}+2\delta\int_{|z|>\delta}\dif z/|z|^{d+\alpha}\leq c\delta^{1-\alpha},
\end{align*}
where the constant $c$ is independent of $\delta$.
Hence, by \eqref{AS2}, for $\delta$ small enough,
$$
(\sL^{(\alpha)}_\kappa +b\cdot\nabla) f(x)\leq  c\delta^{1-\alpha}-c_1\leq -c_1/2,  \ x\in U_\delta.
$$
As above we get \eqref{ES8}.
\end{proof}

Let $\Gamma:=\{z\in\p D: b(z)\cdot \vec {n}(z)=0\}$.
In the case (ii) of the above lemma, it does not tell us that for the boundary point $z_0\in \Gamma$, whether it holds
\begin{align}\label{WP3}
\lim_{D\ni x\to z_0}\sup_{\nu\in(0,1)}\mE^\nu_x\tau_D=0.
\end{align}
Notice that the estimate \eqref{WP2} only implies that the above limit holds for the interior point $z_0$ of closed set $\Gamma$.
However, when $\alpha\geq 1-\frac{\beta}{2}$, we have the following affirmative answer.
\bl\label{Le73}
Let $\frac{1}{2}\leq\alpha<1$ and $z_0\in\p D$. Assume $b(z_0)\cdot \vec {n}(z_0)=0$ and $b\in C^\beta$ with $\beta\in[2(1-\alpha),1]$. Then limit \eqref{WP3} holds.
\el
\begin{proof}
For $\delta>0$, let $U_\delta:= D\cap B_\delta(z_0)$. Let $B=B_r(y), B'=B_{2r}(y')$ be two 
exterior tangent balls of $D$ at point $z_0$  so that $B\cap\bar D=B'\cap\bar D=\{z_0\}$.
Without loss of generality, we assume $z_0=0$ and $y=(r,0)$, $y'=(2r,0)$ (see Figure 3 below). 
For $\delta\leq r$ and $x=(x_1,x^*_1)\in U_\delta\subset (B')^c$, we have
$$
-\delta<x_1<2r-\sqrt{(2r)^2-|x_1^*|^2}\leq \frac{|x_1^*|^2}{3r}. 
$$
Hence,  
\begin{align*}
d_B(x)=&\sqrt{(r-x_1)^2+|x_1^*|^2}-r=r\left[\sqrt{1-\frac{2x_1}{r}+\frac{|x|^2}{r^2}}-1\right]\\
\geq & r\Big[\sqrt{1+|x|^2/(3r^2)}-1\Big]\geq |x|^2/ {(9r)},
\end{align*}
which means that
\be
\label{Dist} |x-z_0|^2\leq (9r)d_B(x),\  x\in U_\delta. 
\ee 
Now by Lemma \ref{Est-F1}, there are $\theta_0\in(0,\alpha/2)$ and $c_0>0$ such that for all $\theta\in(0,\theta_0)$, $\nu\in(0,1)$ and $x\in U_\delta$,
\begin{align}\label{AS91}
\left(\nu\Delta^{1/2}d^\theta_B+\sL^{(\alpha)}_\kappa d^\theta_B\right)(x)\leq-c_0\nu d^{\theta-1}_B(x)-c_0  d^{\theta-\alpha}_B(x).
\end{align}
Since $b(z_0)\cdot \nabla d_B(z_0)=-b(z_0)\cdot \vec {n}(z_0)=0$ and $b\in C^\beta$, there is a constant $c>0$ such that for all $x\in U_\delta$
(see Figure \ref{fig4}),
\begin{figure}[h]\small
\centering
\includegraphics
[height=0.3\textwidth,width=0.95\textwidth]{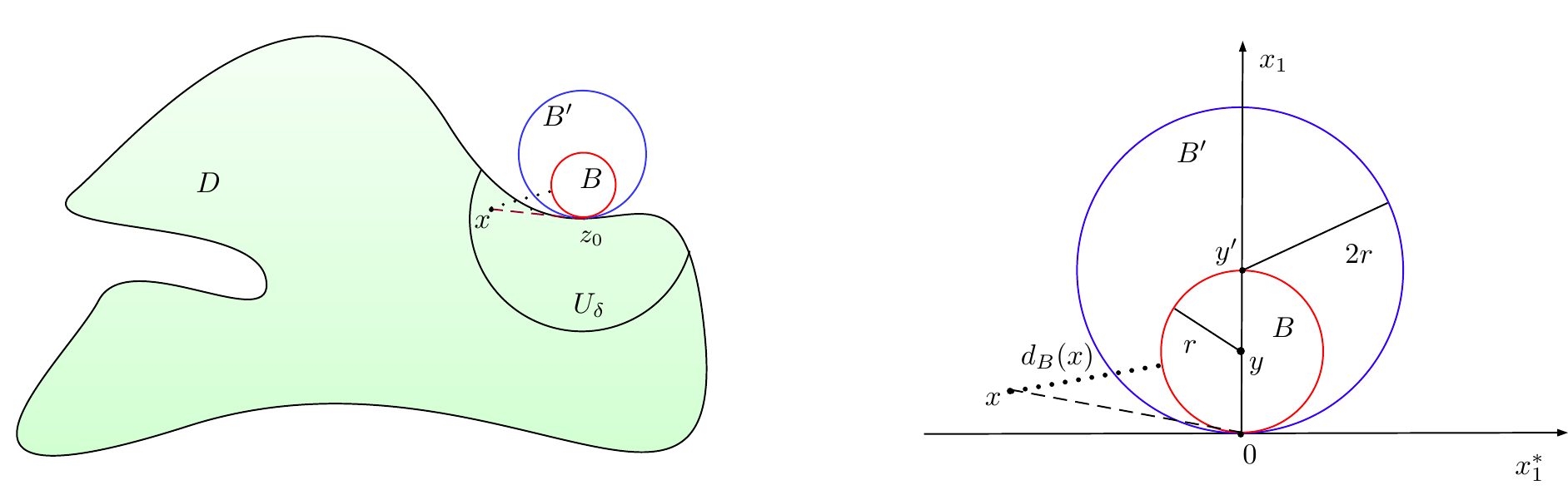}
\caption{Distance function to exterior tangent ball}\label{fig4}
\end{figure}
\begin{align}\label{AS92}
|b(x)\cdot \nabla d_B(x)|=|b(x)\cdot \nabla d_B(x)-b(z_0)\cdot \nabla d_B(z_0)|\leq c|x-z_0|^\beta\overset{\eqref{Dist}}{\leq} cd^{\beta/2}_B(x). 
\end{align}
Since $\beta/2+\alpha\geq 1$, combining \eqref{AS91} and \eqref{AS92}, by choosing $\delta,\theta$ small enough, we get
$$
\left(\sA_\nu d^\theta_B\right)(x)\leq -c_0d^{\theta-\alpha}_B(x)+c\theta d^{\theta-1+\beta/2}_B(x)\leq -c_0d^{\theta-\alpha}_B(x)/2\leq -c_1\delta^{\theta-\alpha},  \ x\in U_\delta.
$$
Hence, by \eqref{Dyn}, we have for all $x\in U_\delta$,
$$
\mE^\nu_x \left(d^\theta_B(X_{\tau_{U_\delta}})\right)\leq d^\theta_B(x)-c_1\delta^{\theta-\alpha}\mE^\nu_x\tau_{U_\delta},
$$
which implies that
$$
\lim_{U_\delta\ni 
x\to z_0}\sup_{\nu\in(0,1)}\mE^\nu_x\tau_{U_\delta}\leq c^{-1}_1\delta^{\alpha-\theta}\lim_{U_\delta\ni 
x\to z_0} d^\theta_B(x)=0.
$$
The proof is complete by Lemma \ref{Le55}.
\end{proof}
\subsection{Proof of Theorem \ref{MAIN1}}
Let $\alpha\in(0,1)$ and $\nu>0$.
Consider the following nonlocal supercritical Dirichlet problem with viscosity term $\nu\Delta^{1/2}$:
\begin{align}\label{IBP1}
\left\{
\begin{aligned}
&\p_t u=\nu\Delta^{1/2} u+\sL^{(\alpha)}_{\kappa} u+b\cdot\nabla u+f\ \mbox{ on }\ \mR_+\times D,\\
&u=0\ \mbox{ on }\ \mR_+\times D^c,\ \ u(0,x)=\varphi(x),\ x\in D.
\end{aligned}
\right.
\end{align}
We first show the unique solvability to the above Dirichlet problem. 
\bt
Let $\alpha,\beta\in(0,1)$ with $\alpha+\beta>1$ and $\gamma\in(1-\alpha,\beta]$. 
Suppose \eqref{Con1}, $b\in C^\beta$ and
$$
|\kappa(x,z)-\kappa(x,z')|\leq \kappa_1|z-z'|^\gamma.
$$
For any $\theta\in(0,\frac{1}{2})$, $f\in \mB^{(0)}_\gamma(D)$ and $\varphi\in\cC^{(-\theta)}_{1+\gamma}(D)$, 
there is a unique solution $u_\nu\in \mB^{(-\theta)}_{1+\gamma}(D)$
solving problem \eqref{IBP1}. Moreover, there are $\theta_0=\theta_0(\alpha,\beta,\gamma,\kappa_0,d)>0$ and $c>0$ 
only depending on $\alpha,\beta,\gamma,\kappa_0,d, \kappa_1,\|b\|_{C^\beta}$
such that for all $\nu>0$ and $T>0$,
\begin{align}\label{Uni0}
\|u_\nu\|_{\mB^{(\theta_0)}_{\alpha+\gamma;T}(D)}&\leq 
c\left((1+T)\|f\|_{\mB^{(0)}_{\gamma;T}(D))}
+\|\varphi\|_{\cC^{(0)}_{\alpha+\gamma}(D)}\right).
\end{align}
\et
\begin{proof}
As for the unique solvability of equation \eqref{IBP1}, 
it follows by \eqref{INSS} and the continuity method as used in proving Theorem \ref{MAIN}. 
We only show the uniform estimate \eqref{Uni0}.
By Theorem \ref{Main44}, there are $\theta_0>0$ and $c>0$ such that for all $\nu\in(0,1)$ and $T>0$,
\begin{align}\label{Uni1} 
\|u\|_{\mB^{(\theta_0)}_{\alpha+\gamma;T}(D)}&\leq 
c\left(\|f\|_{\mB^{(\alpha+\theta_0)}_{\gamma;T}(D)}
+\|\varphi\|_{\cC^{(\theta_0)}_{\alpha+\gamma}(D)}+\|u\|_{\mB^{(0)}_{0;T}(D)}\right).
\end{align}
On the other hand, by the maximum principle \eqref{Max}, we have
$$
\|u\|_{\mB^{(0)}_{0;T}(D)}=\|u\|_{L^\infty_T(C^{0}(D))}\leq \|\varphi\|_{C^0(D)}+T\|f\|_{L^\infty_T(C^0(D))}.
$$
Substituting this into \eqref{Uni1}, we obtain \eqref{Uni0}.
\end{proof}
\begin{proof}[Proof of Theorems \ref{MAIN1} and \ref{MAIN2}]
Let $\varphi\in \cC^{(0)}_{\alpha+\gamma}(D)$. For $\theta\in(0,\frac{1}{2})$ and $\nu\in(0,1)$, 
let $\rho$ be a nonnegative smooth function with support in $B_1$ and
$\int_{\mR^d}\rho=1$. Define
$$
\rho_\nu(x):=\nu^{-d}\rho(x/\nu),\ \ \chi_\nu:={\bf 1}_{D_{3\nu}}*\rho_\nu,\ \varphi_\nu:=(\varphi*\rho_\nu)\cdot\chi_\nu, 
$$
where $D_{3\nu}:=\{x\in D: \mathrm{dist}(x, D^c)>3\nu\}$. Then $\varphi_\nu\in\cC^{(-\theta)}_{1+\gamma}(D)$ and 
for some $c>0$,
\begin{align}\label{AS5}
\sup_{\nu\in(0,1)}\|\varphi_\nu\|^{(0)}_{\alpha+\gamma; D}\leq c\|\varphi\|^{(0)}_{\alpha+\gamma;D},
\end{align}
and for each $x\in D$,
$$
\varphi_\nu(x)\to\varphi(x),\ \nu\to 0.
$$
Let $u_\nu\in \mB^{(-\theta)}_{1+\gamma}(D)$
be the unique solution of \eqref{IBP1} corresponding to  ($f$, $\varphi_\nu$).
By \eqref{Uni0} and \eqref{AS5}, we have the following uniform estimate:
\begin{align}\label{AL1}
\sup_{\nu\in(0,1)}\|u_\nu\|_{\mB^{(\theta_0)}_{\alpha+\gamma;T}(D)}&\leq 
c\left((1+T)\|f\|_{\mB^{(0)}_{\gamma;T}(D)}+\|\varphi\|_{\cC^{(0)}_{\alpha+\gamma}(D)}\right),
\end{align}
and also by the maximum principle \eqref{Max},
\begin{align}\label{AL3}
\sup_{\nu\in(0,1)}\|u_\nu\|_{\mB^{(0)}_{0;T}(D)}=\sup_{\nu\in(0,1)}\|u_\nu\|_{L^\infty_T(C^{0}(D))}\leq T\|f\|_{L^\infty_T(C^0(D))}+\|\varphi\|_{C^0(D)}.
\end{align}
On the other hand, by Theorem \ref{Th0} with $\theta=-\theta_0$, we also have
$$
\|\Delta^{1/2}u_\nu\|^{(1+\theta_0)}_{\alpha+\gamma-1; D}\lesssim \|u_\nu\|^{(\theta_0)}_{\alpha+\gamma; D}+\|u_\nu\|_{0; D},
$$
$$
\|\sL^{(\alpha)}_\kappa u_\nu\|^{(1+\theta_0)}_{\alpha+\gamma-1; D}
\lesssim \|u_\nu\|^{(1-\alpha+\theta_0)}_{2\alpha+\gamma-1; D}+\|u_\nu\|_{0; D}
\lesssim \|u_\nu\|^{(\theta_0)}_{\alpha+\gamma; D}+\|u_\nu\|_{0; D},
$$
where the constant contained in $\lesssim$ is independent of $\nu$,
and by \eqref{HL6},
$$
\|b\cdot\nabla u_\nu\|^{(1+\theta_0)}_{\alpha+\gamma-1; D}\lesssim\|b\|^{(0)}_{\alpha+\gamma-1; D}
\|\nabla u_\nu\|^{(1+\theta_0)}_{\alpha+\gamma-1; D}
\lesssim \|b\|_{\beta; D}\|u_\nu\|^{(\theta_0)}_{\alpha+\gamma; D}.
$$
Hence, by equation \eqref{IBP1}, 
\begin{align}\label{AL2}
\|\p_t u_\nu\|_{\mB^{(1+\theta_0)}_{\alpha+\gamma-1;T}(D)}\lesssim \|u_\nu\|_{\mB^{(\theta_0)}_{\alpha+\gamma;T}(D)}
+\|u_\nu\|_{\mB^{(0)}_{0;T}(D)}+\|f\|_{\mB^{(1+\theta_0)}_{\alpha+\gamma-1;T}(D)}
\leq c\left(\|f\|_{\mB^{(0)}_{\gamma;T}(D)}+\|\varphi\|_{\cC^{(0)}_{\alpha+\gamma}(D)}\right).
\end{align}
Thus by \eqref{AS5}-\eqref{AL2} and Lemma \ref{Le42}, there is a subsequence $\nu_k\to 0$ and  
$$
u\in \mB^{(\theta_0)}_{\alpha+\gamma}(D)\cap \mB^{(0)}_0(D)
$$ 
such that for any $\eps\in(0,\alpha+\gamma-1)$, $T>0$ and $D_0\Subset D$,
\begin{align}\label{YF1}
\lim_{n\to\infty}\|u_{\nu_k}-u\|_{L^\infty_T(C^{\alpha+\gamma-\eps}(D_0))}=0.
\end{align}
Let $\sL^{(\alpha)}_{\kappa,b}:=\sL^{(\alpha)}_{\kappa}+b\cdot\nabla$. Since for any test function $\phi\in C^\infty_c(D)$ and $t\geq 0$,
$$
\int_{\mR^d} u_{\nu_k}(t)\phi=\int_{\mR^d}\varphi \phi+\nu_{\nu_k}\int^t_0\!\!\int_{\mR^d}  u_{\nu_k}\Delta^{1/2} \phi
+\int^t_0\!\!\int_{\mR^d} \phi\sL^{(\alpha)}_{\kappa,b} u_{\nu_k}+\int^t_0\!\!\int_{\mR^d} f\phi,
$$
by \eqref{YF1} and taking limits $k\to\infty$, we obtain
$$
\int_{\mR^d} u(t)\phi=\int_{\mR^d}\varphi \phi+\int^t_0\!\!\int_{\mR^d} \phi\sL^{(\alpha)}_{\kappa,b} u+\int^t_0\!\!\int_{\mR^d} \phi f,
$$
which implies that
$$
u(t,x)=\varphi(x)+\int^t_0\sL^{(\alpha)}_{\kappa,b} u(s,x)\dif s+\int^t_0f(s,x)\dif s,\ (t,x)\in\mR_+\times D.
$$
In particular, $(u,f)\in \bH^{\alpha+\gamma}(D)$ (see \eqref{AL5} for a definition of $\bH^{\alpha+\gamma}(D)$).
\\
\\
{\it Case {\bf (A)} of Theorem \ref{MAIN1}}: Since $\mP_x(X_{\tau_D}\in\p D)=0$ by Lemma \ref{Le71}, the condition (i) in Theorem \ref{Th44} is satisfied.
Thus the probabilistic representation holds and the uniqueness follows from it.
\\
\\
{\it Cases {\bf (B)} and {\bf (C)} of Theorem \ref{MAIN1}}: One can show that the condition (ii) in Theorem \ref{Th44} is satisfied for $u$ and $\Gamma_1=\p D$, that is,
\begin{align}\label{HS85}
u\in C((0,\infty)\times\bar D),\ \ u|_{(0,\infty)\times\p D}=0.
\end{align}
Indeed, by Remark \ref{Re13} and Proposition \ref{Pr512}, one sees that the condition (i) in Theorem \ref{Th44} is satisfied for $\mP^{\nu_k}_x$. So,
we have the following probabilistic representation for $u_{\nu_k}$,
$$
u_{\nu_k}(t,x)=\mE^{\nu_k}_x \Big(\varphi_{\nu_k}(X_{t}){\bf 1}_{\{\tau_{D}>t\}}\Big)+\mE^{\nu_k}_x\left(\int^{t\wedge\tau_{D}}_0f(t-s,X_s)\dif s\right).
$$
From this, one sees that
\begin{align}\label{WP5}
|u_{\nu_k}(t,x)|&\leq\|\varphi\|_\infty\mP^{\nu_k}_x(\tau_{D}>t)+\|f\|_\infty\mE^{\nu_k}_x\tau_D\leq 
\left(\|\varphi\|_\infty /t+\|f\|_\infty\right)\mE^{\nu_k}_x\tau_D
\end{align}
Since $\p D$ is compact, 
by Lemma \ref{Le72} and a standard covering technique, there are $\theta\in(0,\frac{\alpha}{2})$ and $c>0$ such that for all $x\in D$ and $\nu_k\in(0,1)$,
\begin{align*}
|u_{\nu_k}(t,x)|\leq c\left(\|\varphi\|_\infty  /t+\|f\|_\infty \right)d_{D^c}^\theta(x).
\end{align*}
By taking limit $\nu_k\to 0$, we obtain
$$
|u(t,x)|\leq c\left(\|\varphi\|_\infty /t+\|f\|_\infty \right)d_{D^c}^\theta(x) ,
$$
which implies \eqref{HS85}. Thus, the  probabilistic representation \eqref{EP0} holds and the uniqueness follows immediately.
In case {\bf (C)}, by \eqref{EP0} and \eqref{ES8}, as above we have 
$$
|u(t,x)|\leq c\left(\|\varphi\|_\infty /t+\|f\|_\infty \right)d_{D^c}(x).
$$
{\it Case of Theorem \ref{MAIN2}}: If one can show that the condition (ii) in Theorem \ref{Th44} is satisfied for $\Gamma_0=\Gamma_<$
and $\Gamma_1=\Gamma_>\cup\Gamma_=$, then the probabilistic representation holds and the uniqueness follows from it.
First of all, by Lemma \ref{Le71} we have
$$
\mP_x(X_{\tau_D}\in\Gamma_0)=\mP_x(X_{\tau_D}\in\Gamma_<)=0.
$$
Next we need to check that
\begin{align}\label{HS5}
u\in C((0,\infty)\times(D\cup\Gamma_1)),\ \ u|_{(0,\infty)\times\Gamma_1}=0.
\end{align}
By \eqref{WP5}, Lemmas \ref{Le72} and \ref{Le73}, we obtain that for any $0<t_0<t_1<\infty$ and $z_0\in \Gamma_1$,
$$
\lim_{D\ni x\to z_0}\sup_{t\in[t_0,t_1]}|u(t,x)|\leq c\left(\|\varphi\|_\infty /t_0+\|f\|_\infty \right)\lim_{D\ni x\to z_0}\sup_{k}\mE^{\nu_k}_x\tau_D=0,
$$
which implies \eqref{HS5}. Finally, the conclusions (i) and (ii) follows by \eqref{WP5}, \eqref{WP2} and \eqref{ES8}.
The proof is complete.
\end{proof}

\section{Appendix}
\subsection{One dimensional L\'evy processes with a drift}

Fix $\alpha\in(0,1)$ and let $Z_t$ be a one dimensional  rotationally invariant and symmetric $\alpha$-stable process over some probability space $(\tilde\Omega,\tilde\sF,\tilde\mP)$.
For $x\in\mR$, let $\mP_x$ be the law of L\'evy process $X^x_t=x+Z_t+t$ in canonical space $\Omega$.
The following proposition shows that when $\alpha\in(0,1)$, the behavior of a L\'evy process with a drift could be quite different as in the case of $\alpha\in[1,2)$.
\bp\label{Prop81}
Let $D:=(0,1)$. It holds that
$$
(i)\ \inf_{x\in(0,1/4)}\mE_x\tau_D>0;\ \ (ii)\ \sup_{x\in D}\mP_x(X_{\tau_D}=0)=0;\ \ (iii)\sup_{x\in D}\mP_x(X_{\tau_D}=1)>0.
$$
\ep
\begin{proof}
(i) First of all, it is well known that for any $\eps>0$ (see \cite{Be}), 
$$
\lim_{t\to 0+}|Z_t|/ t^{\frac{1}{\alpha}-\eps}=0,\ \ \tilde \mP-a.s.
$$
Thus, since $\alpha\in(0,1)$, one can choose $\eps\in(0,\frac{1}{\alpha}-1)$ and $\delta\in(0,\frac{1}{8})$ small enough such that
$$
\tilde\mP(|Z_t|\leq t^{\frac{1}{\alpha}-\eps},\ \forall t\leq\delta)\geq \tfrac{1}{2}.
$$
For any $x\in(0,\frac{1}{4})$, letting $Q_x:=(\frac{x}{2},\frac{1}{2})$, we have
\begin{align}\label{Del}
\mP_x(\tau_{Q_x}>\delta)&=\tilde\mP\Big(X^x_t\in(\tfrac{x}{2},\tfrac{1}{2}),\forall t\leq\delta\Big)\geq
\tilde\mP\Big(|Z_t|\leq t^{\frac{1}{\alpha}-\eps},\ \forall t\leq\delta\Big)\geq \tfrac{1}{2},
\end{align}
and
\begin{align}\label{ER1}
\mE_x\tau_D\geq\mE_x\tau_{Q_x}\geq\delta\mP_x(\tau_{Q_x}>\delta)\geq\tfrac{\delta}{2}.
\end{align}
(ii) Let $D_0:=(0,\infty)$. Since $\tau_D\leq\tau_{D_0}$,  we prove the following stronger claim: 
\begin{align}\label{ET2}
\mP_x(X_{\tau_{D_0}}=0)=0,\ \ \forall x\in D_0.
\end{align}
For any $x\in(0,\frac{1}{4})$, letting $Q_x:=(\frac{x}{2},\frac{1}{2})$, by formula \eqref{Le} and \eqref{ER1}, we have
$$
\mP_x(X_{\tau_{Q_x}}\in D^c_0)=\mE_x\int^{\tau_{Q_x}}_0\dif s\int_{D^c_0}\frac{\dif z}{|z-X_s|^{1+\alpha}}
\geq\mE_x\tau_{Q_x}\int_{z<-1/2} \frac{\dif z}{|z|^{1+\alpha}}\geq c_1.
$$
For $x\geq \frac{1}{4}$, define $Q_x:=(\frac{x}{2},\frac{3x}{2})$. Noticing that for any $\lambda>0$, $(Z_{\lambda t})_{t\geq 0}$ has the same law as
$(\lambda^{1/\alpha}Z_t)_{t\geq 0}$, we have for $\lambda$ small enough independent of $x\geq \frac{1}{4}$,
\begin{align*}
\mP_x\Big(\tau_{Q_x}<(\lambda x)^\alpha\Big)\leq\tilde\mP\left(\sup_{t\in[0,(\lambda x)^\alpha]}|Z_t+t|\geq\frac{x}{2}\right)
\leq\tilde\mP\left(\sup_{t\in[0,(\lambda x)^\alpha]}|Z_t|\geq\frac{x}{4}\right)=
\tilde\mP\left(\sup_{t\in[0,1]}|Z_t|\geq\frac{1}{4\lambda}\right)<1.
\end{align*}
Hence, for $x\geq\frac{1}{4}$ and $Q_x:=(\frac{x}{2},\frac{3x}{2})$, by the L\'evy system again, we have
\begin{align*}
\mP_x(X_{\tau_{Q_x}}\in D^c_0)&=\mE_x\int^{\tau_{Q_x}}_0\dif s\int_{D^c_0}\frac{\dif z}{|z-X_s|^{1+\alpha}}
\geq\mE_x\tau_{Q_x}\int_{z<-3x/2}\frac{\dif z}{|z|^{1+\alpha}}\\
&\geq  (\lambda x)^\alpha\mP_x\Big(\tau_{Q_x}\geq(\lambda x)^\alpha\Big)\cdot cx^{-\alpha}\geq c_2.
\end{align*}
Thus, by Lemma \ref{Le59} with $c_0=c_1\wedge c_2$, we get \eqref{ET2}.
\\
\\
(iii) Noticing that
$$
\Psi(\xi):=-\log\tilde\mP(\e^{{\rm i}\xi (Z_1+1)})=|\xi|^\alpha-{\rm i}\xi,
$$
we have
$$
\int^\infty_0{\rm Re}\left(\frac{1}{1+\Psi(\xi)}\right)\dif\xi=\int^\infty_0{\rm Re}\left(\frac{1}{1+|\xi|^\alpha-{\rm i}\xi}\right)\dif\xi
=\int^\infty_0\frac{1+|\xi|^\alpha}{(1+|\xi|^\alpha)^2+|\xi|^2}\dif\xi<\infty.
$$
By Kesten's theorem (see \cite{Ke}, \cite{Be}), one-point sets are non-polar sets of $Z_t+t$. To show (iii), we use a contradiction argument.
Fix $\eps\in(0,1)$ and let $U:=[\eps,1)$. We shall show that for some $x\in U$, $\mP_x(X_{\tau_U}=1)>0$, which automatically implies that
$\mP_x(X_{\tau_D}=1)>0$. Suppose now that
\begin{align}\label{FG5}
\mP_x(X_{\tau_U}=1)=0,\ \ \forall x\in U.
\end{align}
Under this assumption we show that the single point set $\{1\}$ is a polar set.
Let $\sigma_0:=0$.  For $n\in\mN$, define stopping times $\sigma_n$ and $\tau_n$ recursively as follows:
$$
\tau_n:=\inf\{s>\sigma_{n-1}: X_t\notin U\},\ \sigma_n:=\inf\{s>\tau_n: X_s\in U\}.
$$
By (ii) and the L\'evy system, we have
$$
\mP_x(X_{\sigma_1}=1)=0,\ \forall x\in U^c,
$$
and furthermore, for any $n\geq 2$,
\begin{align}\label{GH1}
\mP_x(X_{\tau_n}=1)=0,\ \ \forall x\in U,\ \ \mP_x(X_{\sigma_n}=1)=0,\ \forall x\in U^c.
\end{align}
Indeed, since $\tau_n=\sigma_{n-1}+\tau_U\circ\theta_{\sigma_{n-1}}$, 
by the strong Markov property and \eqref{FG5}, we have
\begin{align*}
\mP_x(X_{\tau_n}=1)=\mP_x(X_{\tau_U}\circ\theta_{\sigma_{n-1}}=1)=\mE_x\Big(\mP_{X_{\sigma_{n-1}}}(X_{\tau_U}=1)\Big)=0,
\end{align*}
and similarly,
\begin{align*}
\mP_x(X_{\sigma_n}=1)=\mP_x(X_{\sigma_1}\circ\theta_{\tau_{n}}=1)=\mE_x\Big(\mP_{X_{\tau_{n}}}(X_{\sigma_1}=1)\Big)=0.
\end{align*}
Let $\sF_0:=\{\emptyset,\Omega\}$ and $\sF_n$ the sigma-field generated by $\sigma_1,\cdots,\sigma_n$. Let $\delta$ be as in \eqref{Del} and define
$$
A_n:=\{\sigma_n-\sigma_{n-1}>\delta\}.
$$
Since $\tau_n-\sigma_{n-1}=\tau_U\circ\theta_{\sigma_{n-1}}$, by the strong Markov property again and \eqref{Del}, we have
\begin{align*}
\mP_x(A_n|\sF_{n-1})\geq\mP_x\Big(\tau_n-\sigma_{n-1}>\delta\big|\sF_{n-1}\Big)
=\mP_{X_{\sigma_{n-1}}}(\tau_U>\delta)\geq\tfrac{1}{2},\ \ a.s.
\end{align*}
Thus, by the second Borel-Cantelli's lemma (see \cite[Theorem 5.3.2]{Du}), we have
$$
\mP_x(A_n\ i.o.)=1\Rightarrow\mP_x\Big(\sup_n\sigma_n=\infty\Big)=1,
$$
which together with  \eqref{GH1} implies that 
$$
\mP_x(X_t=1,\exists t\geq 0)=\mP_x\Big(\cup_{n\in\mN}\{X_{\sigma_n}=1\}\cup\{X_{\tau_n}=1\}\Big)=0,\ \ \forall x\in U.
$$
In other words, the single point set $\{1\}$ is a polar set. Thus, we get a contradiction.
\end{proof}

\subsection{Solvability of fractional Dirichlet problems}
The aim of this subsection is to provide a self-contained proof for Theorem \ref{Th52}.
First of all we show the following interior estimate, which is essentially contained in Theorem \ref{Main0}. 
Here the main difference is that $u$ is not necessarily zero outside $D$.
For the reader's convenience, we prove it again.
\bl
Let $D_0\Subset D\Subset D_1$ be bounded domains and $\alpha\in(0,2)$, $\gamma\in(0,1)$. 
Assume that $u\in L^\infty_{loc}(\mR_+; C^{\alpha+\gamma}_{loc}(D_1)\cap L^\infty(\mR^d))$ satisfies
$$
\p_t u=\Delta^{\frac{\alpha}{2}}u+f\mbox{ in $\mR_+\times D$},\ \ u(0)=0.
$$
Then there is a constant $c=c(D_0,D,\alpha,\gamma,d)>0$ such that for all $T>0$,
\begin{align}\label{EY1}
\|u\|_{L^\infty_T(C^{\alpha+\gamma}(D_0))}\leq 
c\left(\|f\|_{L^\infty_T(C^\gamma(D))}+\|u\|_{L^\infty_T(L^\infty(\mR^d))}\right).
\end{align}
\el
\begin{proof}
Since $u\in L^\infty_{loc}(\mR_+; C^{\alpha+\gamma}_{loc}(D_1)\cap L^\infty(\mR^d))$ and $D\Subset D_1$, it is easy to see that
$$
\sup_{t\in[0,T]}\|u(t)\|^{(0)}_{\alpha+\gamma;D}<\infty.
$$
For $x_0\in D$, let $R:=d_{D^c}(x_0)/8$ and define
\begin{align}\label{EY3}
u_R(t,x):=R^{-\alpha}u(R^\alpha t, Rx+x_0),\ \ f_R(t, x):=f(R^\alpha t,Rx+x_0),\ w_R:=u_R\chi_1.
\end{align}
By definitions and scaling, it is easy to see that
$$
\p_t w_R=\Delta^{\frac{\alpha}{2}} w_R
+\chi_1\Delta^{\frac{\alpha}{2}}((1-\chi_3) u_R)+[\chi_1,\Delta^{\frac{\alpha}{2}}](\chi_3 u_R)+\chi_1 f_R\ \ \mbox{in} \ \ \mR_+\times \mR^d.
$$
Noticing that
\begin{align}\label{EY4}
u_R(R^{-\alpha}t,x)=R^{-\alpha} u(t,Rx+x_0)=:R^{-\alpha}u^{x_0}_R(t),
\end{align}
by the global Schauder's estimate \eqref{N8}, we have
\begin{align*}
R^{-\alpha}\|u^{x_0}_R\|_{L^\infty_T(C^{\alpha+\gamma}(B_1))}
\leq\|w_R\|_{\bB^{\alpha+\gamma}_{TR^{-\alpha}}}
&\lesssim \|\chi_1\Delta^{\frac{\alpha}{2}}((1-\chi_3) u_R)\|_{\bB^\gamma_{TR^{-\alpha}}}
+\|[\chi_1,\Delta^{\frac{\alpha}{2}}](\chi_3 u_R)\|_{\bB^\gamma_{TR^{-\alpha}}}\\
&+\|\chi_1 f_R\|_{\bB^\gamma_{TR^{-\alpha}}}+\|w_R\|_{\bB^0_{TR^{-\alpha}}}
=:I_1+I_2+I_3+I_4.
\end{align*}
For $I_1$, since $((1-\chi_3) u_R)(x)=R^{-\alpha}((1-\chi^{x_0}_{3R})u)(Rx+x_0)$, as in \eqref{PY4} we have
\begin{align*}
I_1&\lesssim \left\|\Delta^{\frac{\alpha}{2}}((1-\chi_3) u_R)\right\|_{L^\infty_{TR^{-\alpha}}(C^\gamma(B_2))}
=\left\|\Big(\Delta^{\frac{\alpha}{2}}((1-\chi^{x_0}_{3R}) u)\Big)^{x_0}_{R}\right\|_{L^\infty_{T}(C^\gamma(B_2))}
\lesssim R^{-\alpha}\|u\|_{L^\infty_{T}(L^\infty(\mR^d))}.
\end{align*}
For $I_2$, by Lemma \ref{N7}, we have for all $\eps>0$,
$$
I_2\leq \eps\|\chi_3 u_R\|_{\bB^{\alpha+\gamma}_{TR^{-\alpha}}}+c_\eps\|\chi_3 u_R\|_{\bB^0_{TR^{-\alpha}}}
\lesssim R^{-\alpha}\left(\eps\sup_{t\in[0,T]}\|u(t)\|^{(0)}_{\alpha+\gamma;D}+c_\eps \|u\|_{L^\infty_T(L^\infty(\mR^d))}\right).
$$
For $I_3$ and $I_4$, we have 
\begin{align*}
I_3\leq\|f\|_{L^\infty_T(L^\infty(\mR^d))},\ \ I_4\lesssim R^{-\alpha}\|u\|_{L^\infty_T(L^\infty(\mR^d))}.
\end{align*}
Combining the above calculations, we obtain that for all $\eps>0$,
$$
\|u^{x_0}_R\|_{L^\infty_T(C^{\alpha+\gamma}(B_1))}\leq
\eps\sup_{t\in[0,T]}\|u(t)\|^{(0)}_{\alpha+\gamma;D}+c_\eps \|u\|_{L^\infty_T(L^\infty(\mR^d))}+\lambda^\alpha_D\|f\|_{L^\infty_T(L^\infty(\mR^d))}.
$$
Taking supremum with respect to $x_0\in D$ and by \eqref{HL6} and choosing $\eps$ small enough, we obtain
$$
\sup_{t\in[0,T]}\|u(t)\|^{(0)}_{\alpha+\gamma;D}\leq c \|u\|_{L^\infty_T(L^\infty(\mR^d))}+c\|f\|_{L^\infty_T(L^\infty(\mR^d))},
$$
which gives the desired estimate by $\|u\|_{\alpha+\gamma;D_0}\leq c \|u\|^{(0)}_{\alpha+\gamma;D}$ for some $c>0$.
\end{proof}

The following interior estimate in weighted H\"older spaces is slightly different from Theorem \ref{Main0}.
The key point is that we {\em do not}  assume $u\in \mB^{(-\theta)}_{\alpha+\gamma}(D)$ posteriorily. 
We use a trick from \cite{Ro-Se1}.
\bl\label{Le83}
Let $D$ be a bounded $C^2$-domain and $\alpha\in(0,2)$, $\gamma\in(0,1)$. 
Suppose that $(u,f)\in {\bf H}^{\alpha+\gamma}(D)$ with $f\in \mB^{(\a-\theta)}_{\gamma, T}$ for some $\theta\in (0,\a/2)$,  satisfies
$$
\p_t u=\Delta^{\frac{\alpha}{2}}u+f\mbox{ in $\mR_+\times D$},\ \ u|_{\mR_+\times D^c}=0,\ \ u(0)=0.
$$
Then we have $u\in\mB^{(-\theta)}_{\alpha+\gamma}(D)$ and there is a constant $c=c(\alpha,\theta,\gamma,d, \lambda_D)>0$ such that
\begin{align}\label{EY2}
\|u\|_{\mB^{(-\theta)}_{\alpha+\gamma;T}(D)}\leq c(1+T)\|f\|_{\mB^{(\alpha-\theta)}_{\gamma;T}(D)},\ \ T>0.
\end{align}
\el
\begin{proof}
For $x_0\in D$, let $R:=d_{D^c}(x_0)/8$ and $u_R,f_R, w_R$ be as in \eqref{EY3}. We clearly have
$$
\p_t w_R=\Delta^{\frac{\alpha}{2}} w_R+f_R+\Delta^{\frac{\alpha}{2}}((1-\chi_3) u_R)\ \ in \ \ \mR_+\times B_2.
$$
By \eqref{EY4} and \eqref{EY1}, we have
\begin{align}\label{JH3}
\begin{split}
&R^{-\alpha}\|u^{x_0}_R\|_{L^\infty_T(C^{\alpha+\gamma}(B_1))}
=\|w_R\|_{L^\infty_{TR^{-\alpha}}(C^{\alpha+\gamma}(B_1))}
\lesssim \|f_R\|_{L^\infty_{TR^{-\alpha}}(C^\gamma(B_2))}\\
&\quad+\|\Delta^{\frac{\alpha}{2}}((1-\chi_3) u_R)\|_{L^\infty_{TR^{-\alpha}}(C^\gamma(B_2))}
+\|w_R\|_{L^\infty_{TR^{-\alpha}}(L^\infty(\mR^d))}=I_1+I_2+I_3.
\end{split}
\end{align}
For $I_1$, by \eqref{HL6} we have
\begin{align*}
I_1\leq R^{\theta-\alpha}\|f\|_{\mB^{(\alpha-\theta)}_{\gamma;T}(D)}.
\end{align*}
To estimate $I_2$, we drop the time variable for simplicity. From the calculations in \eqref{PY4}, it is easy to see that
\begin{align*}
\|\Delta^{\frac{\alpha}{2}}((1-\chi_3) u_R)\|_{C^\gamma(B_2)}
=\left\|\Big(\Delta^{\frac{\alpha}{2}}((1-\chi^{x_0}_{3R}) u)\Big)^{x_0}_{R}\right\|_{C^\gamma(B_2)}
\lesssim\int_{B^c_{2R}(x_0)}\frac{|u(z)|}{(|z-x_0|-R)^{d+\alpha}}\dif z.
\end{align*}
Using the change of variable $z-x_0=Ry$, we have
\begin{align*}
\int_{B^c_{2R}(x_0)}\frac{|u(z)|}{(|z-x_0|-R)^{d+\alpha}}\dif z=
R^{-\alpha}\int_{B^c_{2}}\frac{|u(Ry+x_0)|}{(|y|-1)^{d+\alpha}}\dif y.
\end{align*}
Recalling $R=d_{D^c}(x_0)/8$, we have
$$
|u(Ry+x_0)|\leq cR^\theta\|u\|_{\cC^{(-\theta)}_{0}(D)},\ \ |y|\leq 1,
$$
and further,
$$
 |u(Ry+x_0)|\leq cR^\theta(1+|y|^\theta)\|u\|_{\cC^{(-\theta)}_{0}(D)},\ \ y\in\mR^d. 
$$
Thus by Remark \ref{Re68}, we get
$$
I_2\leq cR^{\theta-\alpha}\|u\|_{\mB^{(-\theta)}_{0;T}(D)}\leq cR^{\theta-\alpha}\|f\|_{\mB^{(\alpha-\theta)}_{0;T}(D)},
$$
and also
$$
I_3\leq \|u_R\|_{L^\infty_{TR^{-\alpha}}(L^\infty(B_4))}\leq cR^{\theta-\alpha}\|u\|_{\mB^{(-\theta)}_{0;T}(D)}
\leq cR^{\theta-\alpha}\|f\|_{\mB^{(\alpha-\theta)}_{0;T}(D)}.
$$
Combining the above estimates, we obtain
$$
R^{-\theta}\|u^{x_0}_R(t)\|_{L^\infty_T(C^{\alpha+\gamma}(B_1))}\leq c\|f\|_{\mB^{(\alpha-\theta)}_{\gamma;T}(D)},
$$
which in turn gives the desired estimate by \eqref{HL6}.
\end{proof}

For $x\in\mR^d$, let $\mP_x$ be the law of rotationally invariant and symmetric $\alpha$-stable process $Z$ in canonical space $\Omega$ starting from $x$.
For bounded measurable function $\varphi$, we define
$$
P^D_t\varphi(x):=\mE_x\left( \varphi(X_t); t<\tau_D\right).
$$
It is well known that $P^D_t$ is a strong continuous symmetric Markov semigroup in $L^2(D)$ (see \cite{Ch-Ki-So1}).
\bl
For any $\varphi\in C^2_c(D)$ with compact support in $D$, it holds that
$$
\p_t P^D_t\varphi(x)=P^D_t\Delta^{\frac{\alpha}{2}} \varphi(x),\ \ x\in D.
$$
In particular, for any bounded measurable $f:\mR_+\times D\to\mR$, $u(t,x):=\int^t_0 P^D_{t-s}f(s,x)\dif s$
is a weak solution of \eqref{Lap}. More precisely,
for any $\varphi\in C^2_c(D)$, we have
\begin{align}\label{EY5}
\<u(t),\varphi\>_D=\int^t_0 \<u(s), \Delta^{\frac{\alpha}{2}}\varphi\>_D\dif s+\int^t_0 \<f(s), \varphi\>_D\dif s,
\end{align}
where $\<\cdot,\cdot\>_D$ denotes the inner product in $L^2(D)$.
\el
\begin{proof}
(i) Let supp($\varphi$)$\subset D_0\Subset D$. We have
\begin{align*}
|P^D_s\varphi(x)-\varphi(x)|&
\leq |\mE_x\varphi(X_s)-\varphi(x)|+\mE_x(\varphi(X_s); s\geq \tau_D)\\
&\leq\int^s_0|\mE_x(\Delta^{\frac{\alpha}{2}}\varphi)(X_r)|\dif r+\|\varphi\|_\infty\mP_x(X_s\in D_0; s\geq \tau_D).
\end{align*}
Define $\sigma_{D_0}:=\inf\{t>0: X_t\in D_0\}$,  
$\sigma'_{D_0}:=\inf\{t>\tau_D: X_t\in D_0\}$.
For $x\in D$, we have
\begin{align*}
\mP_x(X_s\in D_0; s\geq \tau_D)\leq \mP_x(s\geq \sigma'_{ D_0}>\tau_D)
&=\mP_x(s\geq \sigma'_{ D_0}=\tau_D+\sigma_{ D_0}\circ\theta_{\tau_D}>\tau_D)\\
&=\mE_x(\mP_{X_{\tau_D}}(s>\sigma_{ D_0}); s>\tau_D).
\end{align*}
Let $\delta:={\rm dist}(D_0, D^c)$. For $y\in D^c$, we clearly have $\tau_{B_\delta(y)}\leq\sigma_{ D_0}$, and so
by Lemma \ref{Le501},
$$
\sup_{y\in D^c}\mP_{y}(\sigma_{ D_0}<s)\leq\sup_{y\in D^c}\mP_{y}\left(\tau_{ B_\delta(y)}<s\right)\leq c_0 s/\delta^\alpha.
$$
Similarly, for $x\in D$, letting $d_x:=d_{D^c}(x)/2$, we also have
$$
\mP_x(\tau_D<s)\leq \mP_x\left(\tau_{B_{d_x}(x)}<s\right)\leq c_0s/d_x^\alpha.
$$
Hence,
$$
\mP_x(X_s\in D_0; s\geq \tau_D)\leq c_0^2s^2/(\delta d_x)^\alpha.
$$
By the dominated convergence theorem, we obtain
\begin{align*}
\lim_{s\downarrow 0}(P^D_{t+s} \varphi(x)-P^D_t\varphi(x))/s=P^D_t\Delta^{\frac{\alpha}{2}} \varphi(x).
\end{align*}
(ii) By definition and Fubini's theorem, we have
\begin{align*}
\int^t_0 \<u(s), \Delta^{\frac{\alpha}{2}}\varphi\>_D\dif s
&=\int^t_0\!\!\! \int^s_0 \<P^D_{s-r}f(r), \Delta^{\frac{\alpha}{2}}\varphi\>_D\dif r\dif s=\int^t_0\!\!\! \int^s_0 \<f(r), P^D_{s-r}\Delta^{\frac{\alpha}{2}}\varphi\>_D\dif r\dif s\\
&=\int^t_0\!\!\! \int^t_r \<f(r), \p_ sP^D_{s-r}\varphi\>_D\dif s\dif r=\int^t_0 \<f(r), P^D_{t-r}\varphi-\varphi\>_D\dif r\\
&=\<u(t),\varphi\>_D-\int^t_0 \<f(r), \varphi\>_D\dif r.
\end{align*}
We complete the proof.
\end{proof}

Now we are ready to give the proof of Theorem \ref{Th52}.
\begin{proof}[Proof of Theorem \ref{Th52}]
For $\eps>0$, define $D_\eps:=\{x\in D: d_{D^c}(x)>\eps\}$. Let $\rho$ be a nonnegative smooth function with support in $B_1$ and
$\int_{\mR^d}\rho=1$. Define
$$
\rho_\eps(x):=\eps^{-d}\rho(x/\eps),\ \ \chi_\eps:={\bf 1}_{D_{2\eps}}*\rho_\eps,\ f_\eps:=f\chi_\eps.
$$
By definition, one sees that 
$$
f_\eps(t,x)=0,\ \ x\notin D_\eps,\ \ f_\eps(t,x)=f(t,x),\ \ x\in D_{3\eps},
$$
and
\begin{align}\label{EY6}
\|f_\eps\|_{L^\infty_T(C^\gamma(D))}\leq C_\eps <\infty,\quad
\sup_{\eps\in(0,1)}\|f_\eps\|_{\mB^{(\alpha-\theta)}_{\gamma;T}(D)}\lesssim\|f\|_{\mB^{(\alpha-\theta)}_{\gamma;T}(D)}.
\end{align}
Let $u_\eps(t,x):=\int^t_0P^D_{t-s}f_\eps(s,x)\dif s$ be the weak solution of \eqref{EY5}.
By a standard smoothing technique and interior estimate \eqref{EY1},  one can show that
$(u_\eps, f_\eps)\in \bH^{\alpha+\gamma}(D)$ (see \eqref{AL5}) solve the following Dirichlet problem:
\begin{align}\label{EY9}
u_\eps(t,x)=\int^t_0\Delta^{\frac{\alpha}{2}}u_\eps(s,x)\dif s+\int^t_0f_\eps(s,x)\dif s,\mbox{\ \ $(t,x)\in\mR_+\times D$}.
\end{align}
By Lemma \ref{Le83} and \eqref{EY6}, there is a constant $c>0$ such that
for all $T>0$ and $\eps\in(0,1)$,
$$
\|u_\eps\|_{\mB^{(-\theta)}_{\alpha+\gamma;T}(D)}\leq c\|f_\eps\|_{\mB^{(\alpha-\theta)}_{\gamma;T}(D)}
\leq c(1+T)\|f\|_{\mB^{(\alpha-\theta)}_{\gamma;T}(D)}.
$$
Moreover, by \eqref{HL5} we also have
$$
\|\p_t u_\eps\|_{\mB^{(\alpha-\theta)}_{\gamma;T}(D)}\leq \|\Delta^{\frac{\alpha}{2}}u_\eps\|_{\mB^{(\alpha-\theta)}_{\gamma;T}(D))}
+\|f_\eps\|_{\mB^{(\alpha-\theta)}_{\gamma;T}(D)}
\leq c(1+T)\|f\|_{\mB^{(\alpha-\theta)}_{\gamma;T}(D)}.
$$
Hence, by Lemma \ref{Le42} there exist a subsequence $\eps_k\to 0$ and a $u\in \mB^{(-\theta)}_{\alpha+\gamma;T}(D)$ such that
for any $\delta\in(0,\alpha\wedge 1)$ and $D_0\Subset D$,
$$
\lim_{k\to\infty}\|u_{\eps_k}-u\|_{L^\infty_T(C^{\alpha+\gamma-\delta}(\bar D_0))}=0.
$$
By taking limits $k\to\infty$ for \eqref{EY9} along $\eps_k$, we find that $u$ satisfies \eqref{Lap0}.
\end{proof}

\subsection{Monte-Carlo simulation of Example 1.3}

In this subsection we use R-language to simulate the Example 1.3.
In the following simulations, we choose $D=(0,1)$, $\alpha=0.5$,  and the points along $t$-axis and $x$-axis equal to $50$.

\begin{figure*}[h]\small
\centering
\includegraphics
[height=0.25\textwidth,width=0.45\textwidth]{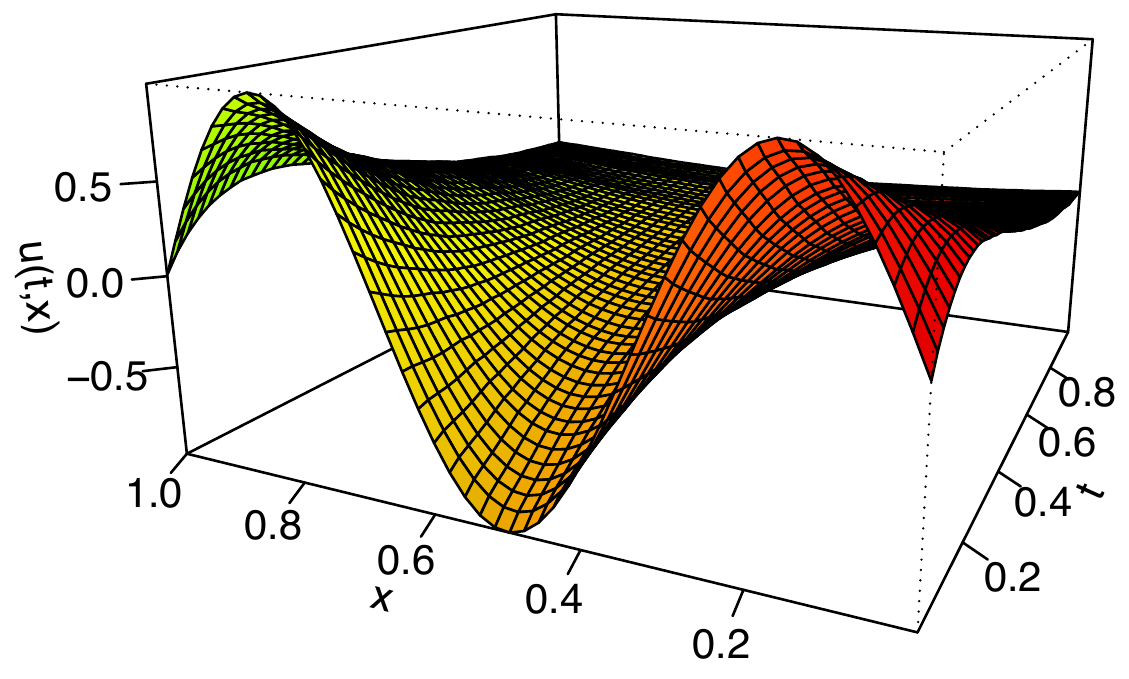}
\includegraphics
[height=0.25\textwidth,width=0.45\textwidth]{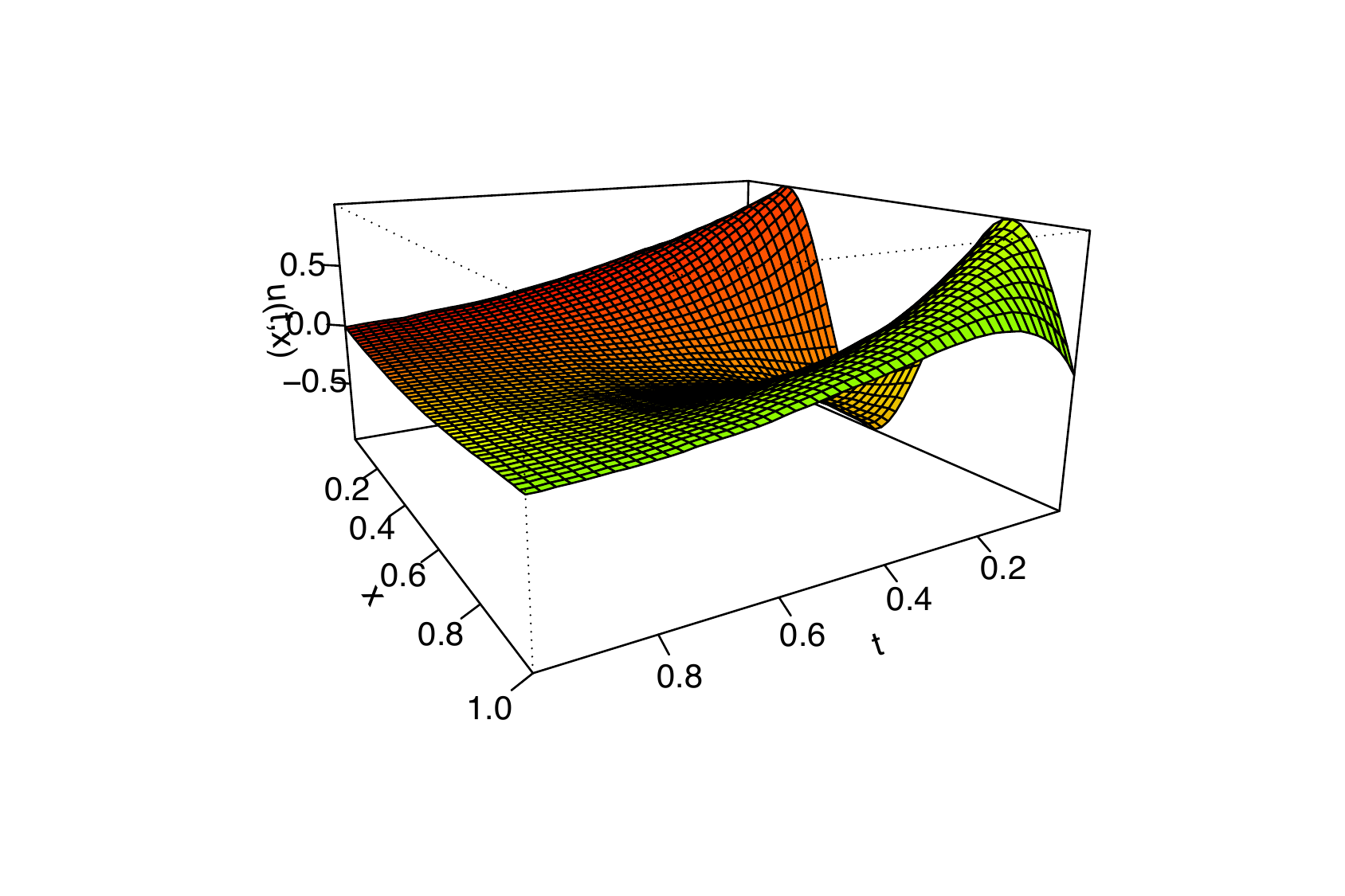}

\caption{$b(x)=\frac{1}{2}-x,\ \varphi(x)=\sin (3\pi x)$}\label{fige3c}

\end{figure*}
\begin{figure*}[h]\small
\centering
\includegraphics
[height=0.25\textwidth,width=0.45\textwidth]{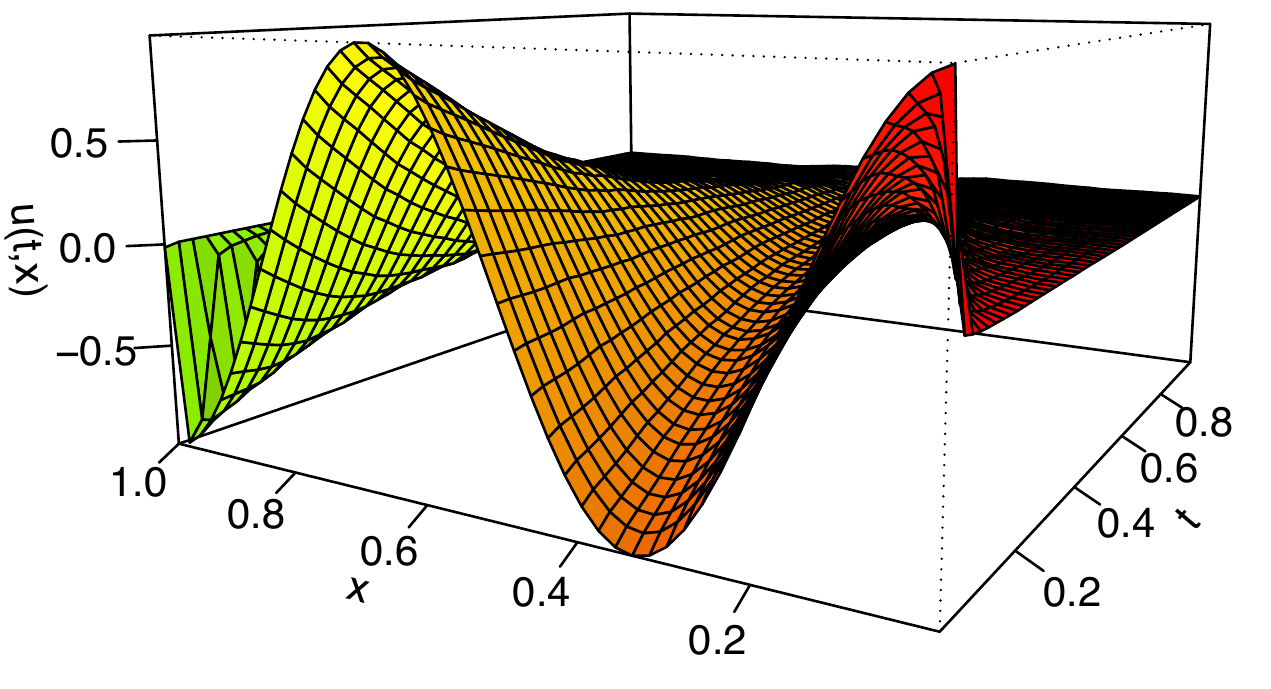}
\includegraphics
[height=0.25\textwidth,width=0.45\textwidth]{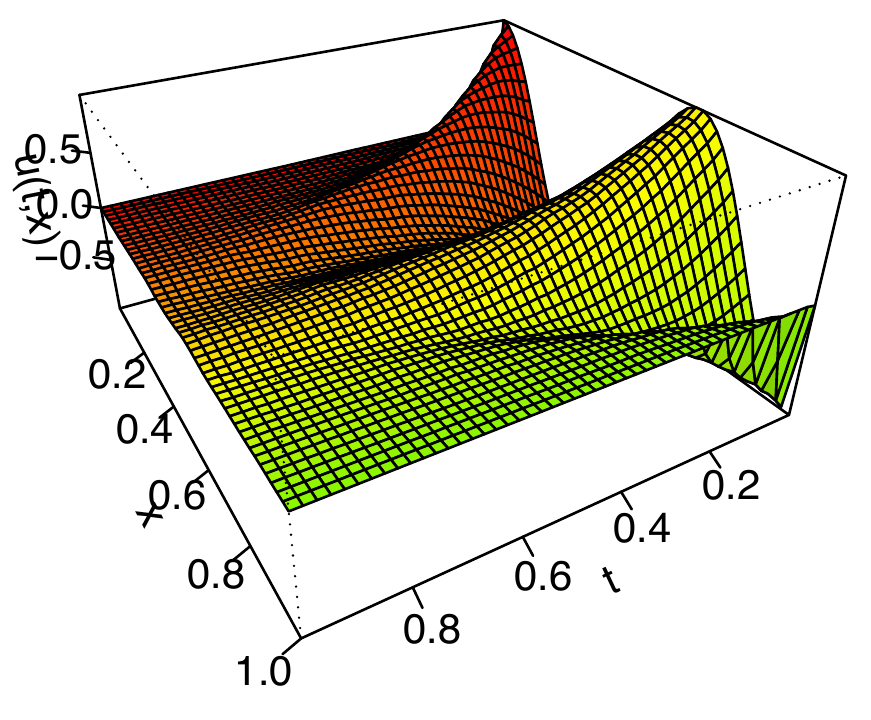}

\caption{$b(x)=x-\frac{1}{2}, \ \varphi(x)=\sin \Big(3\pi x+\pi/2\Big)$}\label{fige4c}

\end{figure*}
\begin{figure*}[h]\small
\centering
\includegraphics
[height=0.25\textwidth,width=0.45\textwidth]{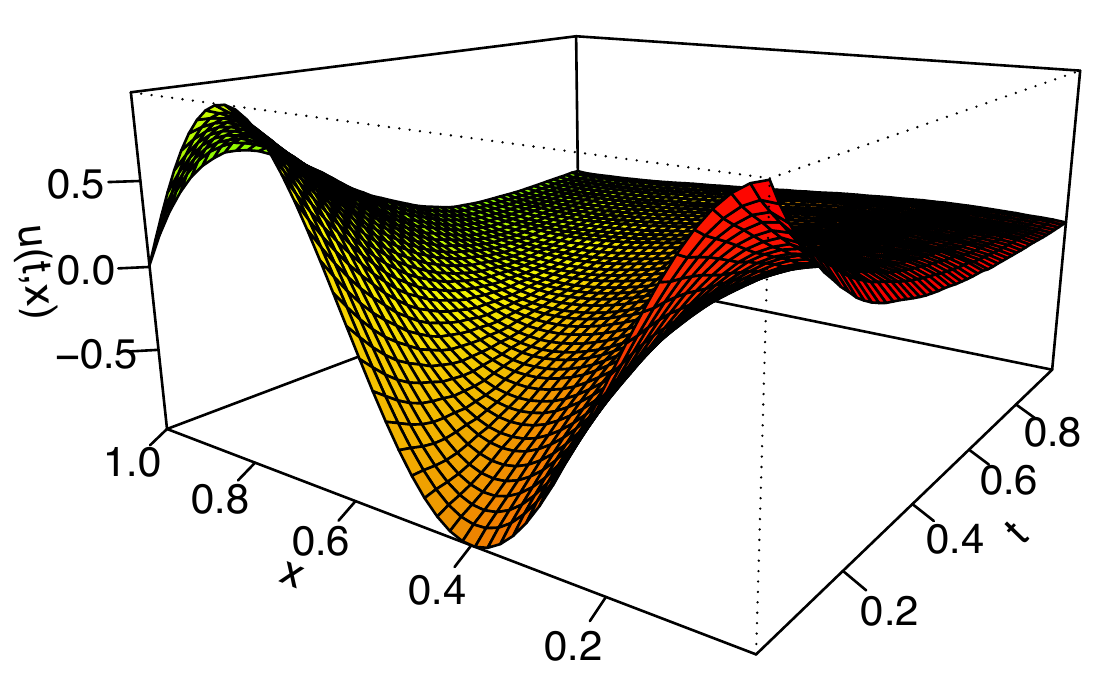}
\includegraphics
[height=0.25\textwidth,width=0.45\textwidth]{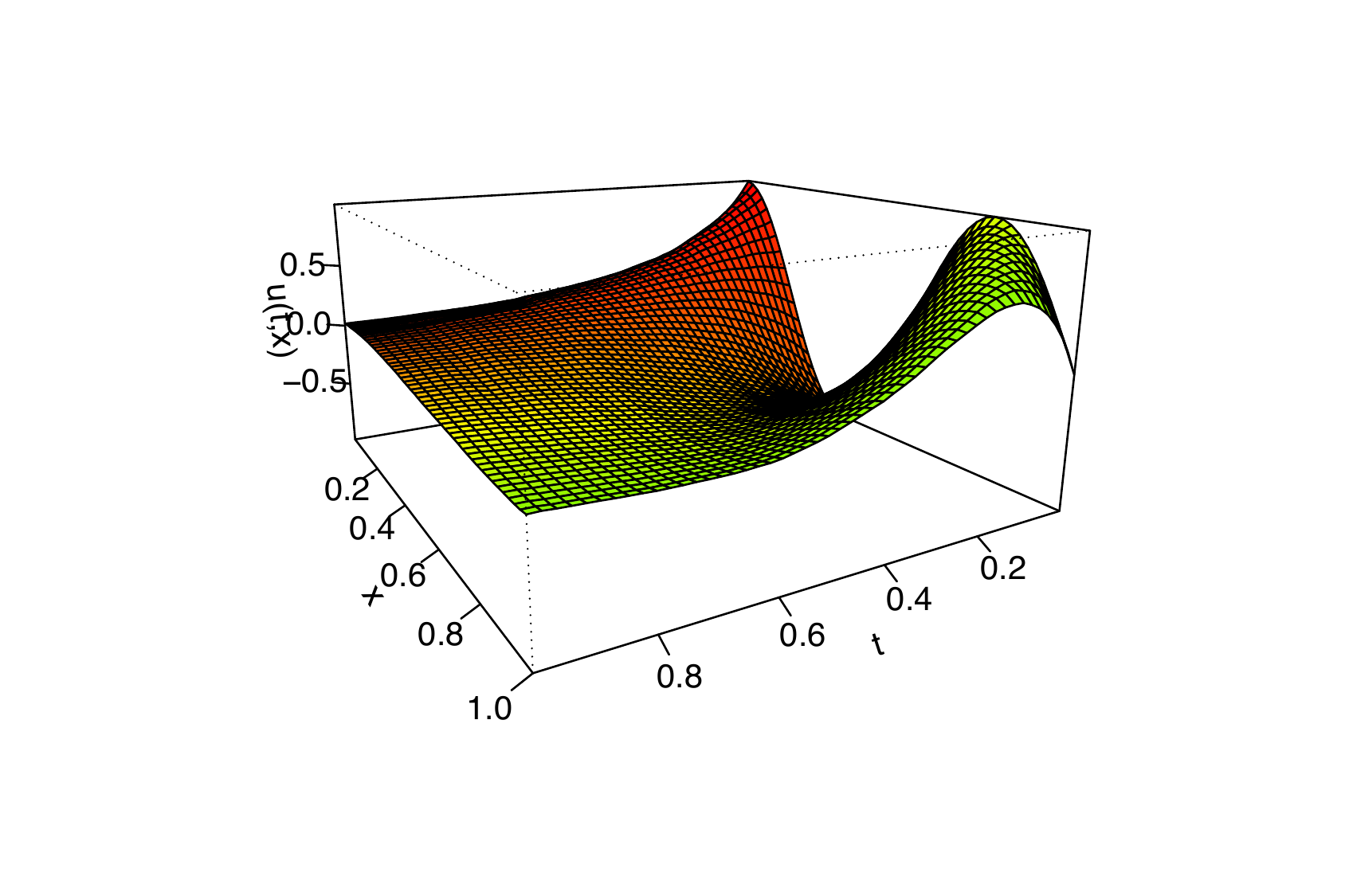}

\caption{$b(x)=-x, \ \varphi(x)=\sin \Big(5\pi (x+1)/2\Big)$}\label{fige5c}
\end{figure*}

In Figure \ref{fige3c}, since $b(z_0)\cdot\vec{n}(z_0)<0$, one sees that for any $t>0$, 
the solution $u(t,x)$ is not continuous up to the boundary.
However, in Figure \ref{fige4c}, since $b(z_0)\cdot\vec{n}(z_0)>0$, one sees that for any $t>0$, the solution $u(t,x)$ is  continuous up to the boundary 
($\lim_{D\ni x\to 0} u(t,x)=\lim_{D\ni x\to 1} u(t,x)=0$), even if the initial value is non-zero at the boundary. 
In Figure \ref{fige5c}, since $b(0)\cdot\vec{n}(0)=0$ and $b(1)\cdot\vec{n}(1)<0$, one sees that 
at the boundary point $1$, $u(t,x)$ behaves like Figure \ref{fige3c},
and at the boundary point $0$, $u(t,x)$ behaves like Figure \ref{fige4c}.
It should be noticed that in all the above figures, when $t$ becomes larger and larger, $u(t,x)$ will be close to zero due to the dissipativity of $\Delta^{\alpha/2}$.

\vspace{5mm}

{\bf Acknowledgement:}
The authors would like to thank Professors Zhen-Qing Chen, Renming Song for their quite useful conversations.

\end{document}